\documentclass[11pt,letterpaper]{amsart}
\usepackage[margin=1in]{geometry}
\usepackage{amsmath,amssymb,amsthm,mathtools,mathrsfs}
\usepackage{microtype,booktabs,tabularx,array,enumitem}
\usepackage[numbers,sort&compress]{natbib}
\usepackage{xcolor,graphicx}
\usepackage[colorlinks=true,linkcolor=blue!45!black,citecolor=blue!45!black,urlcolor=blue!45!black]{hyperref}
\usepackage{doi}
\hypersetup{pdftitle={A closed speed-position SDE for KdV solitons under multiplicative noise},pdfauthor={Justin Holmer and Katherine Zhiyuan Zhang}}
\hypersetup{bookmarksdepth=2}
\newcommand{\R}{\mathbb R}
\newcommand{\E}{\mathbb E}
\newcommand{\Prob}{\mathbb P}
\newcommand{\eps}{\varepsilon}
\newcommand{\dd}{\,\mathrm d}
\newcommand{\Ito}{It\^o}
\newcommand{\sech}{\operatorname{sech}}
\newcommand{\HS}{\mathrm{HS}}
\newcommand{\Xb}{\mathcal X_\beta}
\newcommand{\Kb}{\mathcal K}
\newcommand{\Cs}{\mathcal C}
\newcommand{\Hs}{\mathscr H}
\newcommand{\Ac}{\mathcal A}
\newcommand{\bc}{\bar C}
\newcommand{\ba}{\bar a}
\newcommand{\ip}[2]{\langle #1,#2\rangle}
\newcommand{\norm}[1]{\lVert #1\rVert}
\newcommand{\Op}{O_{\Prob}}

\newtheorem{theorem}{Theorem}[section]
\newtheorem{proposition}[theorem]{Proposition}
\newtheorem{lemma}[theorem]{Lemma}
\newtheorem{corollary}[theorem]{Corollary}
\theoremstyle{definition}

\newtheorem{assumption}[theorem]{Assumption}

\theoremstyle{remark}
\newtheorem{remark}[theorem]{Remark}
\numberwithin{equation}{section}
\title[A closed speed--position SDE for stochastic KdV]{A closed speed--position SDE for KdV solitons under multiplicative noise}
\author{Justin Holmer}
\address{Brown University}
\author{Katherine Zhiyuan Zhang}
\address{Northeastern University}
\date{November 25, 2026}
\begin{document}
\begin{abstract}
We study a Korteweg--de Vries soliton under small, spatially homogeneous multiplicative \Ito\ noise. Building on the first-order speed and position model of de Bouard and Debussche and the closed amplitude equation of Westdorp and Hupkes, we prove a joint approximation for speed and position. We use a one-sided phase and measure speed by $C=-4\lambda$, where $\lambda$ is the isolated Schr\"odinger eigenvalue in the KdV Lax pair. Isospectrality makes the deterministic KdV drift of $C$ vanish. For $T\ge1$ and $S=\eps^2T$ small, a closed two-dimensional SDE approximates the spectral speed and $\eps^2$ times the position with stopped maximal $L^2$ errors $O(\eps\sqrt S)$ and $O(\eps S^{3/2})$. The fixed tube survives with probability tending to one for every $1\le T=o(\eps^{-2})$. At $T=\eps^{-2}/\log^2(1/\eps)$ these errors resolve the positive spectral \Ito\ drift and its accumulated contribution to position. We obtain second-order fluctuations, stopped and survival-conditioned means, and an $O(\eps)$ bounded-Lipschitz comparison for the normalized paths. The estimates transfer to geometric soliton coordinates.
\end{abstract}
\maketitle

\tableofcontents
\section{Introduction}\label{sec:intro}

\subsection{From modulation equations to an independently solvable model}
We consider the \Ito\ equation
\begin{equation}\label{eq:skdv}
 \dd u+\left(u_{xxx}+\frac12\partial_x(u^2)\right)\dd t
       =\eps u\phi\dd\mathcal W_t,
 \qquad u(0,x)=Q_{c_0}(x-a_0),
\end{equation}
on the real line. Here $\mathcal W$ is cylindrical Brownian motion on $L^2(\R)$, $\phi h=k*h$, and $0<\eps\ll1$. Our normalization is
\begin{equation}\label{eq:profile}
 Q_c(y)=3c\sech^2\!\left(\frac{\sqrt c}{2}y\right),
 \qquad Q_c''-cQ_c+\frac12Q_c^2=0.
\end{equation}
The deterministic soliton is $Q_{c_0}(x-a_0-c_0t)$.

Stochastic modulation supplies a decomposition
\begin{equation}\label{eq:decomp-intro}
 u(t,x+a(t))=Q_{c(t)}(x)+v(t,x)
\end{equation}
and exact equations for the parameters whose coefficients depend on the radiation $v$.
For the same multiplicative forcing, de Bouard and Debussche derived a frozen first-order system for speed and position, including their correlated Brownian fluctuations \cite[Sections~4--5]{deBouardDebussche2009}. Westdorp and Hupkes constructed amplitude and position approximations \cite{WestdorpHupkes2024} and proved a quantitative comparison with a closed amplitude SDE \cite[equation~(1.7) and Theorem~1.2]{WestdorpHupkes2026}. We seek a joint comparison accurate enough to resolve the accumulated \Ito\ drift of the speed in both observables.

For the geometric decomposition we use the one-sided orthogonality conditions of \cite[equation~(1.11)]{WestdorpHupkes2026},
\begin{equation}\label{eq:intro-orth}
 \ip{v}{Q_c}=0,\qquad
 \ip{v}{Z_c}=0,\qquad
 Z_c(y)=\int_{-\infty}^y\partial_cQ_c(z)\dd z.
\end{equation}
The primitive $Z_c$ tends to a nonzero constant as $y\to+\infty$. Its pairing with $v$ is continuous in the one-sided weighted space of Section~\ref{sec:geometry}. The adjoint identities for these tests cancel linear radiation terms in the deterministic modulation equations. The remaining deterministic terms are quadratic in localized radiation.

For the speed observable we use the Schr\"odinger operator from KdV inverse scattering and the Lax pair \cite{GardnerGreeneKruskalMiura1967,Lax1968},
\begin{equation}\label{eq:intro-spectral}
 \Hs(u)=-\partial_x^2-\frac16u,
 \qquad C=\Cs(u):=-4\lambda(u).
\end{equation}
Here $\lambda(u)$ is the simple isolated eigenvalue continuing $-c/4$ at $Q_c$. The geometric parameter $c$ defines the decomposition in \eqref{eq:decomp-intro}. The spectral parameter $C$ is then a functional of the solution. Deterministic KdV conserves $C$ exactly, even in the presence of radiation. Its stochastic coefficients still depend on $v$.

The effective equation has the form
\begin{equation}\label{eq:intro-effective}
 \begin{aligned}
 \dd\ba&=\bc\dd t+\eps\sigma_a(\bc)\dd B_t
                  +\eps^2\gamma_a(\bc)\dd t,\\
 \dd\bc&=\eps\sigma_C(\bc)\dd B_t
                  +\eps^2 b_{\rm sp}(\bc)\dd t.
 \end{aligned}
\end{equation}
Here $B$ is two-dimensional Brownian motion. The two diffusion rows $\sigma_a,\sigma_C$ give the position and speed noises their joint covariance. Section~\ref{sec:coefficients} gives the covariance and drift coefficients. Once these deterministic functions have been computed, \eqref{eq:intro-effective} can be simulated independently of \eqref{eq:skdv}. The pathwise estimate uses an adapted coupling. Path statistics of the reduced model may be computed with a fresh Brownian motion.

\subsection{Main theorem and the time window}
Write
\begin{equation}\label{eq:scales}
 S=\eps^2T,\qquad \Lambda_\eps=\log(1/\eps),\qquad
 T_\eps=\frac{\eps^{-2}}{\Lambda_\eps^2},\qquad 0<\eps<e^{-1}.
\end{equation}
We give estimates in terms of $S$ and use the particular choice $T_\eps$ for second-order asymptotics. The distinction between persistence for $S\to0$ and stopped approximation at fixed $S>0$ is made precise below.

Assume that $k$ is real, nonzero, belongs to $L^1(\R)$, and satisfies
\[
 \sum_{j=0}^2\int_\R e^{2\mu|x|}|k^{(j)}(x)|^2\dd x<\infty
\]
for some $\mu>0$, as in Assumption~\ref{ass:noise}. Fix a compact positive speed interval $I=[c_-,c_+]$ with $c_0\in I^\circ$, and choose
$0<\beta<\min(\mu/2,\sqrt{c_-}/4)$. Set
\[
 \norm{f}_{H^j_\beta}=\norm{e^{\beta y}f}_{H^j},\qquad
 r_\beta(v)=\norm{v}_{H^1_\beta}.
\]
The coefficients of \eqref{eq:intro-effective} are computed from $Q_c$ and $k$ in Section~\ref{sec:coefficients}, using the covariance factor \eqref{eq:sigma-factor}. Choose bounded smooth extensions of these coefficients with bounded derivatives through order three and with $\sigma_C=(g,0)$ globally, where $g=\sqrt{q_{CC}}$ on $I$. Together with the linear transport drift $(\bc,0)$ they give globally Lipschitz coefficients, hence a global strong solution. The stopped dynamics are independent of the extension.

For a small fixed $I_0=[c_0-\rho,c_0+\rho]\Subset I$ and a small fixed radius $\delta$, let $\tau$ be the first exit from
\begin{equation}\label{eq:stopping-regime}
 c,C,\bc\in I_0,\qquad \norm{v}_{H^1}+r_\beta(v)<\delta,
\end{equation}
initially intersected with the weighted-chart lifetime. More precisely, speed exit means reaching the boundary of $I_0$. The construction and continuation assertion below make this a well-defined common stop.

\begin{theorem}[Joint speed and position comparison]\label{target:main}
For the solution of \eqref{eq:skdv} with exact soliton initial data, there exist fixed $\rho,\delta,S_0,\eps_0>0$ with the following properties for $0<\eps<\eps_0$.

The decomposition $T_au=Q_c+v$, selected by \eqref{eq:intro-orth}, continues up to the explicit exit \eqref{eq:stopping-regime}. The spectral speed is $C=-4\lambda(u)$ as in \eqref{eq:intro-spectral}. There is an adapted two-dimensional Brownian motion $B$ such that the solution of \eqref{eq:intro-effective}, with $(\ba(0),\bc(0))=(a_0,c_0)$, obeys the following estimates. For $1\le T\le S_0\eps^{-2}$, $S=\eps^2T$, and every $p\ge1$,
\begin{equation}\label{eq:target-radiation}
 \left\|\sup_{t\le T\wedge\tau}\norm{v(t)}_{H^1}^2
       +\int_0^{T\wedge\tau}r_\beta(v(t))^2\dd t\right\|_{L^p(\Omega)}
       \le C_pS,
\end{equation}
\begin{equation}\label{eq:main-weighted-max}
 \left\|\sup_{t\le T\wedge\tau}r_\beta(v(t))^2\right\|_{L^p(\Omega)}
       \le C\eps^2[\log(2+T)+p],
\end{equation}
and
\begin{equation}\label{eq:target-errors}
 \begin{aligned}
 \left\|\sup_{t\le T\wedge\tau}|C(t)-\bc(t)|\right\|_{L^2}
       &\le C\eps\sqrt S,\\
 \left\|\sup_{t\le T\wedge\tau}|\eps^2(a(t)-\ba(t))|\right\|_{L^2}
       &\le C\big(\eps S^{3/2}+\eps^2S+\eps^3\sqrt S\big)
        \le C'\eps S^{3/2}.
 \end{aligned}
\end{equation}
The exit probability satisfies
\begin{equation}\label{eq:main-exit}
 \Prob(\tau\le T)\le K e^{-\kappa/S}+K(1+T)e^{-\kappa_*/\eps^2}.
\end{equation}
The constants $K,\kappa,\kappa_*>0$ and the constants in the comparison are independent of $\eps,T$. Thus $\Prob(\tau>T)\to1$ for every deterministic $1\le T=o(\eps^{-2})$. At $T=T_\eps$,
\begin{equation}\label{eq:target-max-exit}
 \left\|\sup_{t\le T_\eps\wedge\tau}r_\beta(v(t))^2\right\|_{L^2}
       \le C\eps^2\Lambda_\eps,
 \qquad \Prob(\tau>T_\eps)\longrightarrow1,
\end{equation}
and the two errors in \eqref{eq:target-errors} are bounded by $C\eps/\Lambda_\eps$ and $C\eps/\Lambda_\eps^3$, respectively. The symplectic geometric speed also satisfies
\begin{equation}\label{eq:target-transfer}
 \left\|\sup_{t\le T_\eps\wedge\tau}|c(t)-\bc(t)|\right\|_{L^2}
       \le C\left(\frac{\eps}{\Lambda_\eps}+\eps^2\Lambda_\eps\right)
       \lesssim\frac{\eps}{\Lambda_\eps}.
\end{equation}
\end{theorem}
Section~\ref{sec:main-proof} proves the theorem. A direct bound for the time integral of weighted radiation in \eqref{eq:target-radiation} avoids the logarithmic loss that would result from using only the weighted maximum.

For $T=s\eps^{-2}$ with $0<s\le S_0$ fixed, the stopped errors in \eqref{eq:target-errors} are still $O_s(\eps)$. The effective speed on this slow time scale is a nondegenerate diffusion up to the boundary of $I_0$. It has a positive probability of reaching the boundary by $s$, independent of $\eps$. Proposition~\ref{prop:fixed-slow-exit} establishes a positive lower bound for the PDE-only exit as well. For the normalized paths at $T_\eps$, Corollary~\ref{cor:observable-laws} gives an $O(\eps)$ bound for the distance between laws tested against bounded Lipschitz functionals.

\subsection{Resolution of drift and dependence on coordinates}
The virial chart imposes $\ip{v}{Q_c}=\ip{v}{yQ_c}=0$. Its deterministic phase equation contains the linear term $-\ip{v}{Q_c^2}/\norm{Q_c}_2^2$. Proposition~\ref{prop:virial-chart} compares this chart directly with \eqref{eq:intro-orth}. The change in phase is linear to leading order in the localized radiation, and the change in geometric speed is quadratic. The comparison estimates at $T_\eps$ transfer to both virial observables.

At $T=T_\eps$, the model errors are smaller than the respective spectral-drift scales $\Lambda_\eps^{-2}$ and $\Lambda_\eps^{-4}$ in $C$ and $\eps^2a$. More generally, when $T\to\infty$ and $S=\eps^2T\to0$, the ratios of the respective error bounds to $S$ and $S^2$ are $O(T^{-1/2})$. Section~\ref{sec:second-order} gives the stochastic expansions at $T_\eps$. If $g(c)=\sqrt{q_{CC}(c)}$ and $b_0=b_{\rm sp}(c_0)>0$, the second-order speed correction is
\[
 b_0s+\frac12g(c_0)g'(c_0)(W_s^2-s).
\]
The centered term has strictly positive coefficient by Lemma~\ref{lem:increasing-speed-variance}. Proposition~\ref{prop:pure-drift-comparison} also computes the strictly positive difference between the pure spectral and geometric speed drifts, although the two speed coordinates have the same diffusion row on the soliton. Omitting $b_{\rm sp}$ from the closed SDE produces a discrepancy larger than the proved PDE--SDE error; Corollary~\ref{cor:drift-deletion} identifies that discrepancy.

For moments, Section~\ref{sec:means} proves the expansions
\begin{equation}\label{eq:intro-means}
 \begin{aligned}
 \E[C(T_\eps)-c_0\mid\tau_{\rm P}>T_\eps]
       &=b_0\Lambda_\eps^{-2}+O(\Lambda_\eps^{-4}),\\
 \eps^2\E[a(T_\eps)-a_0-c_0T_\eps\mid\tau_{\rm P}>T_\eps]
       &=\tfrac12b_0\Lambda_\eps^{-4}+O(\Lambda_\eps^{-6}).
 \end{aligned}
\end{equation}
Here $\tau_{\rm P}$ is the PDE-only tube exit, obtained by omitting $\bc$ from \eqref{eq:stopping-regime}. Proposition~\ref{prop:means} also gives the means of the explicitly stopped observables, with the ballistic clock stopped together with the phase.

At $T_\eps$ the resolved position bias comes from integrating the speed drift. The direct order-$\eps^2$ phase drift is smaller than the position error bound. In physical distance the endpoint bound on $a-\ba$ is $O(\eps^{-1}\Lambda_\eps^{-3})$, whereas the leading displacement fluctuation has scale $\eps^{-2}\Lambda_\eps^{-3}$. Absolute accuracy of the center and uniform $H^1$ approximation by the soliton at the coupled position follow on the shorter intervals $T=o(\eps^{-4/3})$ (Corollary~\ref{cor:physical-reconstruction}).
\subsection{Stochastic KdV and scattering-based perturbation theory}
For the real-line KdV equation with multiplicative homogeneous noise, de Bouard and Debussche \cite{deBouardDebussche2009} proved random modulation and exponential exit estimates for an unweighted $H^1$ tube under the hypothesis $k\in L^1\cap H^1$. They obtained a fixed-time-interval Gaussian limit for the rescaled radiation and studied its long-time behavior in an exponentially weighted space after a correction of the phase. Their first-order parameter model already contains the two correlated Brownian motions that drive speed and position. The covariance following \cite[(5.2)]{deBouardDebussche2009} includes the cubic-in-time position variance generated by integrating the speed fluctuation. Our small-noise limits in Section~\ref{sec:memory} recover this mechanism on a growing time interval. For additive forcing, de Bouard and Debussche \cite{deBouardDebussche2007} established modulation and control of the radiation, and de Bouard and Gautier \cite{deBouardGautier2010} quantified the gain in persistence obtained by following the modulated soliton.

Westdorp and Hupkes \cite{WestdorpHupkes2024} constructed tractable approximations for the amplitude and position under translation-invariant multiplicative noise, identified leading drift contributions, and compared their predictions with simulations. They also studied amplitude change under slowly decaying deterministic multiplicative forcing \cite{WestdorpHupkes2025}. Their stochastic analysis \cite{WestdorpHupkes2026} permits an additional deterministic forcing that can change the soliton amplitude by order one. Their weighted persistence estimate in Theorem~1.1 covers $T_\eps$ in the noise-only case. Equation~(1.7) of that paper is a closed scalar amplitude SDE, and Theorem~1.2 estimates the probability that the amplitude approximation exits a prescribed tolerance. Their proof in Section~9 uses the spatially translated noise for the coupling. Our geometric chart \eqref{eq:intro-orth} is their chart (1.11), and their Lemma~3.3 proves that both deterministic modulation errors are quadratic in the localized radiation. In their normalization $\varphi_c=Q_c/2$, so their speed diffusion test $(2/9)c^{-1/2}\varphi_c^2$ coincides with $Q_c^2/(18\sqrt c)$ in \eqref{eq:K0}. We use this phase chart together with the discrete eigenvalue $C$. The resulting joint speed--position approximation has the two error bounds in \eqref{eq:target-errors}, both below the accumulated spectral-drift scales at $T_\eps$.

Gardner, Greene, Kruskal, and Miura \cite{GardnerGreeneKruskalMiura1967} constructed the KdV inverse scattering transform using the Schr\"odinger operator, and Lax \cite{Lax1968} expressed its isospectral evolution through a commutator. Karpman and Maslov \cite{KarpmanMaslov1977} derived soliton parameter equations for perturbed KdV from scattering data. Kaup \cite{Kaup1976} developed the corresponding perturbation expansion for the Zakharov--Shabat scattering transform. For random perturbations, Garnier \cite{Garnier2001} used inverse scattering and diffusion approximation to derive effective parameter laws and describe radiation. Section~6 of that paper includes a closed spectral amplitude law for a spatially uniform random multiplier. In \eqref{eq:skdv} the noise is white in time, its spatial covariance is $k*\widetilde{k}$ with $\widetilde{k}(x)=k(-x)$, and the multiplier is interpreted in the \Ito\ sense.

Wadati \cite{Wadati1983} obtained an explicit random soliton displacement for a spatially uniform additive force. Wadati and Akutsu \cite{WadatiAkutsu1984} extended the analysis to include damping. Herman \cite{Herman1990} studied a damped KdV soliton under space- and time-dependent noise by singular perturbation theory. Bass, Kivshar, Konotop, and Sinitsyn \cite{BassKivsharKonotopSinitsyn1988} surveyed stochastic perturbations of soliton parameters and the associated radiation. Scalerandi, Romano, and Condat \cite{ScalerandiRomanoCondat1998} investigated numerically the broadening of the averaged soliton under additive forcing. Debussche and Printems \cite{DebusschePrintems1999} simulated the response of KdV solitons to homogeneous noise. Cartwright and Gottwald \cite{CartwrightGottwald2021} derived independently solvable collective-coordinate SDEs for amplitude, width, and position for several scalar-noise models, and tested them against full SPDE simulations. In the periodic damped-driven problem, Kuksin and Piatnitski \cite{KuksinPiatnitski2008} proved averaging for the KdV integrals of motion, and Kuksin \cite{Kuksin2010} constructed an effective equation for their joint law on the diffusive time scale.

\subsection{Deterministic modulation and stochastic phase reduction}
Benjamin \cite{Benjamin1972} and Bona \cite{Bona1975} established variational stability for the KdV soliton. Our orbital estimate uses the constrained coercivity underlying their work and the general results of Weinstein \cite{Weinstein1986} and of Bona, Souganidis, and Strauss \cite{BonaSouganidisStrauss1987}. Pego and Weinstein \cite{PegoWeinstein1994} established asymptotic stability using the linearized semigroup in exponentially weighted spaces. Mizumachi and Tzvetkov \cite{MizumachiTzvetkov2014} used that method to give another proof of the $L^2$ stability result obtained by Merle and Vega via the Miura transform \cite{MerleVega2003}. The weighted semigroup and its adjoint modes enter our phase chart and radiation estimate. The local virial and monotonicity arguments of Martel and Merle \cite{MartelMerle2001,MartelMerle2005} and Martel's analysis of the associated linear problems \cite{Martel2006} provide a complementary approach to generalized KdV. Dejak and Sigal \cite{DejakSigal2006} addressed the nondecaying symplectic primitive by regularization in the variable-bottom problem. Holmer \cite{Holmer2011} combined a local virial estimate with a Lyapunov functional to control the KdV soliton parameters in a slowly varying potential.

Related effective dynamics for NLS include the Hamiltonian motion derived by Fr\"ohlich, Gustafson, Jonsson, and Sigal \cite{FrohlichGustafsonJonssonSigal2004} for a slowly varying external potential. Abou Salem and Sulem \cite{AbouSalemSulem2009} derived effective particle dynamics and diffusion limits in a slowly varying random potential. De Bouard and Fukuizumi \cite{deBouardFukuizumi2009} studied a quadratic trapping potential with time-white multiplicative fluctuations and proved modulation and a Gaussian radiation limit. Gnann, Westdorp, and van Winden \cite{GnannWestdorpVanWinden2025} proved exponentially long orbital stability for a parametrically forced, damped NLS equation and derived fixed-interval second-order expansions of the phase shift and fluctuations.

Brassesco, De Masi, and Presutti \cite{BrassescoDeMasiPresutti1995} obtained Brownian fluctuations of a Ginzburg--Landau interface. Giacomin, Poquet, and Shapira \cite{GiacominPoquetShapira2018} identified phase diffusion near stable limit cycles. MacLaurin \cite{MacLaurin2023} developed stochastic phase reduction under continuous spatial isometries, Adams \cite{Adams2024} proved regularity and an \Ito\ formula for the isochronal phase, and Adams and MacLaurin \cite{AdamsMacLaurin2025} used that phase to derive effective dynamics on a stable invariant manifold. In reaction--diffusion systems, Hamster and Hupkes \cite{HamsterHupkes2020} proved exponentially long traveling-wave stability, and van den Bosch, Hamster, and Hupkes \cite{vanDenBoschHamsterHupkes2026} computed higher-order conditional speed and shape corrections. For KdV the speed and translation tangents form a Jordan chain, and the unweighted continuous spectrum meets the imaginary axis. The one-sided exponential weight yields decay for the projected linearized evolution.

\subsection{Stochastic background and organization}
Section~\ref{sec:stochastic} gives the stochastic integral and stopping conventions. We use Da Prato and Zabczyk \cite{DaPratoZabczyk2014} and Pr\'ev\^ot and R\"ockner \cite{PrevotRockner2007} for stochastic evolution equations, Protter \cite{Protter2005} for stochastic integration, and Kunita \cite{Kunita1990} and Krylov \cite{Krylov2011} for changes of stochastic coordinates. Sections~\ref{sec:geometry}--\ref{sec:localized-closure} construct the coordinates and compute their exact equations and localized coefficients. Sections~\ref{sec:weighted-energy} and~\ref{sec:orbital-exit} bound the weighted radiation and the orbital energy. Section~\ref{sec:target} proves the strong comparison. Section~\ref{sec:interpretation} derives its consequences for changes of chart, fluctuations, and path statistics, and Section~\ref{sec:numerics} compares reduced-model predictions with simulations.

\subsection{Statement of use of AI}
ChatGPT Astra 6 suggested that we use the classical KdV scattering eigenvalue \(\lambda(u)\) to define the spectral speed \(C=-4\lambda(u)\) in our stochastic modulation argument. 
We also used ChatGPT Astra 6.2 to check the proofs and their attribution, revise the exposition, and implement and verify the numerical experiments for the reduced SDE.

\section{Setting and stochastic preliminaries}\label{sec:stochastic}

\subsection{Noise and solution convention}
Let $U=L^2(\R;\R)$, and let $(\Omega,\mathcal F,(\mathcal F_t)_{t\ge0},\Prob)$ be a filtered probability space satisfying the usual completeness and right-continuity conditions. Fix an orthonormal basis $(e_j)$ of $U$ and independent real Brownian motions $(\beta_j)$. The notation
\[
 \mathcal W_t=\sum_{j\ge1}\beta_j(t)e_j
\]
is cylindrical: the series need not converge in $U$. An integrand makes its stochastic integral well defined by sufficiently smoothing or localizing these infinitely many directions.

\begin{assumption}[Noise kernel]\label{ass:noise}
The convolution kernel $k$ is real, nonzero, belongs to $L^1(\R)$, and satisfies
\begin{equation}\label{eq:kernel}
 \sum_{m=0}^2\int_\R e^{2\mu|x|}|\partial_x^m k(x)|^2\dd x<\infty
\end{equation}
for some $\mu>0$.
\end{assumption}

Write $\phi h=k*h$ and $G(u)h=u\phi h$. Equation~\eqref{eq:skdv} means
\[
 u(t)=u(0)+\int_0^t\left[-u_{xxx}(s)-\tfrac12\partial_x(u(s)^2)\right]\dd s
       +\eps\int_0^tG(u(s))\dd\mathcal W_s,
\]
interpreted in a suitable weak or mild sense. We use the global $H^1$ solution supplied by the homogeneous multiplicative-noise theory \cite{deBouardDebussche2007mult}; see also the stochastic KdV foundations \cite{deBouardDebussche1998,deBouardDebusscheTsutsumi1999}. Lemma~\ref{lem:weighted-propagation} supplies the additional propagation of the one-sided weight, and Section~\ref{sec:weighted-energy} proves the uniform stopped estimates.

The formal field $W=\phi\mathcal W$ is stationary in space, with
\[
 \E[W(t,x)W(s,y)]=(t\wedge s)q(x-y),
 \qquad q(z)=\int_\R k(z+r)k(r)\dd r.
\]
It need not be an $L^2(\R)$-valued Wiener process. The product $u\phi$, rather than $\phi$ alone, is the Hilbert--Schmidt integrand used in the equation.

For a Hilbert space $H$, an operator $G:U\to H$ is Hilbert--Schmidt if
\[
 \norm{G}_{\HS(U,H)}^2=\sum_j\norm{Ge_j}_H^2<\infty.
\]
The sum is independent of the orthonormal basis. Parseval's identity gives
\begin{equation}\label{eq:HS-products}
 \begin{aligned}
 \norm{G(f)}_{\HS(U,L^2)}^2&=\norm{k}_2^2\norm{f}_2^2,\\
 \norm{G(f)}_{\HS(U,H^1)}^2
   &=\norm{k}_2^2\norm{f}_{H^1}^2+\norm{k'}_2^2\norm{f}_2^2.
 \end{aligned}
\end{equation}
Indeed $\sum_j|\phi e_j(x)|^2=\norm{k}_2^2$, the corresponding derivative sum equals $\norm{k'}_2^2$, and the mixed sum is $\int kk'=0$. The same calculation yields
\begin{equation}\label{eq:HS-weighted}
 \norm{G(f)}_{\HS(U,H^1_\beta)}^2
 =\norm{k}_2^2\norm{e^{\beta x}f}_{H^1}^2
       +\norm{k'}_2^2\norm{e^{\beta x}f}_2^2,
\end{equation}
where $\norm{f}_{H^1_\beta}=\norm{e^{\beta x}f}_{H^1}$.

\subsection{\Ito\ isometry and maximal estimates}
If $G(s)$ is predictable and
$\E\int_0^t\norm{G(s)}_{\HS(U,H)}^2\dd s<\infty$, then
\[
 M_t=\int_0^tG(s)\dd\mathcal W_s
       :=\lim_{N\to\infty}\sum_{j=1}^N\int_0^tG(s)e_j\dd\beta_j(s)
\]
exists in $L^2(\Omega;H)$. Predictability is the standard nonanticipation condition ensuring that an integrand uses no future Brownian information. The \emph{\Ito\ isometry} is
\begin{equation}\label{eq:ito-isometry}
 \E\norm{M_t}_H^2
     =\E\int_0^t\norm{G(s)}_{\HS(U,H)}^2\dd s.
\end{equation}
It turns a stochastic integral estimate into a deterministic time integral of a squared operator norm. See \cite[Chapter 4]{DaPratoZabczyk2014} and \cite[Chapter 2]{PrevotRockner2007}.

For a supremum over time we use the \emph{Burkholder--Davis--Gundy inequality} (BDG). For $p\ge1$,
\begin{equation}\label{eq:bdg}
 \E\sup_{s\le t}\norm{M_s}_H^p
 \le C_p\E\left(\int_0^t\norm{G(s)}_{\HS(U,H)}^2\dd s\right)^{p/2}.
\end{equation}
In the $p=2$ case, Doob's maximal inequality and \eqref{eq:ito-isometry} give the constant $4$. For the higher moments used in the exit estimates, see \cite{DaPratoZabczyk2014,PrevotRockner2007,Protter2005}.

For example, if $\E\int_0^T\norm{R(s)}_{\HS}^2\dd s\lesssim S$, then
\begin{equation}\label{eq:bdg-example}
 \left\|\sup_{t\le T}\left|\eps\int_0^t R(s)\dd\mathcal W_s\right|\right\|_{L^2(\Omega)}
       \lesssim\eps\sqrt S.
\end{equation}
This square-root gain is the basic source of the speed-error scale in \eqref{eq:target-errors}. In contrast, a deterministic term bounded by a squared radiation norm typically accumulates to $O(S)$.

A stochastic convolution $\int_0^t e^{(t-s)A}G(s)\dd\mathcal W_s$ is generally not a martingale as a function of $t$. Consequently \eqref{eq:bdg} cannot simply be applied to its time supremum with the semigroup suppressed. Section~\ref{sec:weighted-energy} obtains the maximal estimate from a stochastic energy inequality. For related estimates on growing intervals, see \cite{HamsterHupkes2020,WestdorpHupkes2026}.

\subsection{Stopping, probability, and quadratic variation}
For a stopping time $\tau$, a stopped integral means
\[
 \int_0^{t\wedge\tau}G(s)\dd\mathcal W_s
  =\int_0^t\mathbf1_{\{s\le\tau\}}G(s)\dd\mathcal W_s.
\]
The preceding identities apply after localization. A stopped $L^2$ estimate and a high-probability persistence statement are distinct conclusions. For example,
\[
 \Prob\left(\sup_{t\le T\wedge\tau}|X(t)|>Mr\right)
       \le M^{-2}r^{-2}\E\sup_{t\le T\wedge\tau}|X(t)|^2.
\]
Thus a bound $\norm{\sup|X|}_{L^2}\le Cr$ implies $X=\Op(r)$; it does not imply a fixed bound $|X|\le Cr$ with probability tending to one. To obtain such a probability from this estimate alone, the confidence multiplier must tend to infinity. Removing the stop additionally requires control of $\Prob(\tau\le T)$.

For real stochastic integrals $M=\int g\dd\mathcal W$ and $N=\int h\dd\mathcal W$, the quadratic covariation is
\[
 \dd[M,N]_t=\ip{g(t)}{h(t)}_U\dd t.
\]
This is why the two modulation noises must generally be correlated. Replacing them by independent Brownian terms without matching the covariance changes the proposed model.

\subsection{\Ito\ corrections and moving coordinates}
For a sufficiently regular scalar functional $F$ and a strong Hilbert-space semimartingale $\dd u=K(u)\dd t+\eps G(u)\dd\mathcal W$, the \Ito\ formula reads
\begin{equation}\label{eq:ito-functional}
 \begin{split}
 \dd F(u)={}&DF(u)[K(u)]\dd t+\eps DF(u)[G(u)\dd\mathcal W]\\
 &+\frac{\eps^2}{2}\sum_jD^2F(u)[G(u)e_j,G(u)e_j]\dd t.
 \end{split}
\end{equation}
The last term is a second-derivative trace produced by quadratic variation. It is of order $\eps^2$ per unit time and contributes on the scale $S=\eps^2T$.

Equation~\eqref{eq:skdv} is specified in \Ito\ form throughout. A Stratonovich equation with the same written multiplicative coefficient would have an additional \Ito\ drift $\tfrac12\eps^2\norm{k}_2^2u$.

There is a regularity issue in applying \eqref{eq:ito-functional} to KdV: the drift of an $H^1$ solution is not $H^1$-valued. A functional calculation must therefore be justified in a weak framework or by approximation with convergence of all trace and drift terms. Proposition~\ref{prop:h1-spectral-ito} supplies this justification for the spectral functional by mollifying the actual solution and controlling the resulting nonlinear commutator. An arbitrary Galerkin truncation need not preserve the Lax identity, so its deterministic spectral defect must be controlled rather than set equal to zero.

If one instead differentiates $u(t,x+a(t))$ directly, the random change of variables requires an \Ito--Wentzell formula. It includes both the second-order translation term and cross-variation between $u$ and $a$; see \cite{Kunita1990,Krylov2011}. Treating the phase as a scalar chart functional of $u$ is an alternative way of organizing the same correction. The Hessian formula in Appendix~\ref{app:hessian} fixes the pure-soliton phase drift without importing a drift coefficient from a different orthogonality convention.
\section{The one-sided symplectic chart}\label{sec:geometry}

We use the moving orthogonality conditions of Westdorp and Hupkes
\cite[Section~3, (3.1)--(3.2)]{WestdorpHupkes2026}. They arise from the
adjoint generalized eigenspace in the weighted theory of Pego and
Weinstein \cite{PegoWeinstein1994}; see also \cite{MizumachiTzvetkov2014}.
De Bouard and Debussche used this normalized adjoint chain at the initial
speed to define a corrected phase for the limiting stochastic radiation
\cite[Section~4]{deBouardDebussche2009}. We give the formulas in the
normalization \eqref{eq:profile}, including their dependence on the moving
speed, for use in the stochastic coefficient calculations.

\subsection{Weighted spaces and the nondecaying primitive}
Fix a compact interval $I=[c_-,c_+]\Subset(0,\infty)$ containing $c_0$ in its interior, and choose
\begin{equation}\label{eq:beta-choice}
 0<\beta<\min\{\mu/2,\sqrt{c_-}/4\}.
\end{equation}
Set
\[
 H^1_\beta=\{f:e^{\beta y}f\in H^1\},\qquad
 \Xb=H^1\cap H^1_\beta,\qquad
 r_\beta(v)=\norm{e^{\beta y}v}_{H^1}.
\]
The norm on $\Xb$ is the sum of the two displayed norms. The weight is small behind a right-moving soliton and large ahead of it. Its orientation agrees with the fact that linear KdV radiation is transported to the left relative to the soliton.

Let $\dot Q_c=\partial_cQ_c$, and reserve a prime for $y$-derivatives. Direct integration gives
\begin{equation}\label{eq:Z-explicit}
 Z_c(y)=\frac3{\sqrt c}\left[1+\tanh\!\left(\frac{\sqrt c}{2}y\right)\right]
                    +\frac{y}{2c}Q_c(y).
\end{equation}
In particular $Z_c(-\infty)=0$ and $Z_c(+\infty)=6/\sqrt c$. No choice of integration constant makes both limits vanish, because $\int\dot Q_c=6/\sqrt c\ne0$.

For a smooth test function $F$ with $e^{-\beta y}F\in L^2$, define the pairing by
\begin{equation}\label{eq:weighted-pairing}
 \ip{v}{F}=\ip{e^{\beta y}v}{e^{-\beta y}F},\qquad
 |\ip{v}{F}|\le r_\beta(v)\norm{e^{-\beta y}F}_2.
\end{equation}
The pairing is absolutely convergent. Formula~\eqref{eq:Z-explicit} shows that it applies to $Z_c$ and its speed derivatives, uniformly for $c\in I$. We use the word \emph{symplectic} for these two tangent-cotangent pairings; no bounded symplectic form on all of $H^1$ is being asserted.

\subsection{Normalization and adjoint identities}
Define
\begin{equation}\label{eq:DP}
 D_c=\ip{Q_c}{\dot Q_c}=18\sqrt c,
 \qquad P_c=Z_c-c^{-3/2}Q_c.
\end{equation}
The conditions $\ip{v}{Q_c}=\ip{v}{P_c}=0$ are equivalent to \eqref{eq:intro-orth}. The subtraction in \eqref{eq:DP} simplifies the modulation matrix.

\begin{lemma}[Normalized dual chain]\label{lem:dual}
Let $L_c=-\partial_y^2+c-Q_c$ and $A_c=\partial_yL_c$. Then
\begin{equation}\label{eq:dual-identities}
 \begin{gathered}
 L_cQ_c'=0,\qquad L_c\dot Q_c=-Q_c,\\
 A_c^*Q_c=0,\qquad A_c^*P_c=Q_c,\\
 \ip{Q_c'}{P_c}=-D_c,\qquad
 \ip{\dot Q_c}{P_c}=0,\qquad
 \ip{Q_c'}{Q_c}=0.
 \end{gathered}
\end{equation}
The adjoint identities hold distributionally and under the weighted pairings used here.
\end{lemma}
\begin{proof}
Differentiate the profile equation in $y$ and $c$ to obtain the first line. Since $A_c^*=-L_c\partial_y$, $Z_c'=\dot Q_c$, and $L_cQ_c'=0$, the second line follows. Also
\[
 \int Q_c=12\sqrt c,\qquad \norm{Q_c}_2^2=24c^{3/2},
\]
so $D_c=18\sqrt c$. Integration by parts gives
$\ip{Q_c'}{P_c}=-\ip{Q_c}{\dot Q_c-c^{-3/2}Q_c'}=-D_c$.
Finally,
\[
 \ip{\dot Q_c}{Z_c}=\frac12[Z_c^2]_{-\infty}^{+\infty}
        =\frac{18}{c}=c^{-3/2}D_c,
\]
which proves $\ip{\dot Q_c}{P_c}=0$. Boundary terms vanish whenever a localized profile factor is present; the displayed primitive identity accounts explicitly for the nonzero endpoint of $Z_c$.
\end{proof}

\subsection{Existence and smoothness of the chart}
Write $(T_af)(y)=f(y+a)$.
\begin{proposition}[Local phase and geometric speed]\label{prop:chart}
Uniformly for $c_*\in I$ and $a_*\in\R$, a sufficiently small neighborhood of $Q_{c_*}(\cdot-a_*)$, measured by
$\norm{T_{a_*}u-Q_{c_*}}_{\Xb}$, admits unique nearby parameters $(a,c)$ such that
\begin{equation}\label{eq:orth}
 v=T_au-Q_c,\qquad \ip{v}{P_c}=\ip{v}{Q_c}=0.
\end{equation}
The parameter maps are smooth. At $a_*=0$ and $u=Q_c$, their differentials are
\begin{equation}\label{eq:chart-derivatives}
 D\Ac(Q_c)[h]=D_c^{-1}\ip{h}{P_c},\qquad
 D\mathfrak c(Q_c)[h]=D_c^{-1}\ip{h}{Q_c}.
\end{equation}
Here $\Ac$ denotes the position chart and $\mathfrak c$ the geometric-speed chart. The local choice is translation covariant: translating a physical profile by $s$ changes its position by $s$.
\end{proposition}
\begin{proof}
Consider
\[
 F(a,c;u)=\begin{pmatrix}\ip{T_au-Q_c}{P_c}\\\ip{T_au-Q_c}{Q_c}\end{pmatrix}.
\]
One proves smoothness by moving translations onto the smooth test functions, rather than by differentiating $T_au$ repeatedly as an $H^1$ function. The tests and all parameter derivatives needed locally belong to the weighted dual space in \eqref{eq:weighted-pairing}; the speed interval is compact. At $(a,c;u)=(0,c;Q_c)$, Lemma~\ref{lem:dual} gives
\[
 D_{(a,c)}F=\begin{pmatrix}-D_c&0\\0&-D_c\end{pmatrix}.
\]
The Banach-space implicit function theorem therefore applies. Differentiating $F(\Ac(u),\mathfrak c(u);u)=0$ in $u$ proves \eqref{eq:chart-derivatives}. Uniformity follows from the compact speed interval after translating to $a_*=0$. Local uniqueness gives the translation covariance and allows overlapping chart neighborhoods to be glued along a continuous modulated trajectory.
\end{proof}

\begin{remark}[Orbital control in this chart]
The constraint against $P_c$ uses the one-sided weight. For global
$H^1$ control we use the auxiliary translation fixed by the localized
test $Q'_{c_M}$ in Lemma~\ref{lem:orbital-coercivity}; the cost of
transferring that estimate back to $a$ is controlled by $r_\beta(v)$.
\end{remark}

\subsection{Exact deterministic modulation and quadratic decoupling}
\begin{proposition}[Elimination of the linear radiation drift]\label{prop:det-mod}
Let a sufficiently regular deterministic KdV solution have the decomposition \eqref{eq:orth}, with $c\in I$ and $r=r_\beta(v)$ small. Then
\begin{equation}\label{eq:det-matrix}
 \mathsf M(c,v)\begin{pmatrix}\dot a-c\\\dot c\end{pmatrix}
       =\frac12\begin{pmatrix}\ip{v^2}{P_c'}\\\ip{v^2}{Q_c'}\end{pmatrix},
\end{equation}
where
\begin{equation}\label{eq:M}
 \mathsf M(c,v)=
 \begin{pmatrix}
 D_c-\ip{v'}{P_c}&-\ip{v}{\partial_cP_c}\\
 -\ip{v'}{Q_c}&D_c-\ip{v}{\dot Q_c}
 \end{pmatrix}.
\end{equation}
Consequently,
\begin{equation}\label{eq:det-quadratic}
 \begin{pmatrix}\dot a-c\\\dot c\end{pmatrix}
 =\frac1{2D_c}\begin{pmatrix}\ip{v^2}{P_c'}\\\ip{v^2}{Q_c'}\end{pmatrix}
       +O(r^3),
 \qquad |\dot a-c|+|\dot c|\lesssim r^2.
\end{equation}
The constants are uniform on $I$.
\end{proposition}
\begin{proof}
In the moving frame, with $\alpha=\dot a-c$,
\begin{equation}\label{eq:det-v}
 v_t=A_cv+\alpha(Q_c'+v')-\dot c\dot Q_c-\frac12\partial_y(v^2).
\end{equation}
Differentiate the two constraints. The linear terms are
\[
 \ip{A_cv}{Q_c}=0,\qquad
 \ip{A_cv}{P_c}=\ip{v}{Q_c}=0.
\]
Using Lemma~\ref{lem:dual}, the differentiated $P_c$ constraint is
\[
 0=\alpha(-D_c+\ip{v'}{P_c})
       +\dot c\ip{v}{\partial_cP_c}+\frac12\ip{v^2}{P_c'},
\]
and the differentiated $Q_c$ constraint is
\[
 0=\alpha\ip{v'}{Q_c}-\dot c(D_c-\ip{v}{\dot Q_c})
                     +\frac12\ip{v^2}{Q_c'}.
\]
These are exactly \eqref{eq:det-matrix}.

Weighted duality and integration by parts give
$\norm{\mathsf M-D_c\mathrm{Id}}\lesssim r$. Since $P_c'$ and $Q_c'$ decay exponentially and $2\beta<\sqrt{c_-}$,
\[
 |\ip{v^2}{P_c'}|+|\ip{v^2}{Q_c'}|
 \le\norm{e^{\beta y}v}_2^2
       \norm{e^{-2\beta y}(|P_c'|+|Q_c'|)}_\infty
 \lesssim r^2.
\]
Expanding the inverse matrix yields \eqref{eq:det-quadratic}. All integrations can first be made with compactly supported smooth radiation and then extended by the weighted estimates. This establishes the identity whenever the modulated deterministic trajectory is regular enough to differentiate the pairings.
\end{proof}

The cancellation is exact at first order, but it does not make $\dot c$ independent of radiation. If $\int_0^T r_\beta(v)^2\dd t$ is of order $S$, the quadratic speed drift in \eqref{eq:det-quadratic} can still produce an $O(S)$ discrepancy. The spectral coordinate addresses this remaining deterministic contribution.

\subsection{A cutoff realization of the one-sided condition}\label{sec:cutoff}
The weighted formulation avoids a cutoff in the main analysis. For comparison with an $L^2$ test function or a finite computational domain, fix $\chi\in C^\infty(\R)$ with $\chi=1$ on $(-\infty,0]$ and $\chi=0$ on $[1,\infty)$, and set
\[
 P_{c,R}(y)=\chi(y-R)P_c(y),\qquad R\ge1.
\]
Then
\begin{equation}\label{eq:cutoff-bounds}
 \norm{P_{c,R}}_2\lesssim\sqrt{1+R},\qquad
 \norm{e^{-\beta y}(P_c-P_{c,R})}_2\lesssim e^{-\beta R}.
\end{equation}
Furthermore, if $F_{c,R}=A_c^*P_{c,R}-Q_c$, then
\begin{equation}\label{eq:cutoff-defect}
 \norm{e^{-\beta y}F_{c,R}}_2\lesssim e^{-\beta R}.
\end{equation}
Indeed
$F_{c,R}=(\chi(\cdot-R)-1)Q_c+[A_c^*,\chi(\cdot-R)]P_c$.
The commutator is supported in $[R,R+1]$, where the profile and the necessary derivatives of $P_c$ are uniformly bounded. The weighted factor on this interval is $O(e^{-\beta R})$; the remaining profile tail decays faster. This proves \eqref{eq:cutoff-defect}, and \eqref{eq:cutoff-bounds} follows directly from the bounded positive-side limit of $P_c$.

Thus the linear modulation defect created by the cutoff obeys
\begin{equation}\label{eq:cutoff-integrated}
 \int_0^{T\wedge\tau}|\ip{v(t)}{F_{c(t),R}}|\dd t
 \lesssim e^{-\beta R}\sqrt T
       \left(\int_0^{T\wedge\tau}r_\beta(v(t))^2\dd t\right)^{1/2}.
\end{equation}
If the squared radiation integral is $O_{L^p}(S)$, the right side is $O_{L^{2p}}(e^{-\beta R}S/\eps)$. Taking $R=(M/\beta)\Lambda_\eps$, with any fixed $M>2$, suppresses this deterministic defect on the chosen window. The relevant constants are weighted dual norms, which stay bounded; estimating solely by $\norm{P_{c,R}}_2$ would introduce an avoidable $\sqrt R$ loss. A complete cutoff version would also require comparison of its stochastic traces and chart derivatives. No such cutoff theorem is used in the statements below.
\section{Spectral speed and the closed coefficients}\label{sec:coefficients}

\subsection{The isolated eigenvalue and its localized expansion}
Gardner, Greene, Kruskal, and Miura identified the discrete
Schr\"odinger spectrum in the inverse scattering solution of KdV
\cite{GardnerGreeneKruskalMiura1967}. In the normalization
\eqref{eq:skdv}, the speed corresponding to the ground eigenvalue
$\lambda$ is $-4\lambda$.
At $Q_c$, the operator in \eqref{eq:intro-spectral} has the normalized ground state
\begin{equation}\label{eq:ground-state}
 p_c(y)=\frac{c^{1/4}}2\sech\!\left(\frac{\sqrt c}{2}y\right),\qquad
 \Hs(Q_c)p_c=-\frac c4p_c,\qquad
 p_c^2=\frac{Q_c}{12\sqrt c}.
\end{equation}
Its remaining spectrum is $[0,\infty)$. One way to verify this is the factorization
\[
 \Hs(Q_c)+c/4=\mathcal D_c^*\mathcal D_c,\qquad
 \mathcal D_c\mathcal D_c^*=-\partial_y^2+c/4,\qquad
 \mathcal D_c=\partial_y+\frac{\sqrt c}{2}\tanh\!\left(\frac{\sqrt c}{2}y\right).
\]
The zero mode of $\mathcal D_c$ is $p_c$, and the partner operator is free. Let $\Pi_c$ be the orthogonal projection onto $p_c^\perp$, and define the reduced resolvent
\begin{equation}\label{eq:reduced-resolvent}
 R_c^{\rm sp}=\Pi_c(\Hs(Q_c)+c/4)^{-1}\Pi_c,
\end{equation}
where the inverse is taken on $p_c^\perp$ and the operator is zero on the ground state. It is positive and bounded, with norm at most $4/c$.

For real $u$ close to a translate of $Q_c$ in $H^1$, bounded-potential perturbation theory gives a unique simple eigenvalue near $-c/4$, depending smoothly on $u$. Define $\Cs(u)=-4\lambda(u)$ using this eigenvalue. This is the ground eigenvalue in a sufficiently small neighborhood, although additional eigenvalues near zero are not excluded. Standard perturbation theory is described in \cite{Kato1995}.

\begin{proposition}[Spectral derivative and localized quadratic correction]\label{prop:spectral}
Uniformly for $c\in I$,
\begin{equation}\label{eq:spectral-derivatives}
 D\Cs(Q_c)[h]=D_c^{-1}\ip{Q_c}{h},\qquad
 D^2\Cs(Q_c)[h,h]=\frac29\ip{hp_c}{R_c^{\rm sp}(hp_c)}.
\end{equation}
If $\norm{v}_{H^1}$ is sufficiently small and $\ip{v}{Q_c}=0$, then
\begin{equation}\label{eq:localized-spectral}
 \Cs(Q_c+v)-c=\frac19\ip{vp_c}{R_c^{\rm sp}(vp_c)}+\mathcal E_c(v),
 \qquad
 |\mathcal E_c(v)|\lesssim\norm{v}_{H^1}\norm{vp_c}_2^2.
\end{equation}
In particular, for $v\in\Xb$,
\begin{equation}\label{eq:spectral-geometric}
 0\le\Cs(Q_c+v)-c\lesssim\norm{vp_c}_2^2\lesssim r_\beta(v)^2.
\end{equation}
\end{proposition}
\begin{proof}
The Feynman--Hellmann identity gives
$D\lambda(Q_c)[h]=-\ip{p_c^2}{h}/6$. Using \eqref{eq:ground-state} proves the first formula. Differentiating the normalized eigenvalue equation and projecting onto $p_c^\perp$ gives the second derivative, with the two factors $-h/6$ and the factor $-4$ in $\Cs=-4\lambda$; the result is $\frac29\ip{hp_c}{R_c^{\rm sp}(hp_c)}$.

For the localized remainder, put $H_0=\Hs(Q_c)$, $V=-v/6$, $\lambda_0=-c/4$, and $\lambda=\lambda(Q_c+v)$. Since $\ip{v}{Q_c}=0$, $\ip{p_c}{Vp_c}=0$. Decomposing the perturbed eigenvector along $p_c$ and $p_c^\perp$ gives the exact Schur-complement identity
\[
 \lambda-\lambda_0=-\frac1{36}
 \ip{\Pi_c(vp_c)}{[\Pi_c(H_0+V-\lambda)\Pi_c]^{-1}\Pi_c(vp_c)}.
\]
The complementary inverse is uniformly bounded and positive for $\norm{v}_\infty$ small. The min-max principle gives $|\lambda-\lambda_0|\le\norm{V}_\infty$. The resolvent identity then bounds the difference of this inverse and $(H_0-\lambda_0)^{-1}|_{p_c^\perp}$ by $C\norm{v}_\infty$. Multiplication by $-4$ proves \eqref{eq:localized-spectral} and the first inequality in \eqref{eq:spectral-geometric}. Finally $e^{-\beta y}p_c$ is bounded uniformly on $I$, so $\norm{vp_c}_2\lesssim\norm{e^{\beta y}v}_2$. The embedding $H^1(\R)\hookrightarrow L^\infty(\R)$ supplies the stated remainder bound.
\end{proof}

For the decomposition \eqref{eq:orth}, translation invariance therefore yields
\begin{equation}\label{eq:define-C}
 C(t)=\Cs(u(t))=\Cs(Q_{c(t)}+v(t)),\qquad
 |C(t)-c(t)|\lesssim r_\beta(v(t))^2.
\end{equation}
This local correction is stronger than a bound using only the global $H^1$ radiation norm.

\subsection{Isospectral cancellation and the \Ito\ drift}
For a smooth decaying deterministic KdV trajectory, the Lax equation gives
\begin{equation}\label{eq:lax-cancellation}
 D\Cs(u)\left[-u_{xxx}-\tfrac12\partial_x(u^2)\right]=0.
\end{equation}
Indeed, with $q=-u/6$ the equation becomes $q_t-6qq_x+q_{xxx}=0$, and $-\partial_x^2+q$ evolves by a commutator with the skew-adjoint third-order Lax operator \cite{Lax1968}. Pairing this commutator against an eigenfunction shows that its eigenvalue derivative vanishes.

Proposition~\ref{prop:h1-spectral-ito} justifies the functional \Ito\
formula in $H^1$. Together with \eqref{eq:lax-cancellation} it gives
\begin{equation}\label{eq:spectral-ito-target}
 \begin{split}
 \dd C={}&\eps D\Cs(u)[u\phi\dd\mathcal W_t]\\
 &+\frac{\eps^2}{2}\sum_jD^2\Cs(u)[u\phi e_j,u\phi e_j]\dd t.
 \end{split}
\end{equation}
Both the martingale coefficient and the trace depend on the full
profile. Proposition~\ref{prop:localized-spectral-coefficients}
estimates their difference from their values at the soliton.

At the pure soliton, set $h_{c,j}=Q_c\phi e_j$. Formula~\eqref{eq:spectral-derivatives} fixes the speed drift:
\begin{equation}\label{eq:bsp}
 b_{\rm sp}(c)=\frac19\sum_j
      \ip{h_{c,j}p_c}{R_c^{\rm sp}(h_{c,j}p_c)}\ge0.
\end{equation}
Boundedness of the reduced resolvent and \eqref{eq:HS-products} give
absolute convergence and basis independence. Proposition~\ref{prop:pure-drift-comparison}
identifies the contribution distinguishing this trace from the
geometric-speed drift.

\begin{remark}[Comparison with Garnier's scalar multiplier]\label{rem:scalar-garnier}
For spatially uniform scalar noise the pure coefficients can be computed
explicitly. This is a separate forcing model from Assumption~\ref{ass:noise}.
For $x$ near $1$, let $\nu(\nu+1)=2x$ with $\nu(1)=1$. Direct substitution
of $\sech^\nu(\sqrt c\,y/2)$ in the eigenvalue equation for $\Hs(xQ_c)$
gives a positive eigenfunction with eigenvalue $-c\nu^2/4$. Thus
$\Cs(xQ_c)=c\nu(x)^2$ and differentiation gives
\[
 D\Cs(Q_c)[Q_c]=\tfrac43c,\qquad
 \tfrac12 D^2\Cs(Q_c)[Q_c,Q_c]=\tfrac4{27}c.
\]
For a scalar Stratonovich noise $\eps u\circ\dd W_t$, the conversion
to \Ito\ form contributes an additional
$\frac12 D\Cs(Q_c)[Q_c]=2c/3$ to the drift coefficient on the slow
time $s=\eps^2t$. The resulting pure spectral diffusion and drift
coefficients are $4c/3$ and $22c/27$. Garnier's white-noise formula
\cite[Section~6, (74)--(76)]{Garnier2001}, at unit intensity, gives
$\kappa(s)=\kappa(0)\exp(-s/27+2W_s/3)$ for the limiting
scattering parameter. Applying \Ito's formula to $C=4\kappa^2$
gives exactly these two coefficients.
\end{remark}

\subsection{Pure-soliton covariance and phase drift}
By \eqref{eq:chart-derivatives} and \eqref{eq:spectral-derivatives}, the joint first-order noise map for $(a,C)$ at a centered soliton is
\begin{equation}\label{eq:K0}
 K_0(c)h=\frac1{D_c}
   \begin{pmatrix}\ip{Q_cP_c}{\phi h}\\\ip{Q_c^2}{\phi h}\end{pmatrix}.
\end{equation}
Define
\begin{equation}\label{eq:q-matrix}
 q_0(c)=K_0(c)K_0(c)^*
     =\begin{pmatrix}q_{aa}(c)&q_{aC}(c)\\q_{aC}(c)&q_{CC}(c)\end{pmatrix}.
\end{equation}
In particular,
\begin{equation}\label{eq:g}
 q_{CC}(c)=\frac{\norm{\phi^*(Q_c^2)}_2^2}{D_c^2}.
\end{equation}
This is also the frozen speed variance of de Bouard and Debussche
\cite[(3.12), Sections~4--5]{deBouardDebussche2009}. Their corrected
phase uses the other normalized adjoint test. In our (phase, speed)
ordering and translation convention, their pair $(B_2,B_1)$ has
covariance $q_0(c_0)$ at unit time. In the normalization $\varphi_c=Q_c/2$ of
Westdorp and Hupkes, the geometric-speed diffusion test
$(2/9)c^{-1/2}\varphi_c^2$ in
\cite[(1.7), Section~3]{WestdorpHupkes2026} equals $Q_c^2/D_c$.
The pure phase drift is defined in the same phase chart:
\begin{equation}\label{eq:gamma}
 \gamma_a(c)=\frac12\sum_jD^2\Ac(Q_c)[h_{c,j},h_{c,j}].
\end{equation}
The chart Hessian is bounded on $\Xb$, and the noise is Hilbert--Schmidt into $\Xb$ at $Q_c$ by \eqref{eq:HS-products}--\eqref{eq:HS-weighted}. Hence the sum is absolutely convergent and basis independent. Appendix~\ref{app:hessian} expresses it entirely in profile integrals.

\begin{proposition}[Nondegeneracy and a smooth covariance factor]\label{prop:cov}
Under Assumption~\ref{ass:noise}, $q_0(c)$ is positive definite for every $c\in I$, uniformly on $I$. A smooth factor with the speed noise in its first Brownian direction is
\begin{equation}\label{eq:sigma-factor}
 \Sigma_0(c)=
 \begin{pmatrix}
 q_{aC}/\sqrt{q_{CC}}&\sqrt{q_{aa}-q_{aC}^2/q_{CC}}\\
 \sqrt{q_{CC}}&0
 \end{pmatrix},\qquad \Sigma_0\Sigma_0^*=q_0.
\end{equation}
Its rows are $\sigma_a$ and $\sigma_C$ in \eqref{eq:intro-effective}.
\end{proposition}
\begin{proof}
The functions $Q_cP_c$ and $Q_c^2$ are linearly independent and decay exponentially at both ends. To check independence, divide a proposed relation by the strictly positive $Q_c$: the positive-side limit of $P_c$ is nonzero whereas $Q_c$ tends to zero. Their Fourier transforms, and those of their linear combinations, extend analytically to a strip. Since $k\in L^1$ is nonzero, its continuous Fourier transform is nonzero on some nonempty open interval. If $\phi^*f=0$ for a linear combination $f$ of these two functions, its Fourier transform vanishes on that interval, hence identically by analyticity. Thus $f=0$. This proves positive definiteness of their Gram matrix after applying $\phi^*$. Smooth dependence on $c$ and compactness give uniform positivity. Formula~\eqref{eq:sigma-factor} now follows by multiplication and has smooth entries.
\end{proof}

\begin{lemma}[Strict increase of the speed variance]\label{lem:increasing-speed-variance}
Under Assumption~\ref{ass:noise} we have $q'_{CC}(c)>0$ for every $c>0$.
\end{lemma}
\begin{proof}
Use $\widehat f(\xi)=\int e^{-iy\xi}f(y)\dd y$ and put
$A(c,\xi)=\widehat{Q_c^2}(\xi)/D_c$. The Fourier transform of
$\sech^2$ gives $\widehat Q_c(\xi)=12\pi\xi/\sinh(\pi\xi/\sqrt c)$.
The profile equation gives
$A=2(c+\xi^2)\widehat Q_c/(18\sqrt c)>0$.
With $z=\pi|\xi|/\sqrt c$ and $z\coth z=1$ at $z=0$,
\[
  c\,\partial_c\log A(c,\xi)
       =\frac{c}{c+\xi^2}+\frac12(z\coth z-1)>0 .
\]
On each compact positive speed interval, $A$ and $\partial_c A$
decay exponentially in $|\xi|$, up to polynomial factors.
We may therefore differentiate the Parseval identity
$q_{CC}(c)=(2\pi)^{-1}\int|\widehat k|^2 A(c,\xi)^2\dd\xi$.
The integrand of its derivative is positive wherever
$\widehat k\ne0$, a set of positive measure.
\end{proof}

\subsection{What spatial homogeneity removes}
For a fixed cylindrical noise realization, translating a profile also translates the noise directions. In our convention the exact parameter noise map has the form
\begin{equation}\label{eq:covariance-translation}
 \Kb(a,c,v)=\Kb(0,c,v)T_a,
\end{equation}
because convolution commutes with translations. Since $T_a$ is unitary on $U$, the covariance $\Kb\Kb^*$ is independent of $a$. The operator $\Kb$ itself is not independent of $a$.

Suppose the exact two-parameter martingale has covariance $q(t)=\Kb(t)\Kb(t)^*$, positive definite before a stop. Let $\Sigma(t)$ be its factor in the convention \eqref{eq:sigma-factor}. Then
\begin{equation}\label{eq:adapted-B}
 \dd B_t=\Sigma(t)^{-1}\Kb(t)\dd\mathcal W_t,\qquad t<\tau,
\end{equation}
has quadratic covariation $\dd[B]_t=\mathrm{Id}_2\dd t$. L\'evy's characterization makes it Brownian motion up to the stop. It can be continued after the stop by a fixed two-dimensional projection of the remaining cylindrical increments, as constructed in \eqref{eq:global-adapted-B}. This is a standard martingale construction; see \cite{Protter2005,DaPratoZabczyk2014}.

The effective equation is coupled to the PDE using this $B$, rather than by erasing $a$ from a fixed cylindrical integrand. Independently simulating \eqref{eq:intro-effective} with any standard two-dimensional Brownian motion gives its law. Predicting a particular PDE realization by a pathwise comparison requires the adapted coupling \eqref{eq:adapted-B}. These are two distinct uses of the same reduced equation.
\section{Exact stochastic equations at energy regularity}\label{sec:exact-ito}

We derive the equations for the spectral speed and the moving phase
along the energy-space solution. The spectral equation holds throughout
the unweighted spectral tube. The geometric equations use the additional
one-sided regularity established in Lemma~\ref{lem:weighted-propagation}.

\subsection{Spectral calculus away from the pure soliton}
For a real potential $w$ in a sufficiently small $H^1$ neighborhood of $Q_c$, write
\[
 \lambda_w=\lambda(w),\qquad \psi_w>0,\qquad \norm{\psi_w}_2=1,
 \qquad \mathsf P_w=\psi_w\otimes\psi_w,\qquad \Pi_w=1-\mathsf P_w .
\]
The symbol $\mathsf P_w$ denotes an orthogonal spectral projection; it is distinct from the symplectic test $P_c$. Set
\begin{equation}\label{eq:general-resolvent}
 R_w=\Pi_w(\Hs(w)-\lambda_w)^{-1}\Pi_w,\qquad
 \mathsf S_w=\Hs(w)-\lambda_w+\mathsf P_w .
\end{equation}
The inverse defining $R_w$ is restricted to $\psi_w^\perp$. Equivalently,
$R_w=\mathsf S_w^{-1}-\mathsf P_w$. All these definitions commute with spatial translation.

\begin{lemma}[Differentiation and the weak Lax cancellation]\label{lem:general-spectral}
In this neighborhood $\Cs$ is a smooth real functional on $H^1$, and
\begin{equation}\label{eq:general-spectral-derivatives}
 \begin{split}
 D\Cs(w)[h]&=\frac23\ip{\psi_w^2}{h},\\
 D^2\Cs(w)[h,\ell]&=\frac29\ip{h\psi_w}{R_w(\ell\psi_w)} .
 \end{split}
\end{equation}
The eigenfunction is in $H^3$; its $H^3$ norm and the relevant functional derivatives are bounded uniformly on the tube, including all translates. The first derivative extends to $H^{-2}$ by duality. With
$K(w)=-w'''-\frac12(w^2)'$, it satisfies
\begin{equation}\label{eq:weak-lax}
 D\Cs(w)[K(w)]=0,\qquad w\in H^1
\end{equation}
throughout this spectral neighborhood.
\end{lemma}
\begin{proof}
Multiplication by an $H^1$ function is a bounded perturbation on $L^2$, and the eigenvalue is uniformly separated from the remaining spectrum. The Riesz projection construction and the normalized eigenvalue equation therefore give smooth dependence on $w$; see \cite{Kato1995}. Choosing the positive normalization fixes the sign. Differentiating the eigenvalue equation gives
\[
 D\lambda_w[h]=-\frac16\ip{\psi_w^2}{h},\qquad
 D\psi_w[h]=\frac16R_w(h\psi_w).
\]
The second identity uses $\ip{\psi_w}{D\psi_w[h]}=0$. Differentiating $\Cs=-4\lambda$ proves \eqref{eq:general-spectral-derivatives}; self-adjointness of $R_w$ makes the Hessian symmetric.

The eigenvalue equation is
$\psi_w''=-(w/6+\lambda_w)\psi_w$. Boundedness of $w$, the spectral gap, and normalization first give a uniform $H^2$ bound. Since $H^1(\R)$ is an algebra, the same equation then gives an $H^3$ bound. These arguments also give continuity of $\psi_w$ in $H^3$ when $w$ varies in $H^1$. Hence $\psi_w^2$ is a continuous $H^2$ gradient for the functional, so its pairing with $K(w)\in H^{-2}$ is defined.

For real smooth rapidly decaying $w$, the Lax identity \eqref{eq:lax-cancellation} gives \eqref{eq:weak-lax}. Approximate any potential in the tube by real Schwartz potentials in $H^1$. The map $w\mapsto K(w)$ is continuous from $H^1$ to $H^{-2}$, and $w\mapsto\psi_w^2$ is continuous into $H^2$. The approximants remain in the spectral neighborhood. Passing to the limit proves \eqref{eq:weak-lax}, without requiring the actual stochastic solution to possess three spatial derivatives.
\end{proof}

\subsection{An \Ito\ formula obtained by smoothing the actual solution}
The next proposition concerns the global $H^1$ solution, stopped while its eigenvalue stays in the preceding spectral neighborhood. It does not require a weighted radiation norm.

\begin{proposition}[Spectral \Ito\ equation in $H^1$]\label{prop:h1-spectral-ito}
Let $u$ be the continuous adapted $H^1$ solution of \eqref{eq:skdv}, and stop it in a sufficiently small spectral tube with speeds in a compact positive interval. Then $C=\Cs(u)$ is a continuous semimartingale and
\begin{equation}\label{eq:proved-spectral-ito}
 \dd C
 =\frac{2\eps}{3}\ip{u\psi_u^2}{\phi\dd\mathcal W}
   +\frac{\eps^2}{9}\sum_j
       \ip{u\psi_u\phi e_j}{R_u(u\psi_u\phi e_j)}\dd t .
\end{equation}
The trace converges absolutely and is basis independent. In particular, the drift displayed here is nonnegative and there is no deterministic KdV contribution.
\end{proposition}
\begin{proof}
Fix a nonnegative even $\rho\in C_c^\infty(\R)$ with $\int\rho=1$, and put
$J_n f=\rho_n*f$, $\rho_n(x)=n\rho(nx)$. Smooth the actual solution:
$u_n=J_nu$. For each fixed $n$, this is an $H^1$-valued strong semimartingale, since $J_n:H^{-2}\to H^1$ is bounded. Its equation is
\begin{equation}\label{eq:mollified-equation}
 \dd u_n=\big[K(u_n)-\tfrac12\partial_x d_n\big]\dd t
          +\eps J_nG(u)\dd\mathcal W,\qquad
 d_n=J_n(u^2)-(J_nu)^2 .
\end{equation}
The linear third derivative commutes with $J_n$. On each bounded stopped interval, continuity of the $H^1$ path, compactness of its image, and the algebra property of $H^1$ imply
\begin{equation}\label{eq:commutator-convergence}
 \sup_t\norm{u_n(t)-u(t)}_{H^1}\longrightarrow0,\qquad
 \sup_t\norm{d_n(t)}_{H^1}\longrightarrow0
\end{equation}
almost surely. For each $n$, additionally stop $u_n$ at exit from a slightly larger spectral tube. By this uniform convergence, the auxiliary exit does not precede the original bounded stop for all sufficiently large $n$, almost surely. Thus the formulas for $\Cs(u_n)$ can be passed to the original stopped interval.

Apply the Hilbert-space \Ito\ formula \eqref{eq:ito-functional} to $\Cs(u_n)$. By Lemma~\ref{lem:general-spectral}, the drift $K(u_n)$ cancels exactly. The additional drift is
\begin{equation}\label{eq:spectral-commutator}
 D\Cs(u_n)[-\tfrac12\partial_xd_n]
       =\frac13\ip{(\psi_{u_n}^2)'}{d_n},
\end{equation}
which tends uniformly to zero after localization by \eqref{eq:commutator-convergence} and the uniform eigenfunction bound. Thus no isospectral property of an approximate evolution has been assumed.

The noise converges in the appropriate Hilbert--Schmidt norm:
\[
 J_nG(u)\longrightarrow G(u)
       \quad\hbox{in }\HS(U,H^1).
\]
Indeed $J_n\to1$ strongly and is uniformly bounded on $H^1$; apply this fact to the Hilbert--Schmidt operator $G(u)$, using \eqref{eq:HS-products}. The convergence is uniform along each compact time path. Continuity of the first two derivatives of $\Cs$ now yields convergence of both the martingale integrands and the Hessian traces. For clarity, the trace convergence uses the elementary estimate
\begin{align*}
 \left|\sum_j H_n[A_ne_j,A_ne_j]-\sum_j H[Ae_j,Ae_j]\right|
 \le{}&\norm{H_n-H}\norm{A_n}_{\HS}^2\\
 &+\norm{H}(\norm{A_n}_{\HS}+\norm{A}_{\HS})\norm{A_n-A}_{\HS},
\end{align*}
for bounded bilinear forms $H_n,H$ and Hilbert--Schmidt operators $A_n,A$. Dominated convergence after bounded localization treats the time integrals; the \Ito\ isometry and the $p=2$ BDG inequality treat the martingale suprema. Removing these auxiliary localizations proves the formula up to the stated stop.

Substitution of \eqref{eq:general-spectral-derivatives} gives \eqref{eq:proved-spectral-ito}. If $F_u h=u\psi_u\phi h$, then
\[
 \norm{F_u}_{\HS(U,L^2)}^2=\norm{k}_2^2\norm{u\psi_u}_2^2<\infty,
 \qquad
 \sum_j\ip{F_ue_j}{R_uF_ue_j}
       =\operatorname{Tr}(F_u^*R_uF_u).
\]
This proves absolute convergence and basis independence. The continued eigenvalue is the ground eigenvalue, so $R_u$ is positive. This proves the final assertion.
\end{proof}

The proof uses only energy-space regularity of $u$ and the spatial smoothing present in the prescribed noise. The standard \Ito\ formula and maximal estimates are used after spatial mollification, where their strong-semimartingale hypotheses hold. The commutator \eqref{eq:spectral-commutator} is the explicit check needed to return to the original PDE.

\subsection{First and second derivatives of the full phase chart}
We now work in the centered frame $w=Q_c+v$, with both constraints \eqref{eq:orth}. Write
\[
 f_1=P_c,\qquad f_2=Q_c,\qquad
 \ell_c(h)=\begin{pmatrix}\ip{h}{f_1}\\ \ip{h}{f_2}\end{pmatrix}.
\]
Dots on $f_i$ denote $c$-derivatives. All spatial derivatives in the next formulas can be moved onto the smooth tests. Define the three column vectors
\begin{equation}\label{eq:chart-second-vectors}
 \begin{split}
 \mathsf A_i&=\ip{w''}{f_i}=\ip{w}{f_i''},\\
 \mathsf B_i&=\ip{w'}{\dot f_i}=-\ip{w}{\dot f_i'},\\
 \mathsf D_i&=-\ip{\ddot Q_c}{f_i}
             -2\ip{\dot Q_c}{\dot f_i}+\ip{v}{\ddot f_i},
             \qquad i=1,2 .
 \end{split}
\end{equation}
The column $\mathsf D$ here is unrelated to the scalar $D_c=18\sqrt c$.

\begin{lemma}[Full chart Hessian]\label{lem:full-chart-hessian}
Let $\Theta=(\Ac,\mathfrak c)$, and evaluate its derivatives at the centered profile $w$. With $\mathsf M=\mathsf M(c,v)$ from \eqref{eq:M}, put
\begin{equation}\label{eq:general-chart-first}
 \xi_h=\begin{pmatrix}\alpha_h\\\kappa_h\end{pmatrix}
      =\mathsf M^{-1}\ell_c(h)=D\Theta(w)[h].
\end{equation}
Then
\begin{equation}\label{eq:full-chart-hessian}
 D^2\Theta(w)[h,h]
 =\mathsf M^{-1}\left\{
   \alpha_h^2\mathsf A+2\alpha_h\kappa_h\mathsf B
    +\kappa_h^2\mathsf D
    +2\alpha_h\ell_c(h')+2\kappa_h\dot\ell_c(h)\right\},
\end{equation}
where $\dot\ell_c(h)_i=\ip{h}{\dot f_i}$ and
$\ell_c(h')_i=-\ip{h}{f_i'}$. The formulas extend continuously to $h\in\Xb$.
\end{lemma}
\begin{proof}
For the constraint map $F(a,c;u)$ used in Proposition~\ref{prop:chart},
$F_\theta=-\mathsf M$ and $F_u[h]=\ell_c(h)$ at the centered profile. Its second parameter derivatives are exactly
$F_{aa}=\mathsf A$, $F_{ac}=\mathsf B$, and $F_{cc}=\mathsf D$. Its mixed derivatives are $F_{au}[h]=\ell_c(h')$ and $F_{cu}[h]=\dot\ell_c(h)$, and $F_{uu}=0$. Differentiate
$F(\Theta(u);u)=0$ once and twice. Solving the resulting two linear systems gives \eqref{eq:general-chart-first} and \eqref{eq:full-chart-hessian}. Smoothness and weighted duality justify this calculation without differentiating the $H^1$ profile twice in space.
\end{proof}

To state the trace explicitly, define row operators on $U$ by
\begin{equation}\label{eq:chart-noise-rows}
 \begin{gathered}
 L_i h=\ip{w\phi h}{f_i},\qquad
 L_i^{(x)}h=-\ip{w\phi h}{f_i'},\qquad
 L_i^{(c)}h=\ip{w\phi h}{\dot f_i},\\
 N=\mathsf M^{-1}L,\qquad
 n=NN^*,\qquad
 N=\begin{pmatrix}N_a\\N_c\end{pmatrix}.
 \end{gathered}
\end{equation}
For row operators, $\ip{N_a}{L_i^{(x)}}_U$ means the scalar product of their Riesz representatives in $U$. Set
\begin{equation}\label{eq:chart-gamma-full}
 \begin{split}
 t_i&=\ip{N_a}{L_i^{(x)}}_U+\ip{N_c}{L_i^{(c)}}_U,\\
 \Gamma(c,v)&=\mathsf M^{-1}
       \left\{\tfrac12 n_{aa}\mathsf A+n_{ac}\mathsf B
                     +\tfrac12 n_{cc}\mathsf D+t\right\},\\
 V(c,v)&=\frac12\mathsf M^{-1}
            \begin{pmatrix}\ip{v^2}{P_c'}\\\ip{v^2}{Q_c'}\end{pmatrix}.
 \end{split}
\end{equation}
Formula \eqref{eq:full-chart-hessian}, summed with $h=w\phi e_j$, gives
\begin{equation}\label{eq:full-trace-identification}
 \Gamma(c,v)=\frac12\sum_jD^2\Theta(w)[w\phi e_j,w\phi e_j].
\end{equation}
Thus \eqref{eq:chart-gamma-full} evaluates the trace using finitely many scalar products of row operators. It includes the changes in both phase and geometric speed, including their covariance $n_{ac}$.

\subsection{Exact phase equation and its localized coefficients}
Let $\zeta$ denote a lifetime on which $u$ has a continuous $\Xb$ path and the weighted chart can be followed, after translation to its center. On compact subintervals one may additionally stop $|a|$ to use fixed laboratory weighted spaces; translations only multiply the exponential norm by a finite factor. Existence and long-time continuation of this weighted path will be addressed with the radiation estimates. The following assertions are local statements on this lifetime.

\begin{proposition}[Exact geometric modulation and phase trace]\label{prop:geometric-ito}
Before the exit from a sufficiently small $\Xb$ tube with $c\in I$, the chart parameters are continuous semimartingales satisfying
\begin{equation}\label{eq:exact-geometric-ito}
 \dd\begin{pmatrix}a\\c\end{pmatrix}
   =\left[\begin{pmatrix}c\\0\end{pmatrix}+V(c,v)\right]\dd t
     +\eps N(c,v)T_a\dd\mathcal W
     +\eps^2\Gamma(c,v)\dd t .
\end{equation}
With $r=r_\beta(v)$, uniformly on this tube,
\begin{equation}\label{eq:geometric-coefficient-bounds}
 \begin{split}
 |V(c,v)|&\le Cr^2,\\
 \norm{N(c,v)-K_0(c)}_{\HS(U,\R^2)}
                +|\Gamma(c,v)-\Gamma(c,0)|&\le Cr,\\
 \Gamma_a(c,0)&=\gamma_a(c).
 \end{split}
\end{equation}
Here $K_0$ has the phase and first-order speed rows in \eqref{eq:K0}; those speed rows agree for geometric and spectral speed only at the pure soliton.
\end{proposition}
\begin{proof}
First justify the \Ito\ calculation. The compactly supported mollifiers used above are uniformly bounded and converge strongly on $H^1_\beta$, as follows by conjugating convolution with $e^{\beta x}$. They also smooth $H^{-2}_\beta$ into $H^1_\beta$. On a bounded local chart, the smoothed actual solution is therefore a strong semimartingale in the Hilbert space $\Xb$, equipped with its equivalent squared-sum norm. Products satisfy
\[
 \norm{fg}_{H^1_\beta}
 \le C\big(\norm{f}_{H^1}\norm{g}_{H^1_\beta}
            +\norm{g}_{H^1}\norm{f}_{H^1_\beta}\big),
 \qquad f,g\in\Xb .
\]
Consequently $d_n=J_n(u^2)-(J_nu)^2\to0$ in $\Xb$, uniformly on each compact localized time path. The only deterministic error in the chart \Ito\ formula is $D\Theta(u_n)[-\partial_xd_n/2]$. Formula \eqref{eq:general-chart-first}, with the derivative transferred to $f_i$, bounds it by $C\norm{e^{\beta x}d_n}_2$ on a fixed local chart. It tends to zero. The Hilbert--Schmidt convergence of $J_nG(u)$ in $\Xb$ follows from \eqref{eq:HS-products}--\eqref{eq:HS-weighted}. The same trace estimate and BDG argument as in Proposition~\ref{prop:h1-spectral-ito} pass the first and second derivative terms to the limit. Local uniqueness patches the charts; temporary bounds on $|a|$ can then be removed.

The deterministic contribution is $(c,0)^T+V$, by Proposition~\ref{prop:det-mod}. Translation of the noise gives $N(c,v)T_a$, because $T_aG(u)h=w\phi(T_ah)$. Lemma~\ref{lem:full-chart-hessian} gives the trace $\Gamma$. This proves \eqref{eq:exact-geometric-ito}.

We give the localized estimates in detail because replacing $r$ by $\norm{v}_{H^1}$ would be insufficient for the intended comparison. Each test needed in \eqref{eq:chart-second-vectors} and \eqref{eq:chart-noise-rows} satisfies, uniformly in $c$,
\begin{equation}\label{eq:test-dual-bounds}
 \norm{e^{-\beta y}\partial_y^m\partial_c^\ell f_i}_{L^2\cap L^\infty}<C
 \quad (0\le m,\ell\le2).
\end{equation}
These follow from the explicit profiles: the tests decay exponentially on the negative side and are bounded on the positive side. Higher derivatives used for mollification have the same property. In particular,
\[
 |\mathsf A(c,v)-\mathsf A(c,0)|
  +|\mathsf B(c,v)-\mathsf B(c,0)|
  +|\mathsf D(c,v)-\mathsf D(c,0)|\le Cr .
\]
The corresponding estimate $\norm{\mathsf M-D_c\mathrm{Id}}\le Cr$ was proved in Section~\ref{sec:geometry}; hence $\mathsf M^{-1}$ is uniformly bounded and differs from $D_c^{-1}\mathrm{Id}$ by $O(r)$.

The Riesz representative of $L_i(c,v)-L_i(c,0)$ is $\phi^*(vf_i)$. Young's convolution inequality and \eqref{eq:test-dual-bounds} give
\begin{equation}\label{eq:row-localization}
 \norm{\phi^*(v f_i)}_2
 \le\norm{k}_1\norm{vf_i}_2
 \le C\norm{e^{\beta y}v}_2 .
\end{equation}
The identical argument applies with $f_i'$ or $\dot f_i$. Thus the rows $L,L^{(x)},L^{(c)}$ are bounded and differ from their pure values by $O(r)$. Matrix multiplication gives $N-K_0=O(r)$, $n-n(c,0)=O(r)$, and then $\Gamma-\Gamma(c,0)=O(r)$ in \eqref{eq:chart-gamma-full}. All row scalar products, and hence all trace sums, converge absolutely. The bound on $V$ follows from the quadratic estimate in Proposition~\ref{prop:det-mod}. Finally \eqref{eq:full-trace-identification} at $v=0$ is precisely the definition of $\gamma_a$ in \eqref{eq:gamma}.
\end{proof}

\subsection{The radiation equation, including the cross variation}\label{sec:radiation-equation}
The exact identities also give the equation estimated in Section~\ref{sec:weighted-energy}. We record its stochastic corrections before applying an energy estimate.

Define a centered cylindrical Brownian motion by
\[
 \dd\widetilde{\mathcal W}=T_a\dd\mathcal W .
\]
Precisely, for $h\in U$,
$\widetilde{\mathcal W}_t(h)=\int_0^t\ip{T_{a(s)}^*h}{\dd\mathcal W_s}_U$.
Every finite collection has Brownian covariations $\ip{h,\ell}_Ut$, because $T_a$ is a predictable unitary operator; L\'evy's characterization supplies the cylindrical Brownian property. Set
\begin{equation}\label{eq:radiation-noise}
 \begin{split}
 \chi_a&=\phi N_a^*,\\
 \mathcal G_v h&=w\phi h+w' N_a h-\dot Q_c N_c h,\\
 \mathcal J(c,v)&=\Gamma_a w'-\Gamma_c\dot Q_c
       +\tfrac12n_{aa}w''-\tfrac12n_{cc}\ddot Q_c
       +\partial_y(w\chi_a).
 \end{split}
\end{equation}
Here $N_a^*\in U$ is the Riesz representative of the phase-noise row, and all quantities are evaluated at $(c,v)$.

\begin{proposition}[Weak moving-frame equation]\label{prop:radiation-equation}
On the weighted chart lifetime,
\begin{equation}\label{eq:full-radiation-spde}
 \dd v=\left[A_cv-\tfrac12\partial_y(v^2)
          +V_aw'-V_c\dot Q_c+\eps^2\mathcal J(c,v)\right]\dd t
          +\eps\mathcal G_v\dd\widetilde{\mathcal W}.
\end{equation}
The equation holds after pairing with every smooth compactly supported test. In particular, the term $\partial_y(w\chi_a)$ must be retained.
\end{proposition}
\begin{proof}
Apply the \Ito--Wentzell formula to $w(t,y)=u(t,y+a(t))$,
as in the stochastic KdV modulation calculation of
Westdorp and Hupkes \cite[Section~3]{WestdorpHupkes2026}.
See \cite{Kunita1990} for the classical formula and
\cite{Krylov2011} for its distribution-valued version.
At the present regularity it can be applied to the smooth test-function field
$F(t,z)=\ip{u(t)}{\vartheta(\cdot-z)}$, with $\vartheta\in C_c^\infty$, or justified by the spatial mollification already used above. Besides $K(w)\dd t$, $w'\dd a$, and the translated noise, it gives
\[
 \frac{\eps^2}{2}n_{aa}w''\dd t
   +\eps^2\sum_j N_ae_j\,\partial_y(w\phi e_j)\dd t
 =\frac{\eps^2}{2}n_{aa}w''\dd t
   +\eps^2\partial_y(w\chi_a)\dd t .
\]
The second summand is the cross variation of the field noise and the random translation. Apply the ordinary scalar \Ito\ formula to $Q_{c(t)}$:
\[
 \dd Q_c=\dot Q_c(V_c+\eps^2\Gamma_c)\dd t
       +\eps\dot Q_cN_c\dd\widetilde{\mathcal W}
       +\frac{\eps^2}{2}n_{cc}\ddot Q_c\dd t .
\]
Subtract, insert \eqref{eq:exact-geometric-ito}, and use
$K(w)+cw'=A_cv-\partial_y(v^2)/2$. This proves \eqref{eq:full-radiation-spde}. The noise belongs to $\HS(U,L^2\cap L^2_\beta)$ on bounded chart intervals, and every displayed drift is a distribution of order at most two. The test-function argument is therefore valid at $H^1\cap H^1_\beta$ regularity.
\end{proof}

One useful check is the exact first-order constraint identity
\[
 \ip{\mathcal G_vh}{f_i}+\ip{v}{\dot f_i}N_ch=0,\qquad i=1,2.
\]
It follows from $\mathsf MN=L$. The second-order terms are accounted
for by $\Gamma$. Equation~\eqref{eq:full-radiation-spde} uses the
geometric speed $c$ defining the two orthogonality conditions. Since
$w'$ is only in $L^2$ at energy regularity, we justify the weighted
energy calculation by regularization in Proposition~\ref{prop:weighted-energy}.
\section{Localized spectral coefficients and closed stochastic dynamics}\label{sec:localized-closure}

\subsection{A weighted estimate for the unperturbed reduced resolvent}
The phase estimates used only explicit test functions. The spectral drift contains a resolvent, so a further localization argument is needed. An operator-norm estimate
$\norm{R_{Q_c+v}-R_{Q_c}}_{2\to2}\lesssim\norm{v}_{H^1}$
would not yield the desired bound in $r_\beta(v)$.

\begin{lemma}[Weighted unperturbed inverse]\label{lem:weighted-resolvent}
For $s=\beta$ or $s=-\beta$, uniformly for $c\in I$,
\begin{equation}\label{eq:weighted-resolvent}
 \norm{e^{sy}R_c^{\rm sp}f}_{H^2}
       +\norm{e^{sy}\mathsf P_{Q_c}f}_{H^2}
       \le C\norm{e^{sy}f}_2 .
\end{equation}
Consequently the same bound holds for
$\mathsf S_{Q_c}^{-1}=R_c^{\rm sp}+\mathsf P_{Q_c}$.
\end{lemma}
\begin{proof}
Put $\kappa=\sqrt c/2$,
$\mathcal D=\partial_y+\kappa\tanh(\kappa y)$, and
$\mathcal B=-\partial_y^2+\kappa^2$.
The factorization in Section~\ref{sec:coefficients} gives
$\mathcal D\mathcal D^*=\mathcal B$, and
\begin{equation}\label{eq:explicit-reduced-resolvent}
 R_c^{\rm sp}=\mathcal D^*\mathcal B^{-2}\mathcal D .
\end{equation}
Indeed multiplication by $\mathcal D^*\mathcal D$ gives the projection onto $p_c^\perp$, and the right side vanishes on $p_c$.

Conjugation by $e^{sy}$ replaces $\partial_y$ by $\partial_y-s$. The Fourier symbol of the conjugated $\mathcal B$ is
$\xi^2+2is\xi+\kappa^2-s^2$. Since $|s|=\beta<\sqrt{c_-}/4$, its real part is bounded below by $\xi^2+c_-/4-\beta^2>0$. Its inverse maps $H^m$ boundedly to $H^{m+2}$, uniformly in $c$. The conjugated first-order factors have uniformly bounded smooth coefficients. Reading \eqref{eq:explicit-reduced-resolvent} as
$L^2\to H^{-1}\to H^3\to H^2$ proves the first bound. Finally
$\mathsf P_{Q_c}f=p_c\ip{p_c}{f}$, and $e^{\pm\beta y}p_c\in H^2$ uniformly. This proves the projection bound and the conclusion.
\end{proof}

\subsection{Localized eigenfunction and trace estimates}
Continue to write $w=Q_c+v$ with $\ip{v}{Q_c}=0$ and
$\norm{v}_{H^1}+r_\beta(v)<\delta$. Define
\begin{equation}\label{eq:spectral-noise-drift-full}
 \begin{split}
 F_w h&=w\psi_w\phi h,\\
 g(w)h&=\frac23\ip{w\psi_w^2}{\phi h},\\
 b(w)&=\frac19\operatorname{Tr}(F_w^*R_wF_w).
 \end{split}
\end{equation}
These are, respectively, a Hilbert--Schmidt operator $U\to L^2$, a row operator on $U$, and the spectral drift without its $\eps^2$ prefactor.

\begin{proposition}[Dependence on radiation near the soliton]\label{prop:localized-spectral-coefficients}
Choose the eigenfunction sign so that $\ip{\psi_w}{p_c}>0$. There are uniform constants such that
\begin{equation}\label{eq:eigenfunction-localization}
 \norm{\psi_w-p_c}_{H^2}\le C\norm{vp_c}_2\le Cr_\beta(v).
\end{equation}
Moreover,
\begin{equation}\label{eq:spectral-local-coefficients}
 \norm{g(w)-g(Q_c)}_U+|b(w)-b_{\rm sp}(c)|
       \le Cr_\beta(v).
\end{equation}
The constants require smallness of the global $H^1$ norm to preserve the spectral gap, but the right sides contain only the weighted radiation norm.
\end{proposition}
\begin{proof}
Write $\psi_w=\alpha p_c+z$, with $z\perp p_c$ and $\alpha>0$. The complementary eigenvalue equation is
\begin{equation}\label{eq:eigenfunction-schur}
 z=\frac{\alpha}{6}
  [\Pi_c(\Hs(Q_c)-v/6-\lambda_w)\Pi_c]^{-1}\Pi_c(vp_c).
\end{equation}
The inverse is uniformly bounded $L^2(p_c^\perp)\to H^2$ for $\delta$ sufficiently small. Its $L^2$ bound follows from the gap; the $H^2$ bound follows from the differential equation and boundedness of the potential. Thus $\norm{z}_{H^2}\le C\norm{vp_c}_2$. Normalization gives
$|\alpha-1|\le C\norm{z}_2^2$, proving \eqref{eq:eigenfunction-localization}. Proposition~\ref{prop:spectral} also gives
\begin{equation}\label{eq:eigenvalue-projection-local}
 |\lambda_w+c/4|\le Cr^2,\qquad
 \norm{\mathsf P_w-\mathsf P_{Q_c}}_{2\to2}\le Cr,\qquad r=r_\beta(v).
\end{equation}

We next compare the factors that enter the noise and trace. The product identities
\[
 w\psi_w-Q_cp_c=vp_c+w(\psi_w-p_c),\qquad
 w\psi_w^2-Q_cp_c^2=vp_c^2+w(\psi_w-p_c)(\psi_w+p_c)
\]
and the uniform $H^2$ eigenfunction bounds imply
\begin{equation}\label{eq:spectral-product-local}
 \norm{w\psi_w-Q_cp_c}_2+\norm{w\psi_w^2-Q_cp_c^2}_2\le Cr .
\end{equation}
Here $\norm{vp_c}_2\le Cr$, whereas multiplication of the eigenfunction difference by $w$ uses only the bounded global $H^1$ tube. Since $\phi$ is bounded on $L^2$, the second product estimate proves the claim for $g$. Parseval's identity applied to the first one gives
\begin{equation}\label{eq:F-local}
 \norm{F_w-F_{Q_c}}_{\HS(U,L^2)}\le Cr,\qquad
 \norm{F_w}_{\HS(U,L^2)}+\norm{F_{Q_c}}_{\HS(U,L^2)}\le C .
\end{equation}

It remains to localize the resolvent difference in the trace. Put
$F_0=F_{Q_c}$, $\mathsf S_0=\mathsf S_{Q_c}$,
$\Delta\lambda=\lambda_w+c/4$, and
$\Delta\mathsf P=\mathsf P_w-\mathsf P_{Q_c}$. The resolvent identity is
\begin{equation}\label{eq:sandwich-resolvent}
 (R_w-R_c^{\rm sp})F_0
 =\mathsf S_w^{-1}
       (v/6+\Delta\lambda-\Delta\mathsf P)\mathsf S_0^{-1}F_0
       -\Delta\mathsf P F_0 .
\end{equation}
The inverse $\mathsf S_w^{-1}$ is uniformly bounded on unweighted $L^2$. For the only term that needs localization, Lemma~\ref{lem:weighted-resolvent} and Parseval give
\[
 \norm{e^{-\beta y}\mathsf S_0^{-1}F_0}_{\HS(U,L^2)}
 \le C\norm{e^{-\beta y}F_0}_{\HS(U,L^2)}
 =C\norm{k}_2\norm{e^{-\beta y}Q_cp_c}_2\le C .
\]
Since $e^{\beta y}v\in H^1\hookrightarrow L^\infty$,
\[
 \norm{v\mathsf S_0^{-1}F_0}_{\HS(U,L^2)}
 \le\norm{e^{\beta y}v}_\infty
       \norm{e^{-\beta y}\mathsf S_0^{-1}F_0}_{\HS(U,L^2)}
 \le Cr .
\]
Together with \eqref{eq:eigenvalue-projection-local}, equation \eqref{eq:sandwich-resolvent} yields
\begin{equation}\label{eq:resolvent-local-HS}
 \norm{(R_w-R_c^{\rm sp})F_0}_{\HS(U,L^2)}\le Cr .
\end{equation}
This is a bound after applying the resolvent difference to localized noise directions, rather than a global operator-norm bound.

Finally split the trace difference into
\[
 \operatorname{Tr}(F_w^*R_wF_w-F_0^*R_wF_0)
       +\operatorname{Tr}\big(F_0^*(R_w-R_c^{\rm sp})F_0\big).
\]
The first term is bounded by
$\norm{R_w}(\norm{F_w}_{\HS}+\norm{F_0}_{\HS})\norm{F_w-F_0}_{\HS}$.
The second is bounded by
$\norm{F_0}_{\HS}\norm{(R_w-R_c^{\rm sp})F_0}_{\HS}$.
Equations \eqref{eq:F-local} and \eqref{eq:resolvent-local-HS} bound both by $Cr$. Since $b(Q_c)=b_{\rm sp}(c)$, this proves the drift estimate.
\end{proof}

\begin{proposition}[Pure spectral and geometric speed drifts]\label{prop:pure-drift-comparison}
Let $\Gamma_c(c,0)$ denote the geometric-speed component of the chart
trace in \eqref{eq:chart-gamma-full} at the pure soliton. For
$h_j=Q_c\phi e_j$, put
\[
 \alpha_j=D\Ac(Q_c)[h_j],\qquad
 \kappa_j=D\mathfrak c(Q_c)[h_j],\qquad
 w_j=h_j+\alpha_j Q'_c-\kappa_j\dot Q_c .
\]
Under Assumption~\ref{ass:noise}, $b_{\rm sp}(c)>0$, and
\begin{equation}\label{eq:pure-drift-comparison}
 b_{\rm sp}(c)-\Gamma_c(c,0)
   =\frac19\sum_j\ip{p_c w_j}{R_c^{\rm sp}(p_c w_j)}>0 .
\end{equation}
Both inequalities are uniform on compact positive speed intervals.
\end{proposition}
\begin{proof}
If $b_{\rm sp}(c)=0$, positivity of $R_c^{\rm sp}$ and
$\ker R_c^{\rm sp}=\operatorname{span}\{p_c\}$ would place the range
of $F_{Q_c}$ in $\operatorname{span}\{p_c\}$. As $p_c>0$, this
would make $Q_c\phi h$ constant for every $h\in U$. A constant
in $L^2(\R)$ vanishes. Since $Q_c>0$, this contradicts $k\ne0$.

Fix $j$ and apply the chart to $u_t=Q_c+t h_j$, with speed $c_t$
and center $a_t$. Its radiation $v_t=T_{a_t}u_t-Q_{c_t}$ satisfies
$v_t=t w_j+o(t)$ in $L^2$ by \eqref{eq:chart-derivatives} and
strong continuity of translations. We also have
$\norm{v_t}_{H^1}=O(|t|)$ by smoothness of the soliton family.
Translation invariance and \eqref{eq:localized-spectral} at $c_t$
therefore give
\[
 \lim_{t\to0}\frac{\Cs(u_t)-\mathfrak c(u_t)}{t^2}
     =\frac19\ip{p_c w_j}{R_c^{\rm sp}(p_c w_j)} .
\]
Here $p_{c_t}\to p_c$ in $L^\infty$ and $R_{c_t}^{\rm sp}\to
R_c^{\rm sp}$ in operator norm by the uniform spectral gap.
The two speed functionals have the same value and differential at
$Q_c$, so the left side equals half their Hessian difference applied
to $(h_j,h_j)$. The sums converge absolutely since $Q_c\phi$ is
Hilbert--Schmidt into $\Xb$ and the map $h_j\mapsto w_j$ is bounded
into $L^2$. Taking the trace and using \eqref{eq:full-trace-identification}
proves the identity.

Equality would make each $p_cw_j$ a scalar multiple of $p_c$, and
hence each $w_j$ constant and zero in $L^2$. The range of $Q_c\phi$
would lie in $\operatorname{span}\{Q'_c,\dot Q_c\}$. Convolution by
$k$ has infinite rank since $\widehat k$ is nonzero on an open
interval, and multiplication by $Q_c>0$ is injective. This proves
strictness. Smooth dependence of the pure coefficients on $c$ gives
the uniform assertions.
\end{proof}

In the noise-only specialization of Westdorp and Hupkes
\cite[(1.7), Sections~3 and 9]{WestdorpHupkes2026}, evaluating the
geometric-speed equation on the soliton gives precisely
$\Gamma_c(c,0)$ after profile normalization. Formula~\eqref{eq:pure-drift-comparison}
identifies the additional instantaneous drift of the spectral speed
as a sum over the directions of radiation. Over long intervals, the
geometric speed also accumulates the deterministic quadratic term
$V_c(c,v)$ in \eqref{eq:exact-geometric-ito}. The pathwise estimate
\eqref{eq:spectral-geometric} accounts for agreement of the two
leading mean asymptotics in Proposition~\ref{prop:means}.

\subsection{The exact two-parameter system and the adapted Brownian motion}
Assemble the centered phase and spectral-speed rows:
\begin{equation}\label{eq:joint-full-noise}
 \mathcal K(c,v)=\begin{pmatrix}N_a(c,v)\\g(Q_c+v)\end{pmatrix},
 \qquad q(c,v)=\mathcal K(c,v)\mathcal K(c,v)^* .
\end{equation}
The physical noise map is $\mathcal K(c,v)T_a$. At $v=0$ this is $K_0(c)T_a$, because the first variations of spectral and geometric speed agree.

For any positive definite matrix
$q=\left(\begin{smallmatrix}q_{aa}&q_{aC}\\q_{aC}&q_{CC}\end{smallmatrix}\right)$,
use the factor
\begin{equation}\label{eq:full-factor}
 \mathscr S(q)=
 \begin{pmatrix}
 q_{aC}/\sqrt{q_{CC}}&
       \sqrt{q_{aa}-q_{aC}^2/q_{CC}}\\
 \sqrt{q_{CC}}&0
 \end{pmatrix}.
\end{equation}
This is a smooth function of $q$ on the positive definite cone.

\begin{theorem}[Stopped coefficient structure]\label{thm:coefficient-closure}
Under Assumption~\ref{ass:noise}, choose the tube radius $\delta$ sufficiently small, independently of $\eps$ and time. On the continuous weighted-chart lifetime, before exit from $c,C\in I_0$ and $\norm{v}_{H^1}+r_\beta(v)<\delta$, there is a two-dimensional Brownian motion $B$ such that the exact parameters $(a,C)$ satisfy \eqref{eq:exact-structure} with
\begin{equation}\label{eq:proved-coefficient-bounds}
 |\rho_a|\le Lr_\beta(v)^2,\qquad
 |\eta_C|+|\eta_a|+|\rho_C|+|\widetilde\rho_a|
       \le Lr_\beta(v).
\end{equation}
The pure coefficients are precisely \eqref{eq:bsp}, \eqref{eq:gamma}, and the rows of \eqref{eq:sigma-factor}. The constant $L$ depends on the fixed tube, speed interval, weight, and noise kernel, but not on $\eps$, elapsed time, or the phase $a$. Thus the coefficient assertions in Assumption~\ref{ass:closure} are established; its radiation-integral assertion is not part of this theorem.
\end{theorem}
\begin{proof}
Propositions~\ref{prop:geometric-ito} and \ref{prop:h1-spectral-ito} give the exact system
\begin{equation}\label{eq:exact-before-factor}
 \begin{split}
 \dd a&=(c+V_a)\dd t+\eps N_aT_a\dd\mathcal W
                          +\eps^2\Gamma_a\dd t,\\
 \dd C&=\eps g(Q_c+v)T_a\dd\mathcal W
                          +\eps^2 b(Q_c+v)\dd t .
 \end{split}
\end{equation}
The coefficient estimates proved above imply
\[
 \norm{\mathcal K(c,v)-K_0(c)}_{\HS}\le Cr,\qquad
 |q(c,v)-q_0(c)|\le Cr .
\]
By Proposition~\ref{prop:cov}, $q_0(c)$ has a uniform positive lower eigenvalue bound. Decrease $\delta$ so that $q(c,v)$ remains in a compact set of positive definite matrices. Since $|C-c|\le Cr^2$ and the pure coefficients are smooth on $I$, the smoothness of \eqref{eq:full-factor} gives
\begin{equation}\label{eq:factor-local-bound}
 |\mathscr S(q(c,v))-\Sigma_0(C)|\le Cr .
\end{equation}

For completeness, construct $B$ without making an independence assumption about its relation to the PDE. Up to the PDE chart stop $\tau_0$, set
\[
 H(t)=\mathscr S(q(c,v))^{-1}\mathcal K(c,v)T_a .
\]
It satisfies $H(t)H(t)^*=\mathrm{Id}_2$. Choose any fixed row isometry $J:U\to\R^2$ with $JJ^*=\mathrm{Id}_2$ and put
\begin{equation}\label{eq:global-adapted-B}
 B_t=\int_0^t
       [\mathbf1_{\{s\le\tau_0\}}H(s)
          +\mathbf1_{\{s>\tau_0\}}J]\dd\mathcal W_s .
\end{equation}
This continuous local martingale has quadratic covariation $t\mathrm{Id}_2$ and is therefore a standard two-dimensional Brownian motion in the given filtration, by L\'evy's characterization \cite{Protter2005}. Stochastic integral associativity gives
$\mathcal K(c,v)T_a\dd\mathcal W=\mathscr S(q(c,v))\dd B$ before $\tau_0$. First construct this Brownian motion from the PDE, then solve the effective SDE with $B$, and finally take the common stop with the effective speed. This order avoids any circular definition of the stopping time.

The remainders are explicitly
\begin{equation}\label{eq:explicit-remainders}
 \begin{split}
 \rho_a&=c-C+V_a(c,v),\\
 \begin{pmatrix}\eta_a\\\eta_C\end{pmatrix}
     &=\mathscr S(q(c,v))-\Sigma_0(C),\\
 \rho_C&=b(Q_c+v)-b_{\rm sp}(C),\\
 \widetilde\rho_a&=\Gamma_a(c,v)-\gamma_a(C).
 \end{split}
\end{equation}
The quadratic bound follows from Proposition~\ref{prop:spectral} and \eqref{eq:geometric-coefficient-bounds}. Equation \eqref{eq:factor-local-bound} treats the noise rows. For the drifts, insert their pure values at $c$, apply Propositions~\ref{prop:geometric-ito} and \ref{prop:localized-spectral-coefficients}, and then use smoothness in $c$ and $|C-c|\le Cr^2$. This proves \eqref{eq:proved-coefficient-bounds}. All coefficients are continuous adapted functions of the stopped coordinates, hence predictable; on bounded stopped intervals their boundedness supplies the required integrability.
\end{proof}

\begin{corollary}[Criterion for stopped comparison]\label{cor:reduction-only-radiation}
On the chart lifetime, suppose
\[
 \left\|\int_0^{T\wedge\tau}r_\beta(v(t))^2\dd t\right\|_{L^2(\Omega)}
       \le K\eps^2T,\qquad 1\le T\le S_0\eps^{-2}.
\]
Then the independently specified effective SDE, coupled using \eqref{eq:global-adapted-B}, satisfies the stopped speed and position error estimates of Theorem~\ref{thm:conditional}. No additional coefficient hypothesis is required.
\end{corollary}
\begin{proof}
Theorem~\ref{thm:coefficient-closure} proves the coefficient part of Assumption~\ref{ass:closure}. The displayed estimate is its remaining input. Apply Theorem~\ref{thm:conditional}.
\end{proof}

The exact coefficients retain the radiation remainders in
\eqref{eq:explicit-remainders}. The adapted Brownian motion couples
paths for the comparison. To sample the effective law one needs only
the pure coefficients and a standard two-dimensional Brownian motion.
\section{Weighted regularity and a coercive radiation energy}\label{sec:weighted-energy}

The one-sided weight serves two purposes. A finite-time propagation estimate gives the regularity required to differentiate the modulation chart and the radiation energy. Exponential stability of the projected linearization then yields energy estimates uniform in the length of the stopped interval.

Choose $I_0=[c_0-\rho,c_0+\rho]\Subset I$ with $\rho>0$ sufficiently small. We next choose the tube radius $\delta$ and finally $\eps_0$, all independently of elapsed time. Let $\tau_{\rm P}$ denote the PDE-only exit, initially stopped also at the maximal chart lifetime. Every radiation and orbital estimate below holds up to $\tau_{\rm P}$ and up to any earlier stopping time. In particular it applies at the common stop $\tau$ of \eqref{eq:stopping-regime}. Section~\ref{sec:orbital-exit} proves continuation at every point strictly inside the PDE tube.

\subsection{Finite-time propagation of the one-sided weight}
\begin{lemma}[Weighted regularity of the actual solution]\label{lem:weighted-propagation}
Let $u$ be the global continuous $H^1$ solution of \eqref{eq:skdv}, with $u(0)\in H^1\cap H^1_\beta$. Then, almost surely, for every finite $T$,
\begin{equation}\label{eq:weighted-regularity}
 u\in C([0,T];H^1_\beta)\cap L^2(0,T;H^2_\beta).
\end{equation}
Consequently the radiation in any continuous finite moving chart has the same regularity locally in its lifetime. The constants in this preliminary assertion may depend on $T$ and on a localization of the unweighted $H^1$ norm.
\end{lemma}
\begin{proof}
We give the regularity argument because an unweighted $H^1$ existence theorem alone does not justify applying an $H^1_\beta$ energy identity. Write $D=\partial_x$ and $\kappa_j=\norm{k^{(j)}}_2^2$. First use bounded weights
\[
 m_R(x)=\frac{e^{2\beta x}}{1+e^{2\beta(x-R)}},\qquad R\ge1.
\]
They satisfy $0\le m_R'\le2\beta m_R$, $|m_R'''|\le C_\beta m_R$, and, as $R\to\infty$, $m_R\uparrow e^{2\beta x}$ and $m_R'\uparrow2\beta e^{2\beta x}$. Spatial smoothing of the actual solution, followed by passage to the limit, gives
\begin{align}
 \dd\int m_Ru^2
 ={}&\left[-3\int m_R'u_x^2+\int m_R'''u^2
       +\frac23\int m_R'u^3+\eps^2\kappa_0\int m_Ru^2\right]\dd t\notag\\
 &+2\eps\ip{\phi^*(m_Ru^2)}{\dd\mathcal W_t}_U .\label{eq:capped-weight}
\end{align}
For fixed $R$, this passage only uses convergence in $H^1$: the nonlinear smoothing error tends to zero in $H^1$, since $H^1(\R)$ is an algebra. Derivatives in the cubic drift can be moved onto $m_R$ and the smoothed solution. No third derivative of the limiting solution is required.

Stop where $\norm{u}_{H^1}>N$. Young's convolution inequality gives
\[
 \norm{\phi^*(m_Ru^2)}_U\le\norm{k}_2\int m_Ru^2.
\]
The BDG inequality, followed by Young's and Gronwall's inequalities, applied to \eqref{eq:capped-weight}, bounds moments of
\[
 \sup_{t\le T}\int m_Ru(t)^2
       +\int_0^T\int m_R'u_x^2\dd t
\]
uniformly in $R$, on this localized interval. Monotone convergence and Fatou's lemma show that $z=e^{\beta x}u$ belongs to $L^\infty(0,T;L^2)\cap L^2(0,T;H^1)$ there. It solves, in distributions,
\begin{equation}\label{eq:weighted-linear-equation}
 \dd z=\left[-(D-\beta)^3z-\frac12(D-\beta)(u z)\right]\dd t
           +\eps z\phi\dd\mathcal W_t.
\end{equation}
Localize the adapted coefficient $u$ at $\norm{u}_{H^1}=N$ and extend it by its stopped value. Set
\[
 U_\beta(t)=\exp\{t(-D^3-3\beta^2D)\},\qquad
 N_u z=-\tfrac12(D-\beta)(uz),\qquad Y(t)=U_\beta(-t)z(t).
\]
The group is unitary on every $H^s$. For $z=U_\beta(t)Y$, the conjugated equation is
\[
 \dd Y=\big[3\beta D^2Y+\beta^3Y+U_\beta(-t)N_u z\big]\dd t
       +\eps\,U_\beta(-t)G(z)\dd\mathcal W_t .
\]
This identity follows from \eqref{eq:weighted-linear-equation} by testing against $U_\beta(t)h$ for Schwartz functions $h$. One may first multiply these tests by smooth spatial cutoffs and use the weak energy bounds to pass to the limit. On compact time intervals the group preserves the Schwartz class with uniformly bounded seminorms. Thus the calculation applies to the weak weighted solution obtained above.

In \eqref{eq:variational-coercivity} below, $\norm{N_uz}_2\le C_N\norm z_{H^1}$ and
$\norm{G(z)}_{\HS(U,H^1)}^2=\kappa_0\norm z_{H^1}^2+\kappa_1\norm z_2^2$
by \eqref{eq:HS-products}. For $z\in H^2$ and $\eps\le1$ we have
\begin{equation}\label{eq:variational-coercivity}
 2\ip{(1-D^2)z}{3\beta z''+\beta^3z+N_u z}
       +\eps^2\norm{G(z)}_{\HS(U,H^1)}^2
 \le -\kappa\norm{z}_{H^2}^2+C_{N,\beta}\norm{z}_{H^1}^2.
\end{equation}
Indeed the heat operator contributes $-6\beta(\norm{z'}_2^2+\norm{z''}_2^2)$, which absorbs $C_N\norm{z}_{H^2}\norm{z}_{H^1}$ by Young's inequality. Since $U_\beta$ commutes with $D$, the drift and noise in the conjugated equation satisfy the same estimate. With pivot $H^1$ and form domain $H^2$, the drift belongs to the dual of $H^2$, paired through $1-D^2$. The linear growth bounds and \eqref{eq:variational-coercivity} give a unique variational solution in $C(H^1)\cap L^2(H^2)$; see \cite[Chapter 4]{PrevotRockner2007}. Uniqueness in the weaker class $L^\infty(L^2)\cap L^2(H^1)$ follows by the same variational argument with pivot $L^2$ and form domain $H^1$: in this triple $\|N_uz\|_{H^{-1}}\le C_N\|z\|_2$ and $\|G(z)\|_{\HS(U,L^2)}^2=\kappa_0\|z\|_2^2$. Consequently the more regular solution equals the weak weighted solution up to the localization. Removing it gives \eqref{eq:weighted-regularity}. Translation by a continuous finite $a(t)$ preserves these conclusions on every compact interval of the chart lifetime.
\end{proof}

\subsection{A fixed projection onto the stable subspace}
Put
\[
 \mathbb A_c=e^{\beta y}A_ce^{-\beta y}.
\]
The bounded projection onto the weighted constraint space is
\begin{equation}\label{eq:weighted-projection}
 \mathscr P_c h=h+
 \frac{\ip{h}{e^{-\beta y}P_c}}{D_c}e^{\beta y}Q_c'
 -\frac{\ip{h}{e^{-\beta y}Q_c}}{D_c}e^{\beta y}\dot Q_c.
\end{equation}
The dual identities in Lemma~\ref{lem:dual} show that $\mathscr P_c^2=\mathscr P_c$, that it commutes with $\mathbb A_c$, and that its range consists precisely of the two weighted orthogonality conditions. Its complement has smooth, exponentially localized range. These projections and their adjoints are bounded on every fixed Sobolev space, uniformly for $c\in I$.

Fix
\[
 P=\mathscr P_{c_0},\qquad A=\mathbb A_{c_0},\qquad
 z=e^{\beta y}v,\qquad f=Pz.
\]
Since $\mathscr P_cz=z$, the finite-rank correction $\zeta=z-f$ satisfies, for every fixed nonnegative integer $m$,
\begin{equation}\label{eq:projection-equivalence}
 \norm{\zeta}_{H^m}\le C_m|c-c_0|\norm{z}_2,\qquad
 \norm{f}_{H^j}\asymp\norm{z}_{H^j}\quad (j=0,1,2)
\end{equation}
if $\rho$ is small enough. Indeed, $\zeta=(\mathscr P_c-P)z$, and the explicit finite-rank formula is Lipschitz from $L^2$ to $H^m$. For the second statement, absorb the first bound into $z=f+\zeta$. The use of a fixed projection avoids stochastic differentiation of a moving projection.

The weighted semigroup estimates of Pego and Weinstein \cite{PegoWeinstein1994} and Mizumachi and Tzvetkov \cite{MizumachiTzvetkov2014}, as stated in \cite[Theorem 2.2]{WestdorpHupkes2026}, give
\begin{equation}\label{eq:stable-semigroup}
 T(t)=e^{tA}P,\qquad \norm{T(t)}_{L^2\to L^2}\le Me^{-\nu t},
 \qquad \norm{T(t)}_{L^2\to H^1}\le C(1+t^{-1/2})e^{-\nu t}
\end{equation}
for a fixed $\nu>0$. Our choice \eqref{eq:beta-choice} lies within that theory's weight range. We will also use
\begin{equation}\label{eq:adjoint-smoothing}
 \norm{T(t)^*}_{L^2\to H^1}\le C(1+t^{-1/2})e^{-\nu t}.
\end{equation}
Here the small-time bound follows directly for $e^{tA^*}$: its constant-coefficient part is the sum of an Airy group, $3\beta D^2$, and lower constant-coefficient terms, whose Fourier multiplier has this smoothing bound. The remaining operator has order one with bounded smooth coefficients. Duhamel's formula and the integrable kernel $(t-s)^{-1/2}$ give the same bound by the Volterra inequality. Multiplication by $P^*$ preserves $H^1$. For $t\ge1$, factor $T(t)^*=T(1)^*T(t-1)^*$ and use the exponential $L^2$ bound.

\begin{lemma}[Coercive stable energy and transport cancellation]\label{lem:stable-energy}
There is a bounded, nonnegative, self-adjoint operator $\mathscr R:L^2\to L^2$, also bounded from $L^2$ to $H^1$, and a fixed $\alpha>0$ such that
\begin{equation}\label{eq:stable-energy}
 \mathscr E(f)=\norm{f}_{H^1}^2+\alpha\ip{f}{\mathscr Rf}
       \asymp\norm{f}_{H^1}^2
\end{equation}
on $\operatorname{Ran}P$, and the linear energy form satisfies
\begin{equation}\label{eq:linear-dissipation}
 2\ip{(1-D^2)f+\alpha\mathscr Rf}{Af}\le-\kappa\norm{f}_{H^2}^2.
\end{equation}
The left side is initially defined on smooth functions in $\operatorname{Ran}P$ and extends as a quadratic form to $H^2\cap\operatorname{Ran}P$. Moreover, putting
\[
 \Theta_f=P^*[(1-D^2)f+\alpha\mathscr Rf],
\]
one has
\begin{equation}\label{eq:energy-gradient}
 \Theta_f=(1-D^2)f+s_f,\qquad \norm{s_f}_{H^1}\le C\norm{f}_2,
\end{equation}
and, for $z=f+\zeta$ as above,
\begin{equation}\label{eq:transport-cancellation}
 \left|\ip{\Theta_f}{(D-\beta)z}\right|
       \le C\norm{f}_{H^1}^2.
\end{equation}
Pairing $\Theta_f$ with $e^{\beta y}Q_c'$, $e^{\beta y}\dot Q_c$, or any fixed number of their parameter derivatives is bounded by $C\norm{f}_2$.
\end{lemma}
\begin{proof}
Define the operator-norm convergent integral
\begin{equation}\label{eq:gramian}
 \mathscr R=\int_0^\infty T(t)^*T(t)\dd t.
\end{equation}
Equations~\eqref{eq:stable-semigroup}--\eqref{eq:adjoint-smoothing} imply its asserted bounds, since $t^{-1/2}$ is integrable at zero. For $Pf=f$, differentiation under the integral, first on the generator domain, gives the Lyapunov identity
\begin{equation}\label{eq:lyapunov-identity}
 2\ip{\mathscr Rf}{Af}=-\norm{f}_2^2.
\end{equation}
The endpoint at infinity vanishes by exponential stability. This argument does not require $P$ to be orthogonal.

In coordinates,
\[
 A=-D^3+3\beta D^2+(c_0-3\beta^2)D
       +(\beta^3-c_0\beta)-D(Q_{c_0}\,\cdot)+\beta Q_{c_0}.
\]
Integration by parts and boundedness of the profile derivatives give
\[
 2\ip{(1-D^2)f}{Af}
 \le-6\beta(\norm{f'}_2^2+\norm{f''}_2^2)+C\norm{f}_{H^1}^2.
\]
For example, the odd constant-coefficient derivatives have zero quadratic form; the potential contributes only terms bounded by $C\norm{f}_{H^1}^2$. Interpolate the positive $C\norm{f'}_2^2$ against $\norm{f''}_2^2$ and $\norm{f}_2^2$. Choosing $\alpha$ large in \eqref{eq:lyapunov-identity} then proves \eqref{eq:linear-dissipation}. Positivity and boundedness of $\mathscr R$ prove the norm equivalence.

The operator $P^*-1$ has smooth finite-dimensional range. In
$(P^*-1)(1-D^2)f$, derivatives can be moved onto its smooth testing functions, giving an $H^1$ bound by $C\norm{f}_2$. The term $P^*\mathscr Rf$ has the same bound by the smoothing property of $\mathscr R$. This proves \eqref{eq:energy-gradient}. Finally,
\[
 \ip{(1-D^2)f}{(D-\beta)f}=-\beta\norm{f}_{H^1}^2.
\]
Terms containing $\zeta$ are bounded by $C\norm{f}_2\norm{\zeta}_{H^3}$ after integration by parts, and the $s_f$ term by $C\norm{f}_2\norm{z}_2$. Use \eqref{eq:projection-equivalence}. Pairings with the stated smooth profiles are handled by moving the two derivatives from $f$ to the profile.
\end{proof}

\subsection{The stochastic energy inequality}
Write
\[
 r_j(t)=\norm{e^{\beta y}v(t)}_{H^j},\qquad j=1,2,
 \qquad \mathcal E_t=\mathscr E(Pe^{\beta y}v(t)).
\]
Thus $r_1=r_\beta(v)$ and $\mathcal E\asymp r_1^2$.

\begin{proposition}[Dissipation and quadratic variation]\label{prop:weighted-energy}
For sufficiently small fixed $\rho,\delta,\eps_0$, there are constants $a_*,b_*,C_*>0$, independent of $T$ and $\eps\in(0,\eps_0)$, and a continuous real local martingale $M$ with $M_0=0$, such that, before $\tau$,
\begin{equation}\label{eq:energy-inequality}
 \dd\mathcal E_t\le
       (-a_*\mathcal E_t-b_*r_2(t)^2+C_*\eps^2)\dd t+\dd M_t,
 \qquad
 \dd[M]_t\le C_*\eps^2\mathcal E_t\dd t.
\end{equation}
The first inequality means that the difference of its right-hand side and $\dd\mathcal E$ is a nonnegative finite-variation measure. For exact soliton data, $\mathcal E_0=0$.
\end{proposition}
\begin{proof}
Use the full radiation equation \eqref{eq:full-radiation-spde}, including its cross-variation term, and denote the centered cylindrical Brownian motion there by $\widetilde{\mathcal W}$. Since $P$ commutes with $A$, the equation for $f=Pz$ is
\begin{equation}\label{eq:projected-radiation}
 \dd f=Af\dd t+P(F+\eps^2e^{\beta y}\mathcal J)\dd t
                     +\eps\mathcal H\dd\widetilde{\mathcal W},
 \qquad \mathcal H=Pe^{\beta y}\mathcal G_v,
\end{equation}
where
\[
 F=e^{\beta y}\left[(A_c-A_{c_0})v
       -\tfrac12(v^2)'+V_aw'-V_c\dot Q_c\right],\qquad w=Q_c+v.
\]
All coefficients $V,N,\Gamma,n,\chi_a$ retain the definitions from Section~\ref{sec:radiation-equation}.

We record the estimates that determine the smallness conditions. Smooth dependence on $c$ gives
$\norm{e^{\beta y}(A_c-A_{c_0})v}_2\le C\rho r_1$.
The nonlinear term is bounded by $C\norm{v}_\infty r_1\le C\delta r_1$.
The quadratic modulation bound $|V_a|+|V_c|\le Cr_1^2$ gives the same estimate for its two terms. Hence
\begin{equation}\label{eq:weighted-forcing-bounds}
 \norm{F}_2\le C(\rho+\delta)r_1,\qquad
 \norm{e^{\beta y}\mathcal J}_2\le C(1+r_2),\qquad
 \norm{\mathcal H}_{\HS(U,H^1)}^2\le C(1+r_2^2).
\end{equation}
For the second bound, the only second derivative of the radiation is the displayed $n_{aa}w''/2$ term in $\mathcal J$. Also
$\norm{\chi_a}_{W^{1,\infty}}\le C$ follows from $\chi_a=\phi N_a^*$ and $k,k'\in L^2$, so $\partial_y(w\chi_a)$ obeys the asserted bound. For the third bound, multiplication by $w\phi$ uses \eqref{eq:HS-weighted}; the rank-one term $w'N_a$ requires $r_2$, and the $\dot Q_cN_c$ term is localized and smooth.

The drift contribution of $F$ is bounded by
$C(\rho+\delta)r_2^2$, using $\norm{\Theta_f}_2\le Cr_2$.
The remaining drift and the quadratic trace of the noise are bounded by
\[
 C\eps^2 r_2(1+r_2)+C\eps^2(1+r_2^2).
\]
Choose $\rho,\delta$ small in \eqref{eq:linear-dissipation}, then $\eps_0$ small, and apply Young's inequality. A fixed fraction of the negative $H^2$ term remains and also controls $\mathcal E$. This gives the drift inequality in \eqref{eq:energy-inequality}.

A separate estimate is essential for the martingale. Its integrand is the row
\begin{equation}\label{eq:energy-martingale}
 2\eps\ip{\Theta_f}{e^{\beta y}\mathcal G_v\,\cdot}.
\end{equation}
The multiplicative term $e^{\beta y}w\phi$ has Hilbert--Schmidt $H^1$ norm bounded by $C(1+r_1)$, so its pairing with the energy gradient has row norm at most $Cr_1(1+r_1)$. For the translation term, write
\[
 e^{\beta y}w'=e^{\beta y}Q_c'+(D-\beta)z.
\]
Lemma~\ref{lem:stable-energy} bounds its scalar pairing by $C(r_1+r_1^2)$, and $N_a$ is uniformly bounded. The speed term $e^{\beta y}\dot Q_cN_c$ contributes at most $Cr_1$. Thus the row in \eqref{eq:energy-martingale} has norm at most $C\eps r_1$ in the small tube. The \Ito\ isometry identifies the quadratic-variation density with the squared row norm, proving the second inequality in \eqref{eq:energy-inequality}. In particular, no $H^2$ norm occurs in that quadratic variation.

We finish by justifying the energy calculation at the available regularity. Lemma~\ref{lem:weighted-propagation} gives $f\in C(H^1)\cap L^2(H^2)$ locally. Apply the resolvent regularizer $J_\eta=(1-\eta A)^{-1}$ to \eqref{eq:projected-radiation}. For small $\eta>0$, it commutes with $A$ and $P$, is uniformly bounded and converges strongly to the identity on $H^j$, $j=0,1,2$. These facts follow from the Laplace-transform formula for the $C_0$ semigroup on these spaces. Its gain of derivatives makes $J_\eta f$ an $H^1$ semimartingale to which the ordinary Hilbert-space \Ito\ formula applies.

Use the integrated-by-parts linear quadratic form and \eqref{eq:lyapunov-identity} before taking $\eta\downarrow0$. That form is continuous under $L^2(H^2)$ convergence. The remaining drift is in $L^2(L^2)$ by \eqref{eq:weighted-forcing-bounds}, so its pairing with the regularized gradient converges by strong convergence of $J_\eta$ and its $L^2$ adjoint. The noise converges in $L^2(\HS(U,H^1))$; the trace and martingale therefore converge, the latter by BDG after localization. This proves the identity underlying \eqref{eq:energy-inequality} for the actual solution. All further localizations can be removed using the estimates just proved. This procedure does not treat a stochastic convolution as a martingale or apply the $H^1$ \Ito\ formula directly to a third-derivative drift.
\end{proof}

\subsection{Integrated and maximal bounds from a scalar inequality}
The next lemma states the stochastic step explicitly. Its integral estimate and its maximum estimate are different consequences of the same energy inequality.

\begin{lemma}[A dissipative stochastic energy estimate]\label{lem:scalar-energy}
Let $X\ge0$ be continuous up to a stopping time $\tau$, with $X_0=0$, and suppose
\[
 \dd X\le(-aX-bZ+d\eps^2)\dd t+\dd M,
 \qquad \dd[M]\le d\eps^2X\dd t,
 \qquad 0\le X\le KZ,
\]
where $a,b,d,K>0$, $Z\ge0$ is integrable on localized finite intervals, and $M_0=0$. Then for $T\ge1$, $p\ge1$,
\begin{align}
 \left\|\int_0^{T\wedge\tau}Z(t)\dd t\right\|_{L^p}
       &\le C_p\eps^2T,\label{eq:scalar-integral}\\
 \Prob\left(\sup_{t\le T\wedge\tau}X_t\ge h\right)
       &\le C(1+T)\exp(-ch/\eps^2),\qquad h\ge0,\label{eq:scalar-tail}\\
 \left\|\sup_{t\le T\wedge\tau}X_t\right\|_{L^p}
       &\le C\eps^2[\log(2+T)+p].\label{eq:scalar-max}
\end{align}
The constants depend only on the displayed structural constants, and $C_p$ may also depend on $p$.
\end{lemma}
\begin{proof}
Integrating the inequality and discarding the nonnegative terminal energy gives, with $Y=\int_0^{T\wedge\tau}Z\dd t$,
\[
 bY\le d\eps^2T+\sup_{t\le T\wedge\tau}|M_t|.
\]
The bracket is at most $dK\eps^2Y$. BDG therefore gives
\[
 \norm{Y}_{L^p}\le C\eps^2T+C_p\eps\norm{Y}_{L^{p/2}}^{1/2}
       \le C\eps^2T+C_p\eps\norm{Y}_{L^p}^{1/2}.
\]
For $p<2$, the middle expression denotes the corresponding positive moment, not a Banach norm. Young's inequality proves $\norm{Y}_{L^p}\le C_p\eps^2(T+1)$, and hence \eqref{eq:scalar-integral}. Localization first makes the moments finite; the bound then removes the localization by Fatou's lemma.

For the maximum, choose a small fixed $\theta>0$ and apply the scalar \Ito\ formula to $V=\exp(\theta X/\eps^2)$. The bracket estimate implies that its drift is at most
\[
 V\left(d\theta-\frac{a\theta}{2}\frac{X}{\eps^2}\right)
       \le C_0-V.
\]
The last inequality follows by maximizing a linear polynomial times an exponential on $[0,\infty)$. Use the killed process $V_t\mathbf1_{\{t<\tau\}}$, whose jump at $\tau$ is nonpositive. Taking expectations and using the integrating factor gives
\begin{equation}\label{eq:killed-exponential}
 \E\left[V_t\mathbf1_{\{t<\tau\}}\right]\le C,
       \qquad t\ge0.
\end{equation}
This is a killed estimate, not an estimate for $V_{t\wedge\tau}$: a value attained at exit would otherwise be frozen indefinitely.

For an integer $j\ge0$, start on $\{j<\tau\}$ and stop at the first time $X\ge h$, at $\tau$, or at $j+1$, whichever comes first. Dropping the negative term in the drift bound and applying optional stopping gives
\[
 \Prob\left(j<\tau,\ \sup_{j\le t\le(j+1)\wedge\tau}X_t\ge h\right)
       \le C e^{-\theta h/\eps^2},
\]
by \eqref{eq:killed-exponential}. Bounded stopping and then localization justify the use of optional stopping; see \cite{Protter2005}. The first block uses $V_0=1$. A union bound over at most $1+\lceil T\rceil$ blocks proves \eqref{eq:scalar-tail}. Integrating that tail, split at a constant times $\eps^2\log(2+T)$, proves \eqref{eq:scalar-max}. This is where the logarithm enters the supremum estimate. It does not enter \eqref{eq:scalar-integral}.
\end{proof}

\begin{corollary}[Weighted radiation estimates]\label{cor:weighted-bounds}
For $T\ge1$ and $p\ge1$, the stopped radiation satisfies
\begin{align}
 \left\|\int_0^{T\wedge\tau}\norm{e^{\beta y}v(t)}_{H^2}^2\dd t\right\|_{L^p}
       &\le C_p\eps^2T,\label{eq:proved-weighted-integral}\\
 \left\|\sup_{t\le T\wedge\tau}r_\beta(v(t))^2\right\|_{L^p}
       &\le C\eps^2[\log(2+T)+p],\label{eq:proved-weighted-max}\\
 \Prob\left(\sup_{t\le T\wedge\tau}r_\beta(v(t))^2\ge h\right)
       &\le C(1+T)e^{-ch/\eps^2}.\label{eq:proved-weighted-tail}
\end{align}
In particular, \eqref{eq:r-integral-input} holds without a logarithmic loss, and the maximal bound in \eqref{eq:target-max-exit} holds at $T=T_\eps$.
\end{corollary}
\begin{proof}
Apply Lemma~\ref{lem:scalar-energy} with $X=\mathcal E$ and $Z=r_2^2$, using \eqref{eq:projection-equivalence} and \eqref{eq:stable-energy}. Since $r_1\le r_2$, the integral also controls $r_\beta(v)^2$. At $T=T_\eps$, $\log(2+T)\le C\Lambda_\eps$ for small $\eps$.
\end{proof}
\section{Global orbital control, continuation, and exit probabilities}\label{sec:orbital-exit}

The weighted estimates measure radiation close to and ahead of the soliton. A separate conserved-quantity argument controls radiation behind it. We use an auxiliary mass-matched soliton only for this energy argument; the actual modulation parameters remain the one-sided phase $a$, geometric speed $c$, and spectral speed $C$.

\subsection{A mass-matched orbital energy}
Define
\begin{equation}\label{eq:orbital-deficit}
 \begin{gathered}
 \mathcal M(u)=\frac12\int u^2,\qquad
 \mathcal H_0(u)=\frac12\int u_x^2-\frac16\int u^3,\\
 c_M(u)=\left(\frac{\mathcal M(u)}{12}\right)^{2/3},\qquad
 \mathscr D(u)=\mathcal H_0(u)+\frac{36}{5}c_M(u)^{5/2}.
 \end{gathered}
\end{equation}
Indeed, $\mathcal M(Q_c)=12c^{3/2}$ and $\mathcal H_0(Q_c)=-36c^{5/2}/5$. Both terms in $\mathscr D$ are conserved by deterministic KdV. The derivatives of the second term, viewed as a function of $\mathcal M$, are $c_M$ and $D_{c_M}^{-1}$.

\begin{lemma}[Coercivity in the one-sided chart]\label{lem:orbital-coercivity}
In a sufficiently small fixed tube, $\mathscr D(u)\ge0$ and
\begin{equation}\label{eq:orbital-coercivity}
 \mathscr D(u)\le C\norm{v}_{H^1}^2,\qquad
 \norm{v}_{H^1}^2\le C\left[\mathscr D(u)+r_\beta(v)^2\right].
\end{equation}
Moreover, the noise row satisfies
\begin{equation}\label{eq:orbital-noise}
 \norm{D\mathscr D(u)[G(u)\,\cdot]}_U^2\le C\mathscr D(u).
\end{equation}
All constants are uniform in translation and in $c\in I_0$.
\end{lemma}
\begin{proof}
The constraint $\ip{v}{Q_c}=0$ gives
\[
 \mathcal M(u)=\mathcal M(Q_c)+\tfrac12\norm{v}_2^2,
 \qquad |c_M-c|\le C\norm{v}_2^2.
\]
Choose an auxiliary translation $b$ so that
$\omega=u(\cdot+b)-Q_{c_M}$ is orthogonal to $Q_{c_M}'$. The usual implicit-function argument is uniform in the small tube. Mass matching also gives
$\ip{\omega}{Q_{c_M}}=-\norm{\omega}_2^2/2$.
The standard constrained coercivity of
$L_c=-\partial_y^2+c-Q_c$ on the orthogonal complement of $Q_c,Q_c'$ implies
\[
 \mathscr D(u)=\frac12\ip{L_{c_M}\omega}{\omega}
                       -\frac16\int\omega^3
       \asymp\norm{\omega}_{H^1}^2.
\]
For clarity, subtracting the $Q_{c_M}$ component of $\omega$ changes it by $O(\norm{\omega}_{H^1}^2)$; the constrained quadratic bound then absorbs this change and the cubic term. This is the usual orbital coercivity argument; see \cite{Weinstein1986,BonaSouganidisStrauss1987}.

At the phase $a$, the difference $Q_c-Q_{c_M}$ is even, and so has zero pairing with $Q_{c_M}'$. The translation implicit-function estimate therefore gives
\[
 |b-a|\le C|\ip{v}{Q_{c_M}'}|
       \le C\norm{e^{\beta y}v}_2.
\]
It follows that
$\norm{v}_{H^1}\le C\norm{\omega}_{H^1}+Cr_\beta(v)+C\norm{v}_{H^1}^2$.
Absorbing the last term proves the second bound in \eqref{eq:orbital-coercivity}; the first follows by Taylor expansion at $Q_c$, where $\mathscr D$ and its first derivative vanish. This reasoning uses an auxiliary $L^2$ phase for coercivity, rather than claiming unweighted coercivity from the non-$L^2$ test $P_c$.

Finally $D\mathscr D(Q_d)=0$ for every nearby $d$. Its derivative is locally Lipschitz from $H^1$ to $H^{-1}$, uniformly along the translated soliton family. Since $\norm{G(u)}_{\HS(U,H^1)}$ is bounded in the tube, comparison with the translated $Q_{c_M}$ gives
$\norm{D\mathscr D(u)G(u)}_U\le C\norm{\omega}_{H^1}$. This proves \eqref{eq:orbital-noise}.
\end{proof}

\begin{proposition}[Unweighted orbital moments]\label{prop:orbital-moments}
For $p\ge1$, $T\ge1$, and $S=\eps^2T$, one has
\begin{equation}\label{eq:orbital-moments}
 \left\|\sup_{t\le T\wedge\tau}\mathscr D(u(t))\right\|_{L^p}\le C_pS,
 \qquad
 \left\|\sup_{t\le T\wedge\tau}\norm{v(t)}_{H^1}^2\right\|_{L^p}\le C_pS.
\end{equation}
\end{proposition}
\begin{proof}
Apply the mass and energy balances of de Bouard and Debussche \cite[Lemma 2.3]{deBouardDebussche2009} to the $H^1$ solution, followed by the scalar \Ito\ formula for the function of mass in \eqref{eq:orbital-deficit}. This gives
\begin{equation}\label{eq:orbital-ito}
 \dd\mathscr D(u)=\eps D\mathscr D(u)[G(u)\dd\mathcal W]
                      +\eps^2 B_{\mathscr D}(u)\dd t,
\end{equation}
where in our normalization
\begin{align}
 B_{\mathscr D}(u)=\frac12\bigg[
 &\kappa_0\norm{u_x}_2^2+\kappa_1\norm{u}_2^2
       -\kappa_0\int u^3+c_M\kappa_0\norm{u}_2^2\notag\\
 &+D_{c_M}^{-1}\norm{\phi^*(u^2)}_U^2\bigg].\label{eq:orbital-trace}
\end{align}
For example, the derivative part of the energy trace uses
$\sum_j\norm{(u\phi e_j)'}_2^2=\kappa_0\norm{u'}_2^2+\kappa_1\norm{u}_2^2$.
The last term in \eqref{eq:orbital-trace} is the quadratic variation of mass under the scalar chain rule. In the tube, $c_M$ stays in a compact positive interval and $|B_{\mathscr D}(u)|\le C$.

Write $D_* = \sup_{t\le T\wedge\tau}\mathscr D(u(t))$. Nonnegativity, $\mathscr D(u(0))=0$, and \eqref{eq:orbital-noise} give, by BDG,
\[
 \norm{D_*}_{L^p}\le C\eps^2T
       +C_p\eps\left\|\left(\int_0^{T\wedge\tau}\mathscr D(u(t))\dd t\right)^{1/2}\right\|_{L^p}
 \le CS+C_p\sqrt S\norm{D_*}_{L^p}^{1/2}.
\]
Young's inequality proves the first assertion. Combine \eqref{eq:orbital-coercivity} with \eqref{eq:proved-weighted-max}. Since $\log(2+T)\le CT$ for $T\ge1$, the second assertion follows. The constants may depend on $p$, but not on $T$ or $\eps$.
\end{proof}

\subsection{Continuation and exit on the slow time scale}
\begin{theorem}[Continuation and exponential exit bounds]\label{thm:exit}
The weighted chart continues up to the PDE exit $\tau_{\rm P}$ from the conditions on $(c,C,v)$ in \eqref{eq:stopping-regime}. Let $\tau$ also include the effective-speed exit, as in that equation. There are $S_0,K,\kappa>0$ such that, for $S=\eps^2T\le S_0$, $T\ge1$, and either $\vartheta=\tau_{\rm P}$ or $\vartheta=\tau$,
\begin{equation}\label{eq:proved-exit}
 \Prob(\vartheta\le T)\le K e^{-\kappa/S}+K(1+T)e^{-\kappa/\eps^2}.
\end{equation}
Thus $\Prob(\vartheta>T)\to1$ for every deterministic family $1\le T=o(\eps^{-2})$. In particular,
\begin{equation}\label{eq:near-diffusive-exit}
 \Prob(\vartheta\le T_\eps)
       \le K e^{-\kappa\Lambda_\eps^2}
          +K(1+T_\eps)e^{-\kappa/\eps^2}=o(\eps^m)
       \qquad\text{for each fixed }m>0 .
\end{equation}
\end{theorem}
\begin{proof}
First construct the chart along the continuous $\Xb$ solution given by Lemma~\ref{lem:weighted-propagation}, using Proposition~\ref{prop:chart} on overlapping neighborhoods up to the PDE exit. At any finite endpoint of this construction the drifts of $a,c$ and the quadratic-variation densities of their martingales are bounded while the solution stays in the tube. Consequently $a,c$ have finite limits there. Strong continuity of translation in $\Xb$ gives a limit for $v=T_au-Q_c$, while continuity of the isolated eigenvalue gives a limit for $C$. If the limiting point is in the interior, the same local chart extends the construction; a point on its boundary is an explicit exit. This proves continuation for the PDE-only stop. The Brownian motion of \eqref{eq:global-adapted-B} is now defined up to that stop and extended beyond it as in Section~\ref{sec:localized-closure}. Solve the closed SDE with this driver and intersect the two lifetimes to obtain $\tau$.

Consider first $\vartheta=\tau$. In the stopped spectral and effective speed equations, the absolute time-integrated drifts are bounded by $KS$, and each martingale has bracket at most $KS$ up to $T$. In particular, BDG gives, for $p\ge1$,
\begin{equation}\label{eq:speed-exit-moments}
 \left\|\sup_{t\le T\wedge\tau}|C(t)-c_0|\right\|_{L^{2p}}
 +\left\|\sup_{t\le T\wedge\tau}|\bc(t)-c_0|\right\|_{L^{2p}}
       \le K_p(\sqrt S+S).
\end{equation}
The same drift and bracket bounds apply to the stopped deficit $\mathscr D$. Indeed, \eqref{eq:orbital-trace} bounds its drift density by $K\eps^2$, and \eqref{eq:orbital-noise} bounds its bracket density by $K\eps^2\mathscr D\le K\eps^2\delta^2$. All three processes start with zero increment. For a fixed $h>0$ choose $S_0$ so that $KS_0\le h/2$. For each of these increments $X$, the exponential maximal inequality for continuous martingales gives
\[
 \Prob\!\left(\sup_{t\le T\wedge\tau}|X_t|\ge h\right)
 \le 2\exp\!\left(-\frac{h^2}{8KS}\right).
\]
Here one applies the nonnegative exponential supermartingales of both signs to the first martingale crossing of $h/2$. This stopped-bracket argument was used for the stochastic KdV orbital exit by de Bouard and Debussche \cite{deBouardDebussche2009}.

Choose $\delta$ sufficiently small that $|C-c|\le Kr_1^2\le K\delta^2<\rho/2$. An exit through any of the speeds then requires $|C-c_0|$ or $|\bc-c_0|$ to reach $\rho/2$. By \eqref{eq:orbital-coercivity}, an exit through $\|v\|_{H^1}+r_1=\delta$ requires either $\mathscr D$ or $r_1^2$ to reach a fixed $h_\delta>0$. Apply the preceding exponential bound for the speeds and deficit, using $h=\min(\rho/2,h_\delta)$, and apply \eqref{eq:proved-weighted-tail} to $r_1^2$. The union bound proves \eqref{eq:proved-exit}. For $\tau_{\rm P}$ the same proof uses only the spectral speed, deficit, and radiation; none of their bounds involves $\bc$. The two asserted limits follow by substitution.
\end{proof}

\begin{proposition}[Exit at a fixed positive slow time]\label{prop:fixed-slow-exit}
For every $s>0$ there is a number $p_s>0$, independent of $\eps$, such that
\begin{equation}\label{eq:fixed-slow-exit}
 \Prob(\tau\le s\eps^{-2})\ \ge\ \Prob(\tau_{\rm P}\le s\eps^{-2})\ \ge\ p_s .
\end{equation}
\end{proposition}
\begin{proof}
By \eqref{eq:spectral-local-coefficients}, \eqref{eq:g}, and compactness of $I_0$, the tube can be chosen so that the squared norm of the exact spectral noise row is bounded below by $m>0$. The physical translation of this row preserves its norm. The exact drift $\eps^2 b(Q_c+v)$ is nonnegative: $R_{Q_c+v}$ is the positive reduced resolvent at the lowest eigenvalue, and \eqref{eq:spectral-noise-drift-full} expresses $b$ as a positive trace.

Write $M$ for the martingale part of $C-c_0$ until $\tau_{\rm P}$ and, after that stop, continue its stochastic integral against the original cylindrical noise with any fixed row of norm $\eps\sqrt m$. This defines a continuous local martingale with $[M]_t\ge m\eps^2t$ for all $t$. By the Dambis--Dubins--Schwarz theorem \cite{Protter2005}, $M_t=W_{[M]_t}$ for a standard Brownian motion $W$. Suppose $\tau_{\rm P}>s\eps^{-2}$. Then $C-c_0<\rho$ throughout the interval, and the nonnegative drift gives $M\le C-c_0$. Since the bracket is continuous and reaches at least $ms$, it follows that $\sup_{r\le ms}W_r<\rho$. The reflection principle therefore yields the uniform bound in \eqref{eq:fixed-slow-exit} with
\[
 p_s=2\bigl[1-\Phi(\rho/\sqrt{ms})\bigr]>0 ,
\]
where $\Phi$ is the standard normal distribution function. Finally $\tau\le\tau_{\rm P}$ by construction.
\end{proof}

Proposition~\ref{prop:fixed-slow-exit} identifies the effect of keeping a fixed speed interval. Even the PDE-only chart has a positive chance of leaving this regime by a fixed positive slow time. The proposition concerns exit from the specified tube; the solution may still lie near a soliton with a different speed.
\section{Strong comparison and proof of the main theorem}\label{sec:target}

We first isolate the finite-dimensional comparison argument and then assemble the proof of Theorem~\ref{target:main}. All estimates use the common stop \eqref{eq:stopping-regime}, with the fixed radii constructed in Sections~\ref{sec:weighted-energy}--\ref{sec:orbital-exit}.

\subsection{The exact analytic inputs needed for comparison}
The next assumption records the inputs to the finite-dimensional comparison argument as an abstract criterion. For \eqref{eq:skdv}, Theorem~\ref{thm:coefficient-closure} and Corollary~\ref{cor:weighted-bounds} establish all its hypotheses; Theorem~\ref{thm:exit} permits removal of the stop in probability on the chosen window.

\begin{assumption}[Stopped coefficient structure and radiation integral]\label{ass:closure}
On a common stopped interval, let $a,C$ be continuous semimartingales driven by a two-dimensional Brownian motion $B$, and let $r\ge0$ be adapted. Assume
\begin{equation}\label{eq:exact-structure}
 \begin{aligned}
 \dd C={}&\eps\big[\sigma_C(C)+\eta_C\big]\dd B
                 +\eps^2\big[b_{\rm sp}(C)+\rho_C\big]\dd t,\\
 \dd a={}&C\dd t+\rho_a\dd t
       +\eps\big[\sigma_a(C)+\eta_a\big]\dd B
       +\eps^2\big[\gamma_a(C)+\widetilde\rho_a\big]\dd t,
 \end{aligned}
\end{equation}
where all stochastic integrands are predictable, all drift terms are integrable, and
\begin{equation}\label{eq:coefficient-remainders}
 |\rho_a|\le Lr^2,\qquad
 |\eta_C|+|\eta_a|+|\rho_C|+|\widetilde\rho_a|\le Lr.
\end{equation}
The deterministic coefficient functions are Lipschitz on the region used, and for $1\le T\le S_0\eps^{-2}$,
\begin{equation}\label{eq:r-integral-input}
 \left\|\int_0^{T\wedge\tau}r(t)^2\dd t\right\|_{L^2(\Omega)}\le K S,
       \qquad S=\eps^2T.
\end{equation}
The constants $L,K,S_0$ are independent of $\eps,T$.
\end{assumption}

For the intended application, $r=r_\beta(v)$. Proposition~\ref{prop:det-mod} and \eqref{eq:define-C} explain the quadratic bound on $\rho_a$, while the isospectral identity explains the absence of such a deterministic remainder in the first line of \eqref{eq:exact-structure}. The stochastic trace and covariance-factor estimates needed for \eqref{eq:coefficient-remainders} are proved in Section~\ref{sec:localized-closure}.

\subsection{Conditional strong comparison}
\begin{theorem}[Comparison from the stated inputs]\label{thm:conditional}
Assume Assumption~\ref{ass:closure}. Let $(\ba,\bc)$ solve \eqref{eq:intro-effective} with the same Brownian motion and the same initial parameters. Then
\begin{equation}\label{eq:conditional-speed}
 \left\|\sup_{t\le T\wedge\tau}|C(t)-\bc(t)|\right\|_{L^2}
       \le C\eps\sqrt S.
\end{equation}
Moreover, the compensated position error satisfies
\begin{equation}\label{eq:compensated}
 \left\|\sup_{t\le T\wedge\tau}
 \left|a(t)-\ba(t)-\int_0^t(C(s)-\bc(s))\dd s\right|\right\|_{L^2}
       \le C(S+\eps\sqrt S).
\end{equation}
Consequently \eqref{eq:target-errors} holds. The constants depend on the coefficient bounds, $K$, and $S_0$, but not on $T$ or $\eps$ separately.
\end{theorem}
\begin{proof}
All processes in this proof are stopped at $\tau$. Put $\Delta C=C-\bc$ and
$F(t)=\E\sup_{s\le t\wedge\tau}|\Delta C(s)|^2$.
Subtract the speed equations. The $p=2$ BDG inequality bounds the martingale part by
\[
 C\eps^2\int_0^t F(s)\dd s
       +C\eps^2\E\int_0^{t\wedge\tau}r(s)^2\dd s.
\]
Cauchy--Schwarz in time bounds the squared drift contribution by
\[
 C\eps^4t\int_0^t F(s)\dd s
       +C\eps^4t\E\int_0^{t\wedge\tau}r(s)^2\dd s.
\]
Since $\eps^2t\le S_0$, \eqref{eq:r-integral-input} yields
\[
 F(t)\le C\eps^2\int_0^tF(s)\dd s+C\eps^2(S+S^2).
\]
Gronwall's inequality contributes at most $e^{CS_0}$, proving \eqref{eq:conditional-speed}.

For the compensated position, subtract its integrated speed difference. The quadratic drift is bounded in $L^2$ by $CKS$. The noise remainder has maximal $L^2$ norm at most $C\eps\sqrt S$ by \eqref{eq:bdg}. The stochastic coefficient difference arising from $\Delta C$ contributes at most
\[
 C\eps\sqrt T\norm{\sup_{t\le T\wedge\tau}|\Delta C(t)|}_{L^2}
       \le C\eps S.
\]
The order-$\eps^2$ trace remainder contributes at most
$C\eps^2\sqrt T\sqrt S=C\eps S$.
Finally the phase drift difference contributes at most
$C\eps^2T\norm{\sup|\Delta C|}_{L^2}\le C\eps S^{3/2}$.
As $S\le S_0$, these bounds imply \eqref{eq:compensated}.

Adding back the integrated speed error and multiplying by $\eps^2$ gives
\[
 \norm{\sup_{t\le T\wedge\tau}|\eps^2(a(t)-\ba(t))|}_{L^2}
 \le C\eps^2\big(T\eps\sqrt S+S+\eps\sqrt S\big)
 =C\big(\eps S^{3/2}+\eps^2S+\eps^3\sqrt S\big).
\]
Because $T\ge1$, we have $S\ge\eps^2$, and each of the last two terms is at most the first. This proves the asserted comparison.
\end{proof}

\subsection{Assembly of the proof}\label{sec:main-proof}
\begin{proof}[Proof of Theorem~\ref{target:main}]
The global $H^1$ solution is supplied by the well-posedness theory cited in Section~\ref{sec:stochastic}. Lemma~\ref{lem:weighted-propagation} propagates its one-sided weight on finite intervals. Proposition~\ref{prop:chart} constructs the local geometric coordinates, and Theorem~\ref{thm:exit} continues them up to the explicit tube exit. The spectral coordinate and its $H^1$ \Ito\ equation are supplied by Propositions~\ref{prop:spectral} and \ref{prop:h1-spectral-ito}.

Corollary~\ref{cor:weighted-bounds} and Proposition~\ref{prop:orbital-moments} give \eqref{eq:target-radiation} and \eqref{eq:main-weighted-max}. Theorem~\ref{thm:coefficient-closure} constructs the adapted Brownian motion and proves the exact coefficient structure in Assumption~\ref{ass:closure}. That construction first uses the PDE-only stop and then solves the extended closed SDE, after which the common stop is imposed. There is therefore no circular definition involving $\bc$.

The weighted radiation integral verifies the remaining input \eqref{eq:r-integral-input}. Theorem~\ref{thm:conditional} proves \eqref{eq:target-errors}. Theorem~\ref{thm:exit} gives \eqref{eq:main-exit}. At $T=T_\eps$, $S=\Lambda_\eps^{-2}$ and $\log(2+T_\eps)\lesssim\Lambda_\eps$, yielding \eqref{eq:target-max-exit}. Finally, combine the localized correction $|C-c|\le Cr_\beta(v)^2$ with the spectral-speed error to obtain \eqref{eq:target-transfer}.
\end{proof}

\begin{corollary}[Approximation in physical space on a shorter interval]\label{cor:physical-reconstruction}
In the adapted coupling of Theorem~\ref{target:main}, for $1\le T\le S_0\eps^{-2}$,
\begin{equation}\label{eq:physical-reconstruction}
 \left\|\sup_{t\le T\wedge\tau}
 \left\|u(t)-Q_{\bc(t)}(\,\cdot-\ba(t))\right\|_{H^1}\right\|_{L^2}
 \le K\big(\eps\sqrt T+\eps^2T^{3/2}\big).
\end{equation}
Consequently, if $1\le T=o(\eps^{-4/3})$, for each $\eta>0$ the probability that the chart survives to $T$ and that the $H^1$ error stays below $\eta$ on $[0,T]$ tends to one.
\end{corollary}
\begin{proof}
On $I_0$ the derivatives of the translated profile with respect to speed and position are bounded in $H^1$, uniformly in translation. Hence the norm in \eqref{eq:physical-reconstruction} is bounded pointwise by
$\|v\|_{H^1}+K(|c-\bc|+|a-\ba|)$.
The first summand is $O_{L^2}(\sqrt S)$ after taking the stopped supremum. For the speed, \eqref{eq:target-errors} bounds $C-\bc$ by $O_{L^2}(\eps\sqrt S)$ and \eqref{eq:spectral-geometric}, \eqref{eq:main-weighted-max} bound $c-C$ by $O_{L^2}(\eps^2\log(2+T))$. The latter is at most $K\sqrt S$ since $\eps\le1$ and $\log(2+T)\le K\sqrt T$ for $T\ge1$. Unscaling the phase estimate in \eqref{eq:target-errors} gives $O_{L^2}(\eps^2T^{3/2})$. This proves \eqref{eq:physical-reconstruction}. Both terms tend to zero on the stated shorter interval, and \eqref{eq:proved-exit} removes the stop in probability.
\end{proof}
\section{Parameter observables and stochastic asymptotics}\label{sec:interpretation}

\subsection{Comparison with the virial coordinates}\label{sec:chart-comparison}
The spectral speed depends only on the translation orbit of the solution. To compare two choices of geometric parameters, we use the localized tests $Q_c$ and $yQ_c$ for a second chart. At leading order the change of phase is linear in the radiation, and the change of speed is quadratic.

\begin{proposition}[Localized change of chart]\label{prop:virial-chart}
For $T_au=Q_c+v$ in the small one-sided tube, there is a unique nearby pair $(a_{\rm v},c_{\rm v})$ for which
\begin{equation}\label{eq:virial-chart}
 v_{\rm v}=T_{a_{\rm v}}u-Q_{c_{\rm v}},\qquad
 \ip{v_{\rm v}}{Q_{c_{\rm v}}}
       =\ip{v_{\rm v}}{yQ_{c_{\rm v}}}=0.
\end{equation}
Set $r_0=\norm{e^{\beta y}v}_2$ and
$h_1=2\ip{v}{yQ_c}/\norm{Q_c}_2^2$. Uniformly in the tube,
\begin{equation}\label{eq:chart-change-bounds}
 a_{\rm v}-a=h_1+O(r_0^2),\qquad
 |a_{\rm v}-a|\le Cr_0,\qquad |c_{\rm v}-c|\le Cr_0^2.
\end{equation}
More precisely,
\begin{equation}\label{eq:chart-speed-quadratic}
 c_{\rm v}-c=-\frac1{D_c}\left[
 h_1\ip{v}{Q_c'}+\frac12h_1^2\norm{Q_c'}_2^2\right]+O(r_0^3).
\end{equation}
Consequently, for $T\ge1$, $p\ge1$,
\begin{equation}\label{eq:chart-change-moments}
 \begin{aligned}
 \left\|\sup_{t\le T\wedge\tau}|a_{\rm v}-a|\right\|_{L^{2p}}
       &\le C\eps[\log(2+T)+p]^{1/2},\\
 \left\|\sup_{t\le T\wedge\tau}|c_{\rm v}-c|\right\|_{L^p}
       &\le C\eps^2[\log(2+T)+p].
 \end{aligned}
\end{equation}
At $T=T_\eps$, the error estimates of Theorem~\ref{target:main} hold also with $(a,c)$ replaced by $(a_{\rm v},c_{\rm v})$, at the same stated scales.
\end{proposition}
\begin{proof}
Write $h=a_{\rm v}-a$ and $d=c_{\rm v}-c$. Apply the implicit-function theorem to
\[
 F(h,d;v)=
 \begin{pmatrix}
 \ip{T_h(Q_c+v)-Q_{c+d}}{Q_{c+d}}\\
 \ip{T_h(Q_c+v)-Q_{c+d}}{yQ_{c+d}}
 \end{pmatrix}.
\]
Moving the translation onto the tests makes this a smooth function of $(h,d,v)$ even for $e^{\beta y}v\in L^2$. All its required derivatives are bounded because the translated tests and their derivatives decay exponentially. At $v=0$,
\[
 D_{(h,d)}F(0,0;0)=
 \begin{pmatrix}0&-D_c\\-\norm{Q_c}_2^2/2&0\end{pmatrix}.
\]
On the constraint $\ip{v}{Q_c}=0$, the first-order changes are therefore $h=h_1$ and $d=0$. Uniform Taylor expansion proves \eqref{eq:chart-change-bounds}. Expanding the first constraint one order further gives
\[
 0=-D_cd-h\ip{v}{Q_c'}-\tfrac12h^2\norm{Q_c'}_2^2+O(r_0^3).
\]
Substitute $h=h_1+O(r_0^2)$ to obtain \eqref{eq:chart-speed-quadratic}.

The moment estimates follow from \eqref{eq:main-weighted-max}. Combine them with $|C-c|\le Cr_0^2$ and the theorem's parameter errors. At $T=T_\eps$, the additional speed error is $O_{L^2}(\eps^2\Lambda_\eps)$, and the additional error in $\eps^2a$ is $O_{L^2}(\eps^3\sqrt{\Lambda_\eps})$. Both are smaller than the corresponding errors in \eqref{eq:target-errors} for small $\eps$.
\end{proof}

Thus the two PDE phases approach one another even in absolute distance on this interval. The change of coordinates depends on the radiation, so differentiating it also changes the instantaneous noise and drift in the phase equation.

\subsection{Speed and position fluctuations on the near-diffusive window}\label{sec:memory}
Here $T_\eps\to\infty$ in physical time and $\eps^2T_\eps=\Lambda_\eps^{-2}\to0$ in slow time. For convergence-in-law statements about PDE coordinates we may freeze the centered parameters at the tube exit, stopping the clock in the centering term as well. The probability of this exit tends to zero by Theorem~\ref{thm:exit}.

Let
\[
 g_*^2=q_{CC}(c_0)>0,\qquad
 A_0(t)=a_0+c_0t.
\]
De Bouard and Debussche corrected the phase using a weighted adjoint projection in their analysis of multiplicative KdV noise \cite[Section~4, Proposition~4.1 and (4.4)--(4.7)]{deBouardDebussche2009}. In Section~5 of the same paper they wrote the resulting first-order joint system with correlated Brownian motions in (5.1), and computed its full covariance matrix in the display following (5.2). In particular, their position variance contains the cubic, quadratic, and linear powers of time. De Bouard and Gautier studied how allowing modulation changes exit from a neighborhood of a deterministic trajectory for additive forcing \cite{deBouardGautier2010}.

For our coordinates, freezing the coefficients at $c_0$ and retaining the terms of order $\eps$ gives the Gaussian approximation
\begin{equation}\label{eq:frozen}
 \begin{aligned}
 C^{(1)}(t)-c_0&=\eps\sigma_C(c_0)B_t,\\
 a^{(1)}(t)-A_0(t)&=
       \eps\sigma_C(c_0)\int_0^tB_r\dd r
                      +\eps\sigma_a(c_0)B_t .
 \end{aligned}
\end{equation}
Its covariance has the form above with the noise rows in \eqref{eq:sigma-factor}. In particular,
\begin{equation}\label{eq:phase-variance}
 \operatorname{Var}(a^{(1)}(t)-A_0(t))=
 \eps^2\left(\frac{q_{CC}(c_0)t^3}{3}
                  +q_{aC}(c_0)t^2+q_{aa}(c_0)t\right).
\end{equation}
The following result transfers the leading integrated Brownian fluctuation to the nonlinear PDE on the growing physical interval $[0,T_\eps]$.

\begin{proposition}[Small speed fluctuations and integrated phase fluctuations]\label{prop:limit}
For the smoothly extended effective SDE, as $\eps\to0$,
\begin{equation}\label{eq:joint-limit}
 \left(
 \Lambda_\eps[\bc(T_\eps s)-c_0],\,
 \eps^2\Lambda_\eps^3[\ba(T_\eps s)-a_0-c_0T_\eps s]
 \right)_{0\le s\le1}
 \ \Longrightarrow\
 \left(g_*W_s,\,g_*\int_0^sW_r\dd r\right)_{0\le s\le1}
\end{equation}
in $C([0,1];\R^2)$, where $W$ is standard one-dimensional Brownian motion. The probability that the effective speed leaves $I_0$ before $T_\eps$ tends to zero. By Theorem~\ref{target:main}, the same limit holds with $(\ba,\bc)$ replaced by $(a,C)$, and also with $C$ replaced by $c$ or $c_{\rm v}$ and with $a$ replaced by $a_{\rm v}$.
\end{proposition}
\begin{proof}
Put $\widetilde B_s=T_\eps^{-1/2}B_{T_\eps s}$, which is standard two-dimensional Brownian motion for every $\eps$. On this time scale the speed equation is
\[
 \dd\bc(T_\eps s)=\Lambda_\eps^{-1}\sigma_C(\bc(T_\eps s))\dd\widetilde B_s
                +\Lambda_\eps^{-2}b_{\rm sp}(\bc(T_\eps s))\dd s.
\]
Bounded coefficients and BDG give
$\E\sup_{s\le1}|\bc(T_\eps s)-c_0|^2\lesssim\Lambda_\eps^{-2}$.
This implies the exit assertion. After multiplying by $\Lambda_\eps$, the difference from $\sigma_C(c_0)\widetilde B_s$ tends to zero in maximal $L^2$, by Lipschitz continuity and another use of BDG. The drift contributes $O(\Lambda_\eps^{-1})$ after this normalization.

For position, the normalized integrated speed is exactly
\[
 \eps^2\Lambda_\eps^3
       \int_0^{T_\eps s}(\bc(t)-c_0)\dd t
       =\int_0^s\Lambda_\eps(\bc(T_\eps r)-c_0)\dd r.
\]
The direct phase martingale has maximal $L^2$ size
$O(\eps^2\Lambda_\eps^2)$ in this normalization, and the phase \Ito\ drift has size $O(\eps^2\Lambda_\eps)$. Both vanish. Continuity of the time-integration map proves \eqref{eq:joint-limit}, since $\sigma_C(c_0)\widetilde B$ has the law of $g_*W$.

Theorem~\ref{target:main} gives speed error multiplied by $\Lambda_\eps$ and scaled position error multiplied by $\Lambda_\eps^3$ both of size $O_{L^2}(\eps)$ before the stop. The probability of the stop vanishes. Finally \eqref{eq:target-max-exit} and \eqref{eq:spectral-geometric} give
$\Lambda_\eps\sup|C-c|=O_{L^2}(\eps^2\Lambda_\eps^2)\to0$. These observations transfer the limit. Proposition~\ref{prop:virial-chart} gives the corresponding virial-coordinate conclusion.
\end{proof}

Thus $\bc-c_0=\Op(\Lambda_\eps^{-1})\to0$, uniformly over the interval. Also
\[
 \sup_{s\le1}\frac{|\ba(T_\eps s)-A_0(T_\eps s)|}{T_\eps}
            =\Op(\Lambda_\eps^{-1})\longrightarrow0.
\]
In these speed and relative-position senses the free dynamics is recovered. However, the absolute position fluctuation scale at the endpoint is
\begin{equation}\label{eq:absolute-phase-scale}
 \eps T_\eps^{3/2}=\eps^{-2}\Lambda_\eps^{-3}\longrightarrow\infty.
\end{equation}
The frozen covariance \eqref{eq:phase-variance} places order-one fluctuations of the absolute position at $t\asymp\eps^{-2/3}$. At the later time $T_\eps$, Proposition~\ref{prop:limit} identifies a nondegenerate Gaussian endpoint limit for the normalized position. Fixing a positive slow time $s_0$ instead would give physical time $s_0\eps^{-2}$. On that interval a fixed speed tube already has a positive exit probability by Proposition~\ref{prop:fixed-slow-exit}.

\subsection{Second-order fluctuations and the spectral drift}\label{sec:second-order}
The leading limit in Proposition~\ref{prop:limit} suppresses the spectral drift. A second expansion reveals what remains after the leading Brownian term is removed. This requires retaining the variation of the diffusion coefficient as well as its value at $c_0$.

Use the covariance factor \eqref{eq:sigma-factor}, and choose the smooth extensions so that $\sigma_C=(g,0)$ globally, with $g(c)=\sqrt{q_{CC}(c)}$ on $I$. Write
\begin{equation}\label{eq:second-order-notation}
 d_\eps=\Lambda_\eps^{-1},\quad
 W_s^\eps=T_\eps^{-1/2}B^1_{T_\eps s},\quad
 g_0=g(c_0)=g_*,\quad g'_0=g'(c_0),\quad b_0=b_{\rm sp}(c_0).
\end{equation}
The process $W^\eps$ is standard one-dimensional Brownian motion for each $\eps$; it is chosen on the probability space of the strong coupling. Define
\begin{equation}\label{eq:second-order-Z}
 Z_s^\eps=b_0s+\frac12g_0g'_0[(W_s^\eps)^2-s],\qquad
 J_s^\eps=\int_0^sZ_r^\eps\dd r.
\end{equation}
The coefficient of the centered square equals $q_{CC}'(c_0)/4>0$ by Lemma~\ref{lem:increasing-speed-variance}.

\begin{theorem}[Second-order expansion on the same time window]\label{thm:second-order}
Let $A_\eps=\{\tau>T_\eps\}$ and
$Y_\eps(s)=\eps^2[a(T_\eps s)-a_0-c_0T_\eps s]$ on $A_\eps$. Then
\begin{align}
 \left\|\mathbf1_{A_\eps}\sup_{s\le1}
 \left|C(T_\eps s)-c_0-d_\eps g_0W_s^\eps-d_\eps^2Z_s^\eps\right|\right\|_{L^2}
       &\le C(d_\eps^3+\eps d_\eps),\label{eq:second-order-speed-error}\\
 \left\|\mathbf1_{A_\eps}\sup_{s\le1}
 \left|Y_\eps(s)-d_\eps^3g_0\int_0^sW_r^\eps\dd r
                         -d_\eps^4J_s^\eps\right|\right\|_{L^2}
       &\le C(d_\eps^5+\eps d_\eps^3+\eps^2d_\eps).
       \label{eq:second-order-phase-error}
\end{align}
Consequently, in $C([0,1];\R^2)$,
\begin{equation}\label{eq:second-order-limit}
 \begin{pmatrix}
 \Lambda_\eps^2[C(T_\eps s)-c_0-\Lambda_\eps^{-1}g_0W_s^\eps]\\[2pt]
 \Lambda_\eps^4[Y_\eps(s)-\Lambda_\eps^{-3}g_0\int_0^sW_r^\eps\dd r]
 \end{pmatrix}
 \ \Longrightarrow\
 \begin{pmatrix}Z_s\\ \int_0^sZ_r\dd r\end{pmatrix},
 \quad Z_s=b_0s+\tfrac12g_0g'_0(W_s^2-s).
\end{equation}
Stopped continuous extensions are understood on the exit event. The same limit holds with $C$ replaced by $c$ or $c_{\rm v}$, and with $a$ replaced by $a_{\rm v}$.
\end{theorem}
\begin{proof}
First work with the extended closed SDE and write
$U_s=\bc(T_\eps s)-c_0$, $d=d_\eps$. Its speed equation becomes
\[
 \dd U_s=d\,g(c_0+U_s)\dd W_s^\eps
               +d^2b_{\rm sp}(c_0+U_s)\dd s,\qquad U_0=0.
\]
Bounded coefficients and BDG give $\norm{\sup|U|}_{L^4}\le Cd$. Subtracting $dg_0W^\eps$ and using Lipschitz continuity gives
$\norm{\sup|U-dg_0W^\eps|}_{L^4}\le Cd^2$.
Taylor expansion of the diffusion coefficient and another use of BDG now yield
\[
 \left\|\sup_{s\le1}|U_s-dg_0W_s^\eps-d^2Z_s^\eps|\right\|_{L^2}\le Cd^3.
\]
Indeed, the remaining diffusion integrand, before its outside factor $d$, is
\[
 g'_0(U_s-dg_0W_s^\eps)+O(U_s^2),
\]
which is $O_{L^2}(d^2)$ in the required time norm; the drift remainder is $d^2O_{L^2}(d)$.

The stochastic identity used here is
\[
 \int_0^sW_r^\eps\dd W_r^\eps=\frac12[(W_s^\eps)^2-s].
\]
It follows by applying the scalar \Ito\ formula to $W^2$, whose quadratic variation contributes $\dd s$. In particular this iterated stochastic integral has zero expectation. This is the second-order scalar \Ito--Taylor term \cite[Chapter~5]{KloedenPlaten1992}.

For the closed position, the exact identity is
\begin{align*}
 \eps^2[\ba(T_\eps s)-a_0-c_0T_\eps s]
 ={}&d^2\int_0^sU_r\dd r
    +\eps^2d\int_0^s\sigma_a(c_0+U_r)\dd\widetilde B_r\\
 &+\eps^2d^2\int_0^s\gamma_a(c_0+U_r)\dd r,
\end{align*}
where $\widetilde B=T_\eps^{-1/2}B_{T_\eps\cdot}$. Thus its remainder after the two displayed integrated-speed terms is at most $C(d^5+\eps^2d)$ in maximal $L^2$. Theorem~\ref{target:main} adds errors $C\eps d$ and $C\eps d^3$ on $A_\eps$, proving \eqref{eq:second-order-speed-error}--\eqref{eq:second-order-phase-error}.

After division by $d^2$ and $d^4$, the remainders tend to zero, since
\[
 d_\eps+\eps\Lambda_\eps+\eps^2\Lambda_\eps^3\longrightarrow0.
\]
The law of $(Z^\eps,J^\eps)$ is independent of $\eps$, and $\Prob(A_\eps^c)\to0$. This proves \eqref{eq:second-order-limit}. The localized coordinate corrections obey
\[
 \Lambda_\eps^2\sup(|C-c|+|c_{\rm v}-c|)
       =O_{L^2}(\eps^2\Lambda_\eps^3),\qquad
 \eps^2\Lambda_\eps^4\sup|a_{\rm v}-a|
       =O_{L^2}(\eps^3\Lambda_\eps^{9/2})
\]
before the stop. They vanish and give the final assertion.
\end{proof}

This second-order limit uses the Brownian motion selected by the comparison theorem. It is a coupled path statement. Parameter means, treated below, can be estimated without access to that Brownian path.

\begin{corollary}[The effect of deleting the spectral drift]\label{cor:drift-deletion}
Let $(a^{\rm z},C^{\rm z})$ solve the same closed SDE with the same Brownian motion and initial data, but replace $b_{\rm sp}$ by zero in the speed equation. Keep the phase coefficients unchanged. Then, uniformly for $0\le s\le1$, in probability,
\begin{equation}\label{eq:drift-deletion}
 \begin{aligned}
 \Lambda_\eps^2[C(T_\eps s)-C^{\rm z}(T_\eps s)]&\longrightarrow b_0s,\\
 \eps^2\Lambda_\eps^4[a(T_\eps s)-a^{\rm z}(T_\eps s)]
       &\longrightarrow\tfrac12b_0s^2.
 \end{aligned}
\end{equation}
Use the same stopped convention for the PDE coordinates as in Theorem~\ref{thm:second-order}.
\end{corollary}
\begin{proof}
Repeat the closed-SDE expansion with $b_0=0$. Its leading Brownian term and centered-square correction agree with those in \eqref{eq:second-order-Z}; its second-order mean term is absent. Subtract the expansions and use \eqref{eq:second-order-speed-error}--\eqref{eq:second-order-phase-error}. The exit event has vanishing probability.
\end{proof}

Because $b_0>0$ by Proposition~\ref{prop:pure-drift-comparison}, these limits are nonzero for $s>0$. They express drift resolution more precisely than a comparison of error scales alone. They do not identify the direct phase drift $\gamma_a$ separately at its smaller scale.

\subsection{Stopped and survival-conditioned means}\label{sec:means}
A convergence-in-probability statement does not by itself permit taking expectations: a small exceptional event can carry a large observable. We therefore specify the observables and estimate their exceptional contributions. Let $\tau_{\rm P}$ be the PDE-only tube exit, omitting the effective speed from \eqref{eq:stopping-regime}. The proofs in Sections~\ref{sec:weighted-energy}--\ref{sec:orbital-exit} apply to this stop as well: none of their radiation or orbital estimates uses the effective speed, and the exit proof simply omits its exit event. In particular, the same bound \eqref{eq:main-exit} holds for $\tau_{\rm P}$.

Define, for deterministic $t$,
\begin{equation}\label{eq:stopped-observables}
 C_t^{\rm P}=C(t\wedge\tau_{\rm P}),\qquad
 X_t^{\rm P}=\eps^2[a(t\wedge\tau_{\rm P})-a_0-c_0(t\wedge\tau_{\rm P})],
 \qquad p_T=\Prob(\tau_{\rm P}\le T).
\end{equation}
The deterministic trajectory in $X^{\rm P}$ is stopped at the same time as the phase. This removes an artificial discrepancy caused only by freezing one of the two trajectories.

\begin{proposition}[Mean spectral acceleration and accumulated position bias]\label{prop:means}
For $1\le T\le S_0\eps^{-2}$, $S=\eps^2T$, after taking $S_0\le1$,
\begin{align}
 \left|\E(C_T^{\rm P}-c_0)-b_0S\right|
       &\le C(S^2+\eps S+Sp_T),\label{eq:stopped-speed-mean}\\
 \left|\E X_T^{\rm P}-\tfrac12b_0S^2\right|
       &\le C(S^3+\eps S^2+\eps^2S+Sp_T).
       \label{eq:stopped-phase-mean}
\end{align}
At $T=T_\eps$, these imply
\begin{equation}\label{eq:stopped-means-asymptotic}
 \E(C_{T_\eps}^{\rm P}-c_0)=b_0\Lambda_\eps^{-2}+O(\Lambda_\eps^{-4}),
 \qquad
 \E X_{T_\eps}^{\rm P}=\tfrac12b_0\Lambda_\eps^{-4}+O(\Lambda_\eps^{-6}).
\end{equation}
The same asymptotic formulas hold after conditioning on $\tau_{\rm P}>T_\eps$, where the stopped observables equal the original coordinates. They also hold for the symplectic geometric speed, the virial geometric speed, and the virial phase, with the corresponding substitutions.
\end{proposition}
\begin{proof}
Taking expectations in the dissipative energy inequality, before $\tau_{\rm P}$, gives
\begin{equation}\label{eq:first-moment-radiation}
 \E\int_0^{t\wedge\tau_{\rm P}}r_\beta(v(s))^2\dd s\le C\eps^2t
\end{equation}
for all $t\ge0$. This first-moment estimate also holds for $t<1$; it follows directly by dropping the nonnegative terminal energy. The martingale can be localized and then the localization removed, or its bounded bracket on the tube can be used directly.

The exact spectral equation and its localized coefficients give
\[
 \E(C_t^{\rm P}-c_0)=\eps^2\int_0^t
 \E\left[\mathbf1_{\{s<\tau_{\rm P}\}}
              (b_{\rm sp}(C_s)+\rho_C(s))\right]\dd s.
\]
The stopped martingale has zero expectation because its noise row is bounded. Cauchy--Schwarz and \eqref{eq:first-moment-radiation} bound the integrated $\rho_C$ contribution by $C\eps^3t$. Also bounded drift and diffusion give
\[
 |\E(C_s^{\rm P}-c_0)|\le C\eps^2s,
 \qquad \E|C_s^{\rm P}-c_0|^2\le C\eps^2s
 \quad (\eps^2s\le S_0).
\]
Taylor-expand $b_{\rm sp}$ at $c_0$. The linear term is estimated using
\begin{equation}\label{eq:alive-mean}
 \E[\mathbf1_{\{s<\tau_{\rm P}\}}(C_s-c_0)]
 =\E(C_s^{\rm P}-c_0)
   -\E[\mathbf1_{\{\tau_{\rm P}\le s\}}(C_{\tau_{\rm P}}-c_0)].
\end{equation}
The second term is bounded by $Cp_s$, since the stopped spectral speed stays in a fixed compact interval. The quadratic Taylor term is bounded by $C\eps^2s$. Consequently
\[
 \left|\E[\mathbf1_{\{s<\tau_{\rm P}\}}b_{\rm sp}(C_s)]-b_0\right|
       \le C(\eps^2s+p_s).
\]
Integration proves \eqref{eq:stopped-speed-mean}, in fact for all $0\le t\le T$ with $S$ replaced by $\eps^2t$.

For phase, the exact equations and the zero expectation of the stopped phase martingale give
\[
 \E X_T^{\rm P}=
 \eps^2\int_0^T\E[\mathbf1_{\{s<\tau_{\rm P}\}}(C_s-c_0)]\dd s
 +O(\eps^2S).
\]
The error contains $\eps^2\E\int\rho_a$, bounded by \eqref{eq:first-moment-radiation}, and the bounded phase \Ito\ drift of size $\eps^4T$. Use \eqref{eq:alive-mean} and the speed expansion under the time integral. Integrating $b_0\eps^2s$ gives $b_0S^2/2$. The integrated remainders are bounded by
$C(S^3+\eps S^2+Sp_T)$, proving \eqref{eq:stopped-phase-mean}.

To pass to surviving trajectories, the stopped equations and the radiation moments give
\begin{equation}\label{eq:mean-exception-moments}
 \norm{C_T^{\rm P}-c_0}_{L^2}\le C\sqrt S,
 \qquad \norm{X_T^{\rm P}}_{L^2}\le CS^{3/2},\qquad T\ge1.
\end{equation}
For the second bound, the integrated speed contributes $CS^{3/2}$, the phase martingale contributes $C\eps^2\sqrt S$, and the remaining drifts contribute $C\eps^2S$. These are at most $CS^{3/2}$ because $\eps^2\le S\le S_0\le1$.
If $p_T\le1/2$, Cauchy--Schwarz therefore bounds the change from the expectations in \eqref{eq:stopped-speed-mean}--\eqref{eq:stopped-phase-mean} to the conditional expectations by
\[
 C\sqrt S\,p_T^{1/2}\quad\hbox{and}\quad CS^{3/2}p_T^{1/2},
\]
respectively. At $T=T_\eps$, the exit estimate \eqref{eq:proved-exit} gives $p_T=O_N(\eps^N)$ for every fixed $N>0$, since $1/S=\Lambda_\eps^2$. Thus the displayed corrections are $O(S^2)$ and $O(S^3)$. The terms containing $\eps$ in \eqref{eq:stopped-speed-mean}--\eqref{eq:stopped-phase-mean} are also smaller than these respective orders. This proves the asserted unconditional stopped and conditional asymptotics.

Finally, the maximal coordinate differences are $O_{L^2}(\eps^2\Lambda_\eps)$ in speed and $O_{L^2}(\eps^3\sqrt{\Lambda_\eps})$ in diffusively scaled phase. Their expectations, with or without conditioning on the survival event, are negligible at the stated remainder scales.
\end{proof}

The bias in \eqref{eq:stopped-means-asymptotic} is consistent with $\E Z_s=b_0s$: the centered-square term has zero mean. The direct phase drift cannot be extracted from this calculation, since the bound $O(\eps^2S)$ in \eqref{eq:stopped-phase-mean} also contains the quadratic radiation contribution. These moment statements specify what can be compared with an ensemble after recording its exit rule. They require no definition of an unrestricted spectral coordinate after the eigenvalue or chart leaves the controlled neighborhood.

\subsection{Parameter statistics from an independently driven SDE}\label{sec:observable-laws}
The strong coupling supplies a bound between laws, although constructing that coupling is unnecessary when sampling the reduced model independently. This is useful for distribution plots, time-dependent statistics, and the numerical comparisons in Section~\ref{sec:numerics}.

For $0\le s\le1$, define the canonically stopped PDE path
\[
 \mathcal Z_\eps(s)=
 \left(\Lambda_\eps[C_{T_\eps s}^{\rm P}-c_0],\,
       \Lambda_\eps^3X_{T_\eps s}^{\rm P}\right),
\]
and let $\overline{\mathcal Z}_\eps$ denote the corresponding unstopped normalized path of the extended effective SDE. Give $C([0,1];\R^2)$ its supremum norm, using the Euclidean norm in $\R^2$.

\begin{corollary}[A bound for path statistics]\label{cor:observable-laws}
For every functional $F:C([0,1];\R^2)\to\R$ with $|F|\le1$ and Lipschitz constant at most one,
\begin{equation}\label{eq:observable-law-bound}
 \left|\E F(\mathcal Z_\eps)-\E F(\overline{\mathcal Z}_\eps)\right|
 \le C\eps+2\Prob(\tau\le T_\eps)
 \le C'\eps .
\end{equation}
The effective-SDE expectation may be computed using Brownian samples independent of the PDE ensemble. The same conclusion, with an adjusted constant, holds for the virial coordinates.
\end{corollary}
\begin{proof}
Under the coupling of Theorem~\ref{target:main},
\[
 \left\|\mathbf 1_{A_\eps}
    \sup_{s\le1}|\mathcal Z_\eps(s)-\overline{\mathcal Z}_\eps(s)|\right\|_{L^2}
 \le C\eps,\qquad A_\eps=\{\tau>T_\eps\}.
\]
This follows from \eqref{eq:target-errors} after multiplication by $\Lambda_\eps$ and $\Lambda_\eps^3$ in the respective coordinates. Lipschitz continuity bounds the contribution on $A_\eps$ by $C\eps$, and boundedness gives $2\Prob(A_\eps^c)$ on the complementary event. By \eqref{eq:proved-exit} this probability is bounded by $Ce^{-c\Lambda_\eps^2}+C(1+T_\eps)e^{-c_*/\eps^2}=o(\eps)$. The expectation of the effective path functional depends only on its marginal law. For the virial coordinates, \eqref{eq:chart-change-moments} and \eqref{eq:spectral-geometric} add at most $C(\eps^2\Lambda_\eps^2+\eps^3\Lambda_\eps^{7/2})=o(\eps)$ in the normalized path norm on $A_\eps$.
\end{proof}

For an unbounded observable, an exit-probability estimate alone leaves its contribution on the exit event uncontrolled. Proposition~\ref{prop:means} supplies the moment bounds for the mean drift.

\subsection{Drift resolution and position accuracy}
The drift $\eps^2b_{\rm sp}$ enters the second-order expansion. At $T=T_\eps$, the ratios of the proved model errors to the scales generated by this speed drift are
\begin{equation}\label{eq:drift-resolution}
 \frac{\eps/\Lambda_\eps}{\Lambda_\eps^{-2}}=\eps\Lambda_\eps\to0,
 \qquad
 \frac{\eps/\Lambda_\eps^3}{\Lambda_\eps^{-4}}=\eps\Lambda_\eps\to0.
\end{equation}
The second denominator is the scale in $\eps^2a$ obtained by integrating the speed drift. Corollary~\ref{cor:drift-deletion} identifies the corresponding model discrepancy under the specified coupling. Proposition~\ref{prop:means} supplies the additional moment estimates needed for stopped and survival-conditioned expectation statements.

Corollary~\ref{cor:physical-reconstruction} gives absolute phase tracking and $H^1$ reconstruction with the coupled model for $1\le T=o(\eps^{-4/3})$. At $T_\eps$ the available bound for the unscaled phase error grows like $\eps^{-1}\Lambda_\eps^{-3}$. It is a factor $\eps$ smaller than the leading fluctuation scale in \eqref{eq:absolute-phase-scale}. Absolute position tracking at this later endpoint remains open.

The table gives the endpoint scales. The mean column refers to the stopped or survival-conditioned quantities in Proposition~\ref{prop:means}. The last column gives upper bounds from the strong comparison under the adapted coupling.
\begin{center}
\begin{tabular}{@{}llll@{}}
\toprule
Observable & Fluctuation & Mean bias & Model error bound \\
\midrule
$C-c_0$ & $\Lambda_\eps^{-1}$ & $\Lambda_\eps^{-2}$ & $\eps\Lambda_\eps^{-1}$\\[3pt]
$\eps^2(a-a_0-c_0T_\eps)$ & $\Lambda_\eps^{-3}$ & $\Lambda_\eps^{-4}$ & $\eps\Lambda_\eps^{-3}$\\[3pt]
$a-a_0-c_0T_\eps$ & $\eps^{-2}\Lambda_\eps^{-3}$ & $\eps^{-2}\Lambda_\eps^{-4}$ & $\eps^{-1}\Lambda_\eps^{-3}$\\
\bottomrule
\end{tabular}
\end{center}
The coefficients of the leading speed and position biases are $b_0$ and $b_0/2$, respectively. The fluctuation coefficients follow from \eqref{eq:joint-limit}.
\section{Numerical experiments}\label{sec:numerics}

Debussche and Printems \cite{DebusschePrintems1999} simulated KdV solitons under homogeneous noise. Here an independent reduced-model ensemble tests parameter laws, and a separate experiment projects the discretized PDE noise onto two adapted Brownian directions to compare individual paths. These experiments address different conclusions of Theorem~\ref{target:main}.

The supplement regenerates the coefficient figure and table, the reduced-SDE endpoint distributions and moment statistics, and the paired drift-deletion experiment. For the PDE, the archived numerical material consists of the reported tables and vector figures. In particular the four ordered endpoint marginals used below can be recovered from their plotted empirical distribution functions. The original noise seeds and full trajectories are unavailable for independently rerunning the other PDE diagnostics. The estimates from 48 PDE paths at each of two noise strengths serve as illustrations at finite resolution.

\subsection{Coefficients computed before any PDE simulation}\label{sec:numerical-coefficients}
We take $c_0=1$, $a_0=0$, and the Gaussian kernel
\[
 k(x)=\pi^{-1/4}e^{-x^2/2},\qquad
 R(z):=\int_\R k(x)k(x+z)\dd x=e^{-z^2/4}.
\]
Thus $\norm{k}_2=1$, and the covariance of the field increment at two points is $R(x-y)\dd t$. Write
\[
 \mathcal V(f,g)=\iint_{\R^2}f(x)R(x-y)g(y)\dd x\dd y.
\]
The entries of $q_0$ are obtained by applying this bilinear form to $Q_cP_c/D_c$ and $Q_c^2/D_c$. To evaluate the phase trace, Appendix~\ref{app:hessian} gives the entirely deterministic formula
\begin{equation}\label{eq:numerical-gamma}
\begin{split}
2D_c\gamma_a(c)={}&-\frac2{D_c}\mathcal V(Q_cP_c,Q_cP_c')
 +q_{aa}\ip{Q_c''}{P_c}-q_{CC}\ip{\partial_c^2Q_c}{P_c}\\
&+\frac2{D_c}\mathcal V(Q_c^2,Q_c\partial_cP_c)
 +2q_{aC}\ip{Q_c'}{\partial_cP_c}
 -2q_{CC}\ip{\partial_cQ_c}{\partial_cP_c}.
\end{split}
\end{equation}
All pairings are absolutely convergent. In particular, no stochastic samples or PDE data are used to estimate this trace.

The spectral factorization also yields a convenient kernel formula for $b_{\rm sp}$. Put $\kappa=\sqrt c/2$ and $L_0=-\partial_x^2+\kappa^2$. With $\mathcal D_c$ as in Section~\ref{sec:coefficients},
\[
 R_c^{\rm sp}=\mathcal D_c^*L_0^{-2}\mathcal D_c,
 \qquad \mathcal D_c(hp_c)=p_ch',
 \qquad L_0^{-2}(x,y)=G_c(x-y),
\]
where
\[
 G_c(z)=\frac{(1+\kappa|z|)e^{-\kappa|z|}}{4\kappa^3}.
\]
The first identity follows from $\mathcal D_c\mathcal D_c^*=L_0$ and is zero on $p_c$; on $p_c^\perp$ it is the inverse of $\mathcal D_c^*\mathcal D_c$. Substituting $h=Q_c\phi e_j$ and summing the derivative covariances gives
\begin{equation}\label{eq:numerical-bsp}
\begin{split}
b_{\rm sp}(c)=\frac19\iint G_c(x-y)p_c(x)p_c(y)
\big[&Q_c'(x)Q_c'(y)R(x-y)\\
 &+\{Q_c(x)Q_c'(y)-Q_c'(x)Q_c(y)\}R'(x-y)\\
 &-Q_c(x)Q_c(y)R''(x-y)\big]\dd x\dd y.
\end{split}
\end{equation}
The original coefficient calculation compared this formula with an independently assembled discrete Schr\"odinger reduced resolvent.

\begin{figure}[htbp]
\centering\includegraphics[width=\textwidth]{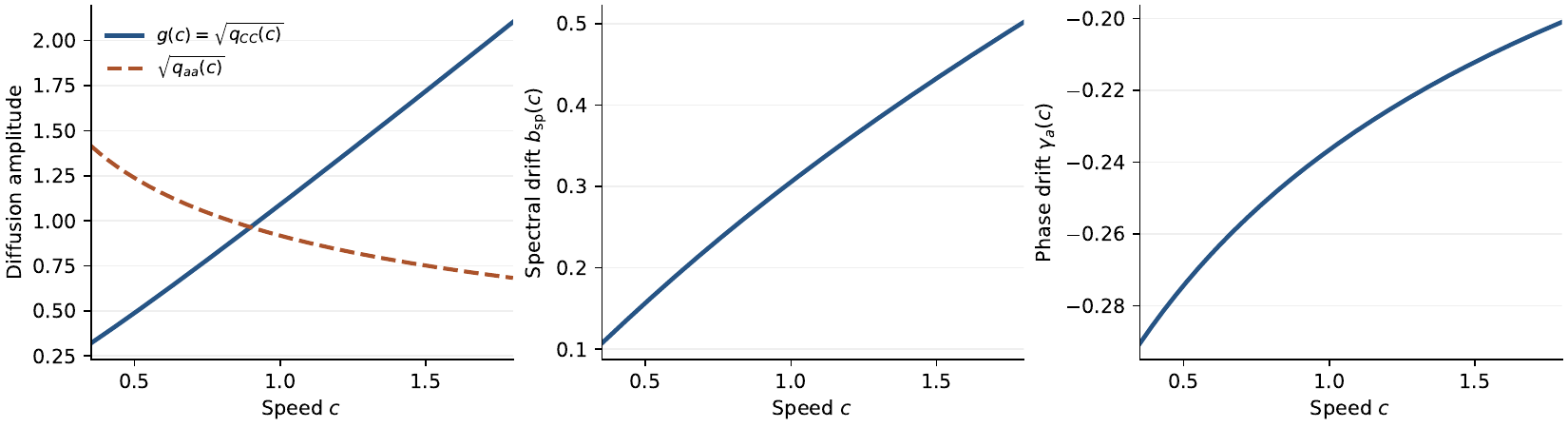}
\caption{Deterministic coefficient functions computed before simulating either equation. The positive spectral drift and the phase trace belong to the specified spectral speed and one-sided phase; their signs should not be transferred to a different chart.}\label{fig:coefficients}
\end{figure}

The production coefficient table has 146 speed nodes in $[0.35,1.8]$, quadrature spacing $1/8$, and spatial half-width 32. Cubic interpolation is used inside the simulation interval. Table~\ref{tab:numerical-coefficients} reports the more finely computed values at $c=1$; Figure~\ref{fig:coefficients} shows their dependence on speed. Halving the quadrature spacing changes $b_{\rm sp}(1)$ by about $1.1\times10^{-7}$. The reduced-resolvent calculation with 192, 256, and 384 Fourier grid points on a length-48 interval gives $0.3056573249$ in each case to the displayed precision. The $1.4\times10^{-8}$ difference from the finest convolution quadrature is much smaller than the stochastic effects examined below. In the original three off-grid checks, the largest interpolation discrepancy is $4.5\times10^{-9}$ and the largest change under quadrature refinement is $2.5\times10^{-7}$.

\begin{table}[htbp]
\centering
\caption{Pure-soliton coefficients at $c=1$, Gaussian correlation length one. The displayed values use spacing $1/16$ and spatial half-width 32.}\label{tab:numerical-coefficients}
\begingroup\small\setlength{\tabcolsep}{5pt}
\begin{tabular}{@{}ccccc@{}}
\toprule
$q_{aa}$ & $q_{aC}$ & $q_{CC}$ & $b_{\rm sp}$ & $\gamma_a$ \\
\midrule
$0.841880725$ & $0.332328762$ & $1.187628582$ & $0.305657311$ & $-0.236580841$ \\
\bottomrule
\end{tabular}
\endgroup

\end{table}

The original directional check of the phase Hessian solves the implicit chart for $Q_c\pm\eta h$, where $h$ is a smooth localized direction, and compares its centered second difference with Appendix~\ref{app:hessian}. For $\eta=0.004,0.002,0.001$, the absolute errors are approximately $8.95\times10^{-6}$, $2.24\times10^{-6}$, and $5.60\times10^{-7}$. The chart residual is tightened to $10^{-13}$ for this check, since division by $\eta^2$ would otherwise amplify the root-solving tolerance. This tests the phase convention as well as the signs in \eqref{eq:numerical-gamma}.

\subsection{SDE and PDE discretizations}\label{sec:numerical-methods}
For the reduced model, use $s=t/T_\eps$, $d=\Lambda_\eps^{-1}$, and $Y(s)=\eps^2[\ba(T_\eps s)-T_\eps s]$. The equations integrated on $0\le s\le1$ are
\begin{equation}\label{eq:numerical-scaled-sde}
\begin{aligned}
 \dd\bc&=d\,g(\bc)\dd B_1+d^2b_{\rm sp}(\bc)\dd s,\\
 \dd Y&=d^2(\bc-1)\dd s+\eps^2d\,\sigma_a(\bc)\dd B
                    +\eps^2d^2\gamma_a(\bc)\dd s.
\end{aligned}
\end{equation}
This formulation avoids subtracting large ballistic positions and makes the computational cost of an SDE path essentially independent of the physical value of $T_\eps$. The scalar speed update is Milstein's method: in addition to its Euler--Maruyama terms it includes
\[
 \tfrac12d^2g(\bc_n)g'(\bc_n)[(\Delta B_{1,n})^2-\Delta s].
\]
The phase uses a trapezoidal speed integral and left-endpoint \Ito\ noise. No higher strong-order assertion is made for this combined update. For background on these methods and the distinction between strong and weak numerical errors, see \cite{Higham2001,KloedenPlaten1992}. The ensembles in Tables~\ref{tab:independent-laws} and~\ref{tab:sde-means} use 1024 steps and 8192 independent antithetic pairs per noise strength. We generate one ensemble at each $\eps$ and reuse the same endpoint draws whenever that $\eps$ appears in either table or its corresponding figure.

For the PDE, set $w(t,y)=u(t,y+t)$. The Fourier linear multiplier of its deterministic equation is $i(k^3+k)$. We use a dealiased Fourier discretization and the fourth-order exponential time-differencing Runge--Kutta method of Cox and Matthews \cite{CoxMatthews2002}. Following Kassam and Trefethen \cite{KassamTrefethen2005}, we evaluate its coefficients by contour quadrature, using the full circle for the complex Fourier multipliers. This deterministic step is composed by Lie splitting with the exact multiplicative-noise substep
\begin{equation}\label{eq:numerical-noise-step}
 w\longmapsto w\exp\left(\eps\Delta W_{k,h}
                 -\tfrac12\eps^2R_h(0)\Delta t\right).
\end{equation}
Here $R_h(0)$ is the actual discrete field variance per unit time. Omitting the compensating term would change the \Ito\ drift. The Fourier covariance multiplier is $2\sqrt\pi e^{-k^2}$, and spatial white Gaussian increments carry the factor $(\Delta x)^{-1/2}$. This normalization is retained in every covariance projection.

The archived PDE runs used a domain of length 384, mesh width $1/4$, and $\Delta t\le0.01$, with the final step count adjusted to end exactly at $T_\eps$. A smooth sponge begins at $|y|=140$ and reaches strength one at the endpoints. Its purpose is to suppress periodic return of outgoing radiation; it changes the finite-domain evolution. The phase test equals $P_c$ up to relative position 40 and is smoothly set to zero by 48. A sixth-order finite-difference Schr\"odinger operator on a window of half-width 24 supplies the eigenvalue and its eigenfunction. Inverse iteration is continued until the eigen residual is below $2\times10^{-10}$. The chart is solved by Newton iteration with residual below $2\times10^{-10}$. The archived refinement results, including the influence of the sponge and cutoff, are reported below.

The original numerical study reported two deterministic checks. A unit-speed soliton is stationary in this frame and is preserved to near roundoff on $0\le t\le163.36$. Since that check is unusually favorable, a second check propagates a soliton of speed $1.15$ over the same interval. At $\Delta t\le0.02$, its phase error is about $9.0\times10^{-6}$ and its $L^2$ profile error is about $2.4\times10^{-5}$. For an exact unit-speed soliton, the mesh-$1/4$ eigenvalue extraction has speed bias $5.6\times10^{-8}$; at mesh $1/8$ that bias is $6.3\times10^{-10}$ on the half-width-24 window. An independent banded Hermitian eigensolve at all 96 saved PDE endpoints agrees with inverse iteration in spectral speed to within $1.8\times10^{-13}$; the smallest gap between the first two discrete eigenvalues is $0.132$. These are tests of the numerical implementation, not error bounds for the stochastic ensemble.

\subsection{Diagnostic stops and statistical uncertainty}\label{sec:numerical-stops}
Every displayed ensemble records its stopping rule. The reduced model stops at the first detected exit of its speed from $(0.5,1.5)$. The PDE diagnostics require $c,C\in(0.5,1.5)$, chart and eigen residuals below $10^{-7}$, and
\[
 \norm{v}_{H^1(\text{domain})}
   +\norm{e^{0.02y}v}_{H^1(\text{domain})}<2,
\]
with $y$ measured from the current phase. The soliton and phase-test window must also remain inside the undamped region. In the coupled experiment the reduced speed is included in the common stop. Diagnostic checks occur at every step in coupled runs and at intervals at most $0.2$ in independent PDE runs. The state at the first detected exit is retained, including its small discrete overshoot; both its state and clock are then frozen.

The numerical radius 2, speed interval, and weight 0.02 are practical choices. They are not certified values of the constants in Theorem~\ref{target:main}. Moreover, a sponge removes far-field radiation and can change a global radiation diagnostic. Accordingly the reported exit fractions are numerical diagnostics, not estimates of the exact real-line exit probability in \eqref{eq:main-exit}. Different stopping criteria for the standalone SDE and PDE must also be kept in mind when comparing finite-noise marginal laws. The definitions agree with the use of a stopped ballistic clock: the compensated phase at exit is $a(t\wedge\tau)-t\wedge\tau$, not $a(t\wedge\tau)-t$.

Antithetic sampling means that, for each Brownian path $B$, the SDE is also driven by $-B$. Each has the correct marginal law, but the pair is dependent. For standard errors we therefore average within each pair and regard the 8192 pair averages as independent observations. Error bars for means are $1.96$ times the resulting Monte Carlo standard error; they exclude discretization and reduction error. The original PDE ensemble uses independent cylindrical noises. The SDE samples displayed in the distribution and scaling figures are generated afresh from the seeds listed with the code.

For the empirical PDE distribution functions, we display the band
\[
 \widehat F_N(x)\ \pm\sqrt{\frac{\log(40)}{2N}},
\]
clipped to $[0,1]$. The Dvoretzky--Kiefer--Wolfowitz inequality with Massart's constant \cite{Massart1990} gives 95\% simultaneous coverage in $x$ for one marginal distribution of the implemented, stopped numerical model. It is not a simultaneous band for all panels, nor a confidence band around the exact SPDE law. At $N=48$ its half-width is about $0.196$, which makes the limited distributional resolution visible.

\subsection{Independent laws and perturbation memory}\label{sec:numerical-laws}
\begin{figure}[htbp]\centering
\includegraphics[width=.98\textwidth]{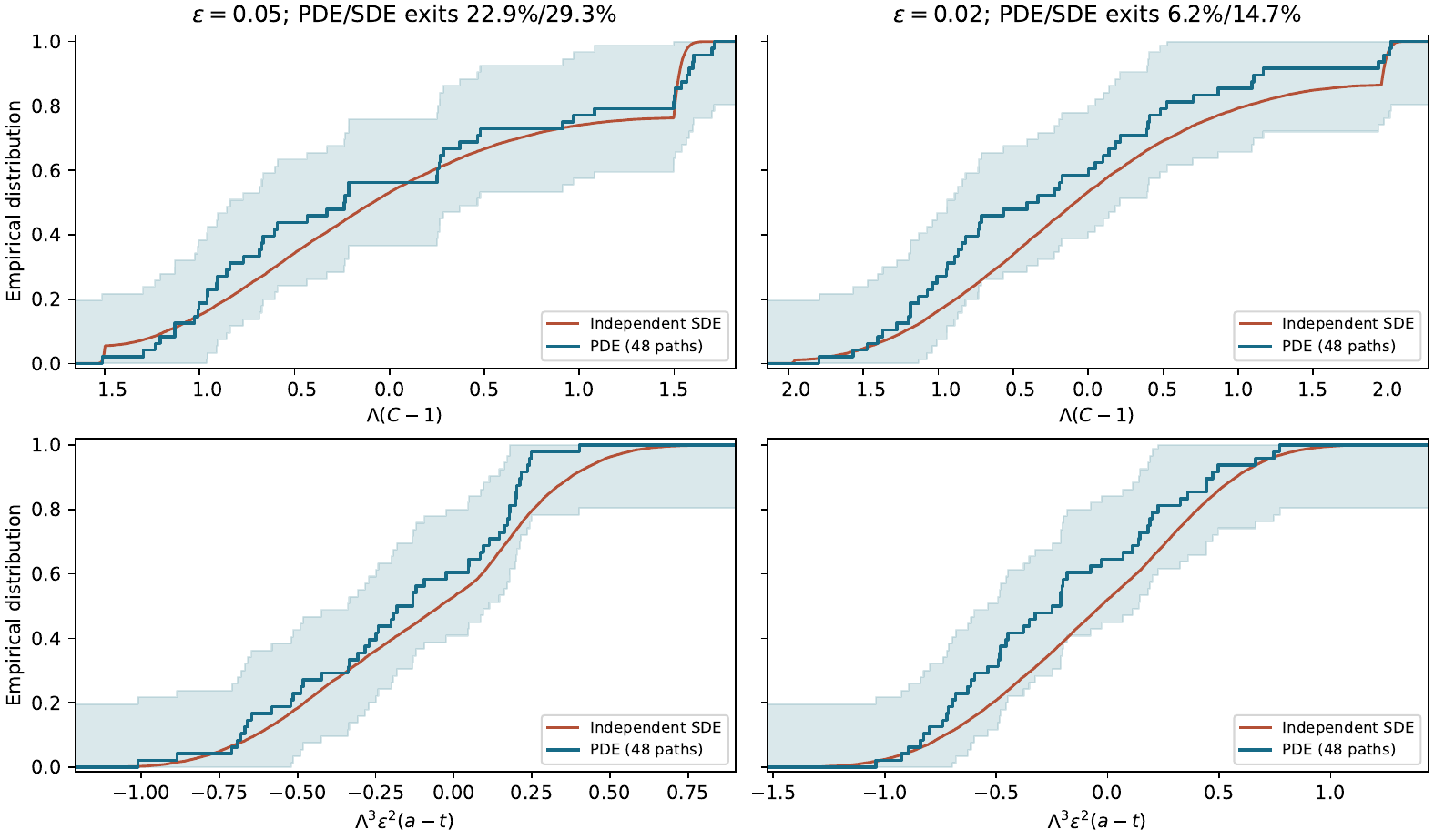}
\caption{Stopped parameter distributions. Each blue staircase uses the 48 ordered PDE endpoints recovered from one marginal curve in the archived vector figure. These marginal samples have no trajectory identifiers. The red curves use the same independent SDE samples as Table~\ref{tab:independent-laws}. Shading gives the one-marginal PDE sampling band.}\label{fig:independent-laws}
\end{figure}

The original PDE experiment used 48 trajectories at each of $\eps=0.05$ and $0.02$, over the respective intervals $T_\eps\approx44.57$ and $163.36$, subject to the diagnostic stops. The independent SDE ensemble has 16384 paths at each noise strength. Figure~\ref{fig:independent-laws} compares the stopped endpoint observables
\[
 \Lambda_\eps(C-1),\qquad
 \Lambda_\eps^3\eps^2(a-t).
\]
The coefficient functions are never fitted to these PDE trajectories. Table~\ref{tab:independent-laws} reports means and Monte Carlo standard errors in these same normalized variables, together with the diagnostic exit fraction. These small PDE ensembles provide a distributional illustration. Their mean uncertainty is too large to establish the predicted small drift, and differences at the scale of the displayed sampling bands cannot be used as precise estimates of the reduction error.

The empirical speed/phase CDF discrepancies are $0.133/0.185$ at $\eps=0.05$ and $0.199/0.203$ at $\eps=0.02$, comparable with the one-marginal PDE sampling-band half-width $0.196$. In particular, the curves do not all lie within that band. These observations must be read with the finite SDE sampling error and the differing diagnostic stops; we do not infer a precise goodness-of-fit conclusion from these panels.

\begin{table}[htbp]\centering
\caption{Independent stopped ensembles. PDE entries are the original reported means and standard errors. SDE entries use the samples underlying Figures~\ref{fig:independent-laws} and~\ref{fig:sde-scaling}. Entries after $\pm$ give one Monte Carlo standard error.}\label{tab:independent-laws}
\begingroup\small\setlength{\tabcolsep}{5pt}
\begin{tabular}{@{}llrrrr@{}}
\toprule
$\eps$ & Model & $N$ & Exit & $\E[\Lambda(C-1)]$ & $\E[\Lambda^3Y]$ \\
\midrule
0.05 & PDE & 48 & 22.9\% & $-0.005\pm0.148$ & $-0.181\pm0.052$ \\
0.05 & SDE & 16384 & 29.3\% & $0.088\pm0.003$ & $-0.094\pm0.001$ \\
0.02 & PDE & 48 & 6.2\% & $-0.215\pm0.147$ & $-0.209\pm0.071$ \\
0.02 & SDE & 16384 & 14.7\% & $0.075\pm0.002$ & $-0.053\pm0.001$ \\
\bottomrule
\end{tabular}
\endgroup

\end{table}

\begin{figure}[htbp]\centering
\includegraphics[width=\textwidth]{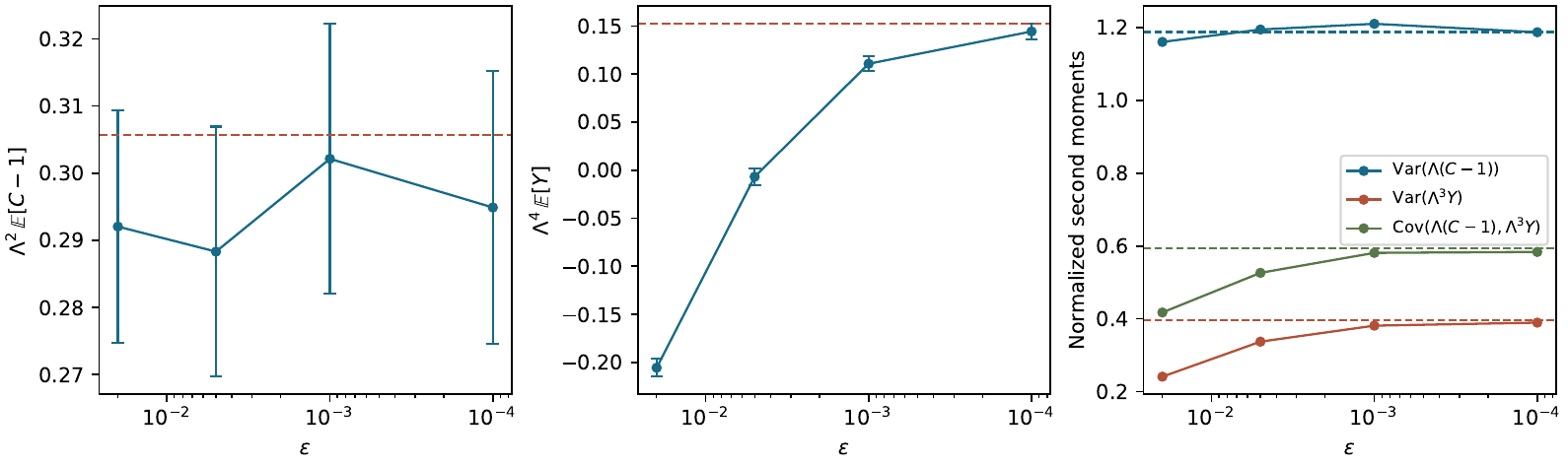}
\caption{Mean and variance scaling of the standalone SDE. Mean bars are pointwise 95\% Monte Carlo intervals. Dashed lines give the analytic asymptotic coefficients. The finite-noise effect of stopping is visible, particularly in the phase mean.}\label{fig:sde-scaling}
\end{figure}

The reduced model is also simulated at $\eps=0.005,0.001,10^{-4}$, where full PDE ensembles would be much more expensive. The endpoint $T_{10^{-4}}$ is about $1.18\times10^6$, although the rescaled SDE still runs on $[0,1]$. Figure~\ref{fig:sde-scaling} and Table~\ref{tab:sde-means} compare its normalized means and variances with the asymptotic coefficients. The limiting covariance is
\begin{equation}\label{eq:limit-cov}
 q_{CC}(1)\begin{pmatrix}1&1/2\\1/2&1/3\end{pmatrix},
 \qquad q_{CC}(1)\approx1.187628582.
\end{equation}
The speed and phase mean targets are $b_0$ and $b_0/2$ after multiplication by $\Lambda_\eps^2$ and $\Lambda_\eps^4$, respectively. These lines are asymptotic targets, not exact finite-noise identities. In particular, a stopped phase mean can be negative at moderate noise while its asymptotic coefficient is positive. Formula~\eqref{eq:alive-mean} explains why conditioning on continued motion affects the accumulated speed contribution. Exit fractions and the actual stopping rule are essential to interpreting that sign.

\begin{table}[htbp]\centering
\caption{Stopped reduced-SDE statistics. Mean entries are normalized as indicated; $\pm$ denotes one standard error computed from antithetic pair averages.}\label{tab:sde-means}
\begingroup\small\setlength{\tabcolsep}{5pt}
\begin{tabular}{@{}lrrr@{}}
\toprule
$\eps$ & Exit & $\Lambda^2\E[C-1]$ & $\Lambda^4\E[Y]$ \\
\midrule
0.02 & 14.70\% & $0.292\pm0.009$ & $-0.206\pm0.005$ \\
0.005 & 4.31\% & $0.288\pm0.009$ & $-0.007\pm0.005$ \\
0.001 & 0.96\% & $0.302\pm0.010$ & $0.111\pm0.004$ \\
0.0001 & 0.07\% & $0.295\pm0.010$ & $0.144\pm0.004$ \\
\bottomrule
\end{tabular}
\endgroup

\end{table}

In physical variables, the speed deviations vanish as $\eps\to0$, whereas typical absolute position fluctuations have scale $\eps^{-2}\Lambda_\eps^{-3}$ and grow. The rescaled phase plots show the displacement generated by accumulated speed fluctuations. Corollary~\ref{cor:physical-reconstruction} gives the shorter time window for approximation in physical space.

\subsection{Second order and the effect of deleting the spectral drift}\label{sec:numerical-second}
\begin{figure}[htbp]\centering
\includegraphics[width=\textwidth]{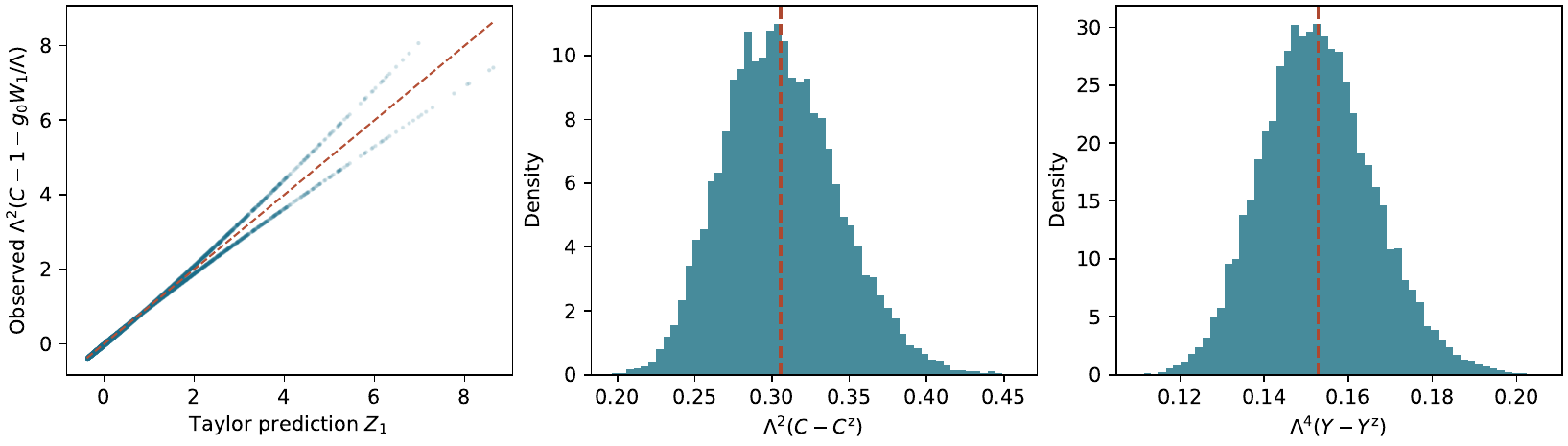}
\caption{Second-order speed residual and paired drift deletion at $\eps=10^{-4}$, using the same 16373 joint survivors in all panels. The diagonal gives the path-dependent stochastic Taylor prediction. The dashed vertical lines give $b_0$ and $b_0/2$ for the normalized speed and scaled-phase differences.}\label{fig:second-order}
\end{figure}

At $\eps=10^{-4}$ we generate a separate sample with 8192 antithetic pairs and 1024 time steps. Each path drives two models with identical Brownian increments. The second model deletes $b_{\rm sp}$ from the speed equation and evaluates all remaining coefficients at its own speed. We apply the diagnostic speed stop to each model. Of 16384 marginal paths, 16373 survive both stops through $s=1$; all three panels of Figure~\ref{fig:second-order} use these joint survivors.

Let $W_1$ be the sum of the first Brownian increments over the rescaled interval. For the full SDE the left panel compares the observed residual $\Lambda_\eps^2(\bc-1-\Lambda_\eps^{-1}g_0W_1)$ with
\[
 Z_1=b_0+\tfrac12g_0g'_0(W_1^2-1).
\]
The other panels compare the terminal differences of the full and drift-deleted equations, with the normalization of Corollary~\ref{cor:drift-deletion}. The empirical medians of the speed and phase differences are $0.3035$ and $0.1522$, respectively; the limiting values are $b_0\simeq0.3057$ and $b_0/2\simeq0.1528$. These medians describe the paired SDE experiment at this noise level.

\subsection{Adapted path comparison and coordinate transfer}\label{sec:numerical-coupling}
For the separate coupled experiment, the discrete implicit phase gradient is denoted by $G_a$ and the normalized spectral eigenfunction by $\psi$. At the left endpoint of a step set
\[
 f_a=uG_a,\qquad f_C=\tfrac23u\psi^2,
 \qquad q_{ij}=\mathcal V_h(f_i,f_j).
\]
The Brownian increments used in the reduced equation are
\begin{equation}\label{eq:numerical-B}
 \Delta B_1=\frac{\ip{f_C}{\Delta W_{k,h}}_h}{\sqrt{q_{CC}}},
 \qquad
 \Delta B_2=
 \frac{\ip{f_a}{\Delta W_{k,h}}_h-(q_{aC}/\sqrt{q_{CC}})\Delta B_1}
 {\sqrt{q_{aa}-q_{aC}^2/q_{CC}}}.
\end{equation}
With the covariance and quadrature used to generate $\Delta W_{k,h}$, these increments have conditional covariance $\Delta t\,\mathrm{Id}$. The construction uses a predictable numerical state and is the finite-dimensional version of \eqref{eq:adapted-B}. Merely giving the solvers the same random seed would not construct \eqref{eq:numerical-B}.

\begin{figure}[htbp]\centering
\includegraphics[width=.98\textwidth]{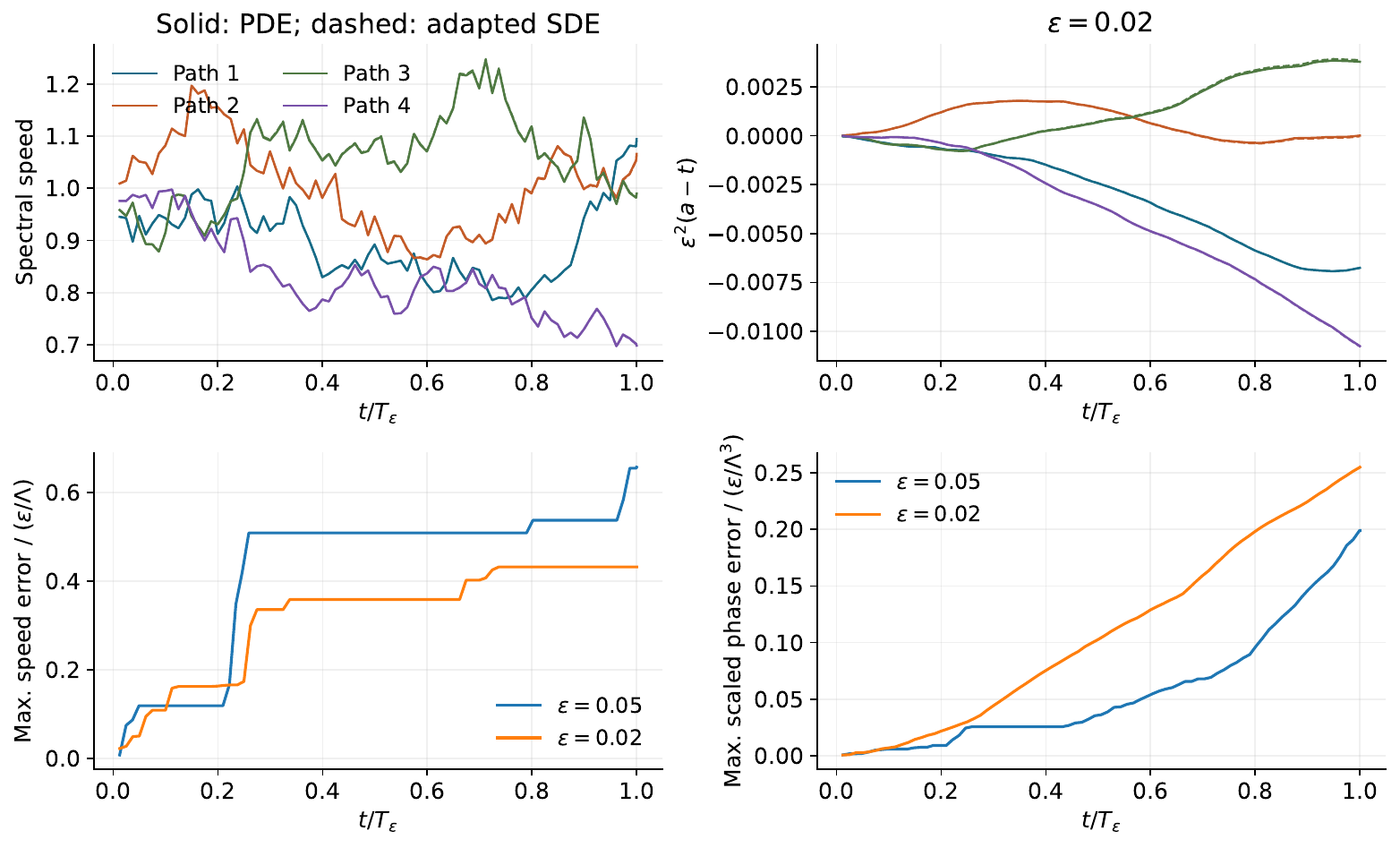}
\caption{Adapted PDE--SDE coupling. Top: four paths at $\eps=0.02$, solid for the PDE and dashed for its coupled SDE. Bottom: largest recorded errors relative to the analytic orders, for the two noise strengths. All stopped trajectories use their stopped clocks.}\label{fig:adapted-coupling}
\end{figure}

The archived coupled experiment used four paths at each of $\eps=0.05$ and $0.02$. Figure~\ref{fig:adapted-coupling} compares their stopped trajectories. Table~\ref{tab:coupled-errors} gives the largest recorded errors over these paths and approximately 80 output times, divided by the theoretical orders. These finite-sample maxima are not estimates of the maximal $L^2(\Omega)$ norm in Theorem~\ref{target:main}. Two noise strengths and four paths are insufficient for fitting an asymptotic rate, and the time-refinement differences below remain relevant to their interpretation.

\begin{table}[htbp]\centering
\caption{Largest recorded coupled errors. $E_C=\max|C-\bc|$ and $E_Y=\max\eps^2|a-\ba|$ use the sampled output times and four paths.}\label{tab:coupled-errors}
\begingroup\small\setlength{\tabcolsep}{5pt}
\begin{tabular}{@{}rrrrrr@{}}
\toprule
$\eps$ & $E_C$ & $E_Y$ & $E_C/(\eps/\Lambda)$ & $E_Y/(\eps/\Lambda^3)$ & Exit \\
\midrule
0.05 & $1.10\times10^{-2}$ & $3.70\times10^{-4}$ & 0.658 & 0.199 & 50\% \\
0.02 & $2.21\times10^{-3}$ & $8.52\times10^{-5}$ & 0.432 & 0.255 & 0\% \\
\bottomrule
\end{tabular}
\endgroup

\end{table}

\begin{figure}[htbp]\centering
\includegraphics[width=.9\textwidth]{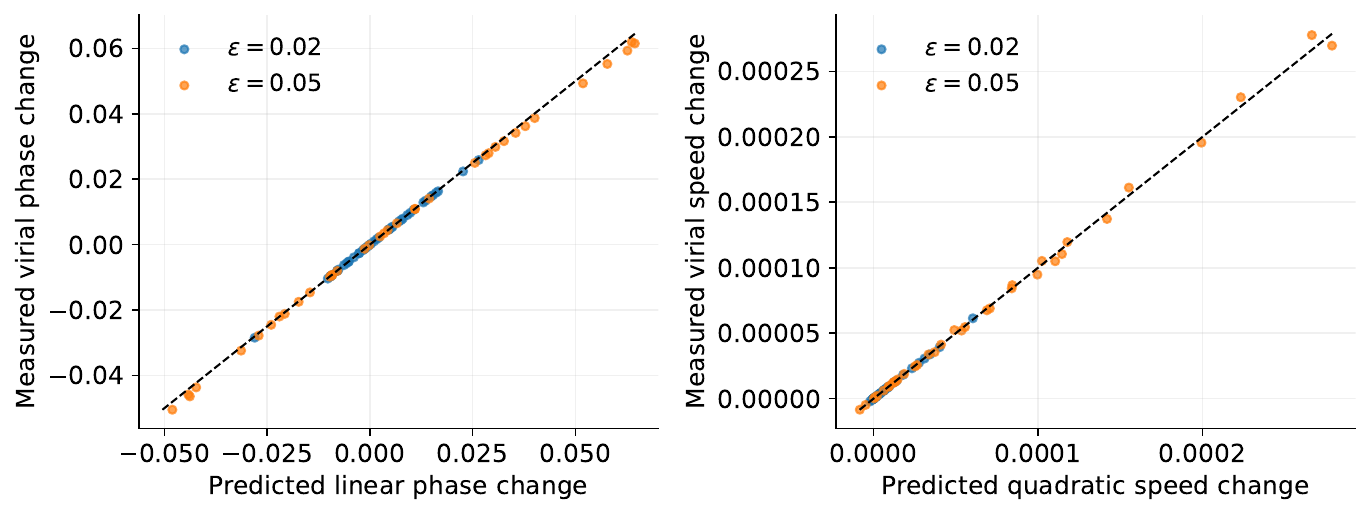}
\caption{Virial-coordinate changes measured on surviving PDE trajectories versus the linear phase and quadratic speed predictions of Proposition~\ref{prop:virial-chart}. Dashed lines are equality lines.}\label{fig:coordinate-transfer}
\end{figure}

Finally the virial chart is extracted from the same PDE states, without changing the primary phase or the simulation. Figure~\ref{fig:coordinate-transfer} compares its phase and speed changes with the linear and quadratic expressions in Proposition~\ref{prop:virial-chart}. Only trajectories surviving the PDE diagnostic stop are plotted. These checks test the coordinate mechanism at finite amplitude; they are not a substitute for the localized estimates.

\subsection{Refinement, reproducibility, and remaining numerical limits}\label{sec:remaining}
SDE step refinement uses coarse increments obtained by summing fine increments of the same Brownian motion. For $\eps=0.005$, 256 antithetic pairs and a 4096-step reference, the original experiment reported decreasing endpoint differences under refinement. In particular, at 1024 steps the speed RMS difference is about $4.7\times10^{-6}$ and the $\Lambda_\eps^3Y$ RMS difference is about $3.0\times10^{-4}$, restricted to joint survivors. These are differences to a numerical reference, not errors to an exact SDE solution.

For the PDE, Table~\ref{tab:pde-refinement} uses two fixed noises throughout the full $T_{0.02}$ interval. The reference has length 512, mesh $1/8$, time step at most $0.005$, sponge onset 196, and phase cutoff starting at 48. A common finest space--time noise grid is used for every case; spatial restrictions and explicit sums of time increments implement the coupling. Each row changes one setting. The errors shown are endpoint RMS differences over paths surviving both corresponding diagnostic stops.

\begin{table}[htbp]\centering
\caption{Paired PDE refinement at $\eps=0.02$, $T_\eps\approx163.36$. $N_{\rm joint}$ counts joint survivors out of two. Phase differences include the factor $\eps^2$.}\label{tab:pde-refinement}
\begingroup\small\setlength{\tabcolsep}{5pt}
\begin{tabular}{@{}lrrr@{}}
\toprule
Changed setting & $N_{\rm joint}$ & Speed RMS & Scaled phase RMS \\
\midrule
Time step $0.01$ & 2 & $1.24\times10^{-4}$ & $6.47\times10^{-6}$ \\
Time step $0.02$ & 2 & $5.09\times10^{-4}$ & $1.52\times10^{-5}$ \\
Mesh width $1/4$ & 2 & $1.32\times10^{-7}$ & $6.74\times10^{-14}$ \\
Length 384, sponge onset 140 & 2 & $1.37\times10^{-6}$ & $5.40\times10^{-8}$ \\
Phase cutoff starts at 32 & 2 & $0$ & $5.42\times10^{-10}$ \\
Sponge strength $1/2$ & 2 & $9.07\times10^{-6}$ & $2.37\times10^{-7}$ \\
\bottomrule
\end{tabular}
\endgroup

\end{table}

Time discretization is the largest measured numerical sensitivity in these tests. At the production time step, the difference to the finer reference is $1.24\times10^{-4}$ in speed and $6.47\times10^{-6}$ in scaled phase. The spatial-mesh difference is $1.32\times10^{-7}$ in speed. Shortening the domain and moving its sponge changes endpoint speed by $1.37\times10^{-6}$, while the phase-cutoff change contributes $5.42\times10^{-10}$ in scaled phase. These comparisons use different paths from the four-path model comparison, so their errors cannot simply be subtracted from that comparison.

The time-step sensitivity and the small number of PDE paths limit the accuracy with which one can compare the two models empirically. At each noise strength the four coupled paths illustrate strong comparison, while the independent ensembles illustrate agreement of marginal laws. The PDE data do not resolve the small mean speed drift or an asymptotic strong-error exponent.

The numerical supplement gives the commands and dependencies for recomputing the five pure-soliton coefficients and the five standalone SDE ensembles. Its plotting program uses the saved SDE endpoint arrays to generate Tables~\ref{tab:independent-laws} and~\ref{tab:sde-means} and Figures~\ref{fig:independent-laws} and~\ref{fig:sde-scaling}. It also rebuilds the coefficient curves in Figure~\ref{fig:coefficients} and the reference values in Table~\ref{tab:numerical-coefficients}. For the PDE columns of Table~\ref{tab:independent-laws} it retains the original reported numbers. Recovering the ordered marginals from the archived vector figure yields all four means and their standard errors to the precision printed in that table. The supplement also regenerates Figure~\ref{fig:second-order} from its paired SDE seed. The original PDE coupling and coordinate-transfer plots and the coupled and refinement tables are preserved as reported; the underlying PDE trajectories are unavailable.

\clearpage
\appendix
\section{An explicit formula for the pure phase Hessian}\label{app:hessian}
This appendix makes \eqref{eq:gamma} computable from profile data. Fix $c$ and a direction $h\in\Xb$, and use the phase and geometric-speed chart along $u_s=Q_c+sh$. Define
\begin{equation}\label{eq:hessian-alpha}
 \alpha=\frac{\ip{h}{P_c}}{D_c},\qquad
 \kappa=\frac{\ip{h}{Q_c}}{D_c},\qquad
 w=h+\alpha Q_c'-\kappa\dot Q_c.
\end{equation}
Then
\begin{equation}\label{eq:phase-hessian}
 \begin{split}
 D_cD^2\Ac(Q_c)[h,h]
 ={}&2\alpha\ip{h'}{P_c}+\alpha^2\ip{Q_c''}{P_c}
           -\kappa^2\ip{\partial_c^2Q_c}{P_c}\\
    &+2\kappa\ip{w}{\partial_cP_c}.
 \end{split}
\end{equation}
The pairing with $h'$ can equivalently be written $-\ip{h}{P_c'}$.

To prove this, write $a(s)=\Ac(u_s)$, $c(s)=\mathfrak c(u_s)$, and
$v_s=T_{a(s)}u_s-Q_{c(s)}$. At $s=0$, the first derivative of $v_s$ is $w$, and its second derivative is
\[
 2\alpha h'+\alpha^2Q_c''+a''(0)Q_c'
       -\kappa^2\partial_c^2Q_c-c''(0)\dot Q_c.
\]
Differentiate $\ip{v_s}{P_{c(s)}}=0$ twice. The term with $c''(0)$ vanishes because $\ip{\dot Q_c}{P_c}=0$, while the term with $a''(0)$ is $-D_ca''(0)$. The remaining terms give \eqref{eq:phase-hessian}. For general $h\in\Xb$, the identity follows by smooth approximation and continuity of all pairings; differentiating the translated tests is another direct justification.

Inserting $h=h_{c,j}=Q_c\phi e_j$ into
\eqref{eq:hessian-alpha}--\eqref{eq:phase-hessian} and taking one
half of the sum gives \eqref{eq:gamma}. Absolute convergence follows
from the Hilbert--Schmidt bounds and Cauchy--Schwarz in the noise
index. The terms involving $P_c$ and $\partial_cP_c$ specify the
dependence on the choice of phase.

\medskip
\bibliographystyle{amsalpha}
{\raggedright\bibliography{stochastic_kdv}}
\end{document}